\documentclass[12pt,a4paper]{article}
\usepackage{amssymb,amsfonts,amsmath,latexsym,amsthm,graphics}
\usepackage[T2A]{fontenc}
\usepackage[cp1251]{inputenc}

\usepackage{epsfig} 

\usepackage{graphicx} 

\numberwithin{equation}{section}

\title{{\bf Estimation under
uncertainties of acoustic and electromagnetic fields from noisy
observations
}}

\author{Yury Podlipenko,~$^1$ Yury Shestopalov,~$^2$ and Vladimir Prishlyak}

\date{$^1$ Kyiv National University, Kyiv, Ukraine\\
$^2$ Karlstad University, Karlstad, Sweden }

\begin{document}

\maketitle





\newcommand{\n}{^{(1)}}
\newcommand{\m}{^{(2)}}
\newcommand{\g}{^{(3)}}
\newcommand{\p}{^{(1,1)}}
\newcommand{\q}{^{(1,2)}}
\newcommand{\s}{^{(2,1)}}
\newcommand{\w}{^{(2,2)}}

\newtheorem{predlll}{Corollary}[section]
\newtheorem{pred}{Theorem}[section]
\newtheorem{predl}{Lemma}[section]
\newtheorem*{predllll}{Lemma}
\newtheorem{predll}{Definition}[section]
\newtheorem{predlllll}{Remark}

\renewcommand{\proofname}%
{\bf Proof}

\renewcommand{\contentsname}%
{\begin{center}Contents\end{center}}

\newpage

\addtocontents{toc}{\large}

\tableofcontents

\newpage
\begin{abstract}
\noindent The creation and justification of the methods for
minimax estimation of parameters of the external boundary value
problems  for the Helmholtz equation in unbounded domains are
considered. When observations are distributed in subdomains, the
determination of minimax estimates is reduced to the solution of
integro-differential equations in bounded domains. When
observations are distributed on a system of surfaces
the problem is reduced to solving integral equations on an
unclosed bounded surface which is a union of the boundary of the
domain and this system of surfaces. Minimax estimation of the
solutions to the boundary value  problems from point observations
is also studied.
\end{abstract}
\bigskip
MSC-class: 35J25, 45Fxx, 45Kxx, 49xx 93E10 (Primary), 78M50
(Secondary)

\section*{}

\addcontentsline{toc}{section}{Introduction}

\begin{center}
{\bf \Large {Introduction}}
\end{center}

\large

In the system analysis of complex processes described by partial
differential equations (PDEs), an important problem is the optimal
reconstruction (estimation) of parameters of the equations, like
values of some functionals on their solutions or right-hand sides,
from observations, which depend on the same solutions.

These problems play an important role in various areas of science
and engineering. Depending on the character of a priori
information, stochastic or deterministic approaches are possible.
The choice is determined by the nature of the problem parameters
which can be random or not. Moreover, the optimality of
estimations depends on a criterion with respect to which a given
value is evaluated.

The field of optimal control of PDEs
 has been strongly influenced by the work of J.L. Lions, who started the systematic study of optimal control
 problems for PDEs in \cite{BIBLlio}, in particular, singular perturbation problems in \cite{BIBLlio-1}
 and ill-posed problems in \cite{BIBLlio-2}. A possible direction of research
 in this field consists in extending results from the finite-dimensional case such as Pontryagin's
 principle, second-order conditions, structure of bang-bang controls, singular arcs and so on. On the other hand partial differential equations  have specific features such as finiteness of propagation for hyperbolic systems, or the smoothing effect of parabolic systems, so that they may present qualitative properties that are deeply different from the ones in the finite-dimensional case.
 The present study is devoted to a class of problems of optimal control and estimation for a specific family of
 PDEs of mathematical physics.

In practice, the data of boundary value problems (BVPs) for differential equations
that simulate a physical or technological object are always given with uncertainty.
For example, the right-hand sides of the equations, initial or boundary conditions may be known approximately;
that is, they belong to certain bounded sets in the corresponding functional spaces.

For solving the
estimation problems we must have supplementary data (observations)
$$
y=C\varphi+\eta,
$$
where $C$ is an operator that specifies the method of measuring and $\eta$ is the measurement error. As a rule, this error is not known and belongs to a certain given set and the operator is not invertible. Therefore, in general, from given $y,$ it is not possible to uniquely reconstruct the sought-for solution $\varphi$ of a BVP and, consequently,  quantity $l(\varphi)$, where $l$ is a given linear continuous functional.
We see that a natural problem arises: to determine a quantity $\widehat{l(\varphi)}$ which would provide the best (in a certain sense) approximation  to the sought-for $l(\varphi)$.

Let us briefly characterize the minimax approach to the solution
of this problem. We are looking for linear with respect to
observations optimal estimates of functionals of solutions and
right-hand sides of BVPs based upon the condition of minimum of
the maximal mean square error of estimation taken over the subsets
mentioned above.

These estimates were called minimax a priori or minimax program estimates (see \cite{BIBLkras}, \cite{BIBLkurzh}).

The situation when the unknown parameters of equations and
observations are perturbed by noise whose statistical
characteristics are not known completely constitutes the case of
special interest.

In the absence of true information about distribution of random
perturbations, the minimax approach proved to be a useful solution
technique. This approach initiated and developed by N.N.
Krasovskii \cite{BIBLkras}, A.B. Kurzhanskii \cite{BIBLkurzh},
O.G. Nakonechnyi \cite{BIBLnak22}, N.F. Kirichenko
\cite{BIBLbubl1}, and B.M. Pshenichnyi enabled one to find optimal
estimates of the BVP parameters for ordinary differential
equations corresponding to the worst realizations of random
perturbations.

The present work is devoted to the creation and rigorous justification of constructive
minimax  estimation methods of parameters of the external BVPs  for the Helmholtz equation in
arbitrary unbounded domains with finite boundaries. We reduce the determination of minimax
estimates to the solution of certain integro-differential equations in bounded domains
when observations are distributed in subdomains.
When observations are distributed on a system of surfaces (that simulate e.g. antennas)
the problem is reduced to solving some integral equations on an unclosed bounded surface which
is a union of the boundary of the domain and this system of surfaces.

These estimation problem are of tremendous significance in many areas of
applied electromagnetics,
acoustics, contact mechanics. Therefore, comprehensive theoretical analysis of estimation techniques is an urgent task.

{\bf Methods and objectives.} The study is aimed at elaboration of the methods of guaranteed estimation of the values of linear functionals defined on solutions to external BVPs for the Helmholtz equation and their right-hand sides.

This task can be fulfilled if the following problems are solved:
\begin{itemize}
\item To reduce estimation of the values of functionals defined on the solutions to external BVPs and the right-hand sides of equations that enter the problem statement to certain  problems of optimal control of systems that are described by certain conjugate BVPs for the Helmholtz equation in bounded domains with a quadratic quality criterion.

\item To obtain, for given restrictions on the unknown second moments of observation noise and unknown deterministic data of the BVPs under study, the systems of integro-differential and integral equations such that the minimax estimates of functionals are expressed in terms of their solutions

\item To prove unique solvability of the obtained systems of integro-differential and integral equations.

\end{itemize}

{\it The object of study} is observation problems under uncertainty when the functions that are observed on a system of subdomains or surfaces are coupled with the solutions to the considered BVPs via linear operators with additive measurement errors.

{\it The method of study.} The systems of integro-differential and integral equations obtained in this work whose solutions are used to express minimax estimates are based on the theory of generalized solutions to BVPs for the Helmholtz equation, utilization of the so-called Dirichlet-to-Neumann (DtN) data-transforming operators, and the the theory of potential in Sobolev spaces.

{\bf A remark on novelty.} For the first time we consider the statement of the problem of minimax estimation of
the parameters of external BVPs for the Helmholtz equation with general boundary
conditions that arise in the mathematical theory of wave diffraction.

For the systems described by such BVPs, we obtain representations for minimax estimates of the values of functionals from the observed solutions and right-hand sides that enter the problem statement; quadratic restrictions are imposed  on unknown deterministic data and second moment of observation noise. We also obtain representations for the estimation errors. The representations are obtained in terms of the solutions to certain uniquely solvable systems of integro-differential and integral equations in bounded domains.

When the unknown solutions of the system states are observed that are described by external BVPs for the Helmholtz equation on a system of surfaces, we obtain systems of integro-differential equations in unbounded domains; the required minimax estimates are expressed via the solutions to these systems using integral operators of the potential theory in Sobolev spaces; and the BVPs are reduced to equivalent integral equation systems on multi-connected surfaces (or contours), the latter being a union of the obstacle boundary and the surfaces on which the observations are made.

We prove the unique solvability of the obtained integral equations
for any values of the wave number $k$ such that $\mbox{\rm
Im\,$k$}\geq 0,$ $k\neq 0.$

{\bf Practical importance.}  The estimation techniques elaborated in this work are of big importance
for the development of the theory of inverse acoustic and electromagnetic
wave scattering  by bounded obstacles.

The methods and results of this study may be used for estimating under uncertain conditions of the system states described by BVPs for Helmholtz equation in more complicated domains with the boundaries that stretch to infinity (for example, in a domain $K\setminus \bar\Omega$ where $K$ is a layer between two parallel planes and $\Omega$ is a bounded domain).
In general, the developed estimation methods can be applied to obtaining minimax estimates of parameters for a wide class of problems of mathematical physics.

\newpage
\makeatletter
\renewcommand{\section}{\@startsection{section}{1}%
{\parindent}{3.5ex plus 1ex minus .2ex}%
{2.3ex plus.2ex}{\normalfont\Large{\bf PART\ \ }}} \makeatother

\begin{center}
\section[
Minimax estimation of the solutions to the Helmholtz problems from
observations distributed in subdomains
]{}
\end{center}

\begin{quote}
{\bf Minimax estimation of the solutions to the Helmholtz problems
from observations distributed in subdomains
}
\end{quote}

\subsection{Notations and definitions}

Let us introduce the notations and definitions that will be used in this work.

$x=(x_1,\ldots,x_n)$ denotes a spatial variable that is varied in an open domain $D \subset\mathbb R^n;$

$dx=dx_1\ldots dx_n$ is a Lebesgue measure in $\mathbb R^n;$

$\chi(M)$ is a characteristic function of the set
$M\subset \mathbb R^n$;

$H^s(\mathbb R^n)$ is a Sobolev space of index~$s\!\!:$
$$
H^s(\mathbb R^n)=\{u\in L^2(\mathbb R^n):\,(1+|y|^2)^{s/2}\mathcal
Fu(y)\in L^2(\mathbb R^n)\},
$$
where $s\geq 0,$ $L^2(\mathbb R^n)$ is a space of square
integrable functions in $\mathbb R^2$ and $\mathcal Fu(y)$ denotes
the Fourier transform of function $u(x).$ If $s<0,$ then
$H^{-s}(\mathbb R^2)$ denotes the space dual to $H^s(\mathbb
R^2)$. Let $D$ be a domain in $\mathbb R^2$ (not necessarily
bounded) with the Lipschitz boundary $\partial D.$ Then $d\partial
D$ denotes the element of measure on contour $\partial D.$
$L^2(\partial D)$ is a space of square integrable
functions on $\partial D$; the function space
\begin{multline*}
L^{2}_{\mbox{\rm \scriptsize loc}}(D)=\left\{u\in \mathcal D'(D): \,u|_{D\cap\Omega_R}\in L^2(D\cap\Omega_R)\right.
\\\left.
\mbox{for every}\,\,R>0\,\,\mbox{such that}\,\,D\cap\Omega_R\neq \emptyset
\right\},
\end{multline*}


Introduce also the Sobolev spaces with the corresponding norms:
$$
H^{s}(D)=\{u|_{D}:\, u \in H^s(\mathbb R^2)\}\,\, (s\in\mathbb R),
$$
$$
H^{s}(\partial D)=\left\{
    \begin{array}{lc}
    u|_{\partial D}:\, u \in H^{s+1/2}(\mathbb R^2) &(s>0),\\
L^2(\partial D)& (s=0),\\
    H^{-s}(\partial D)'\,(\mbox{dual space with respect to}
    \,\, H^{-s}(\Gamma))&
    (s<0)
    \end{array}
 \right.,
$$

\begin{multline*}
H^{1}_{\mbox{\rm\scriptsize loc}}(D)=\left\{u\in\mathcal D'(D):\, u|_{D\cap\Omega_R} \in H^1(D\cap\Omega_R)\right.
\\\left.
\mbox{for every}\,\,R>0\,\,\mbox{such that}\,\,D\cap\Omega_R\neq \emptyset
\right\},
\end{multline*}
\begin{equation}\label{g21p'}
H^1(D,\Delta):=\{u\in D'(D),\,\, u\in
 H^1(D),\,\,\Delta u\in L^2(D)\},
\end{equation}
\begin{equation}\label{g21p}
H^1_{\mbox{\rm\scriptsize loc}}(D,\Delta):=\{u\in D'(D),\,\, u\in
 H^1_{\mbox{\rm\scriptsize loc}}(D),\,\,\Delta u\in L^2_{\mbox{\rm\scriptsize loc}}(D)\},
\end{equation}
\begin{multline}\label{g21p''}
H^{s}_{\mbox{\scriptsize comp}}(D):=\left\{u:\,u\in H^{s}(D),\right.\\
  \left.u \,\,\mbox{is identically zero outside some ball centered at the origin}\right\},
\end{multline}
where $\mathcal D'(D)$
is the space of distributions in $D$; here and below by $\Omega_R$ we denote the ball $\Omega_R:=\{x:|x|<R\}$; the Laplacian is taken in the sense of distributions in $D$; and $s\in \mathbb R.$

{\bf Theorem} (The trace theorem for $H^1(D)$, \cite{BIBLMcLean}, p. 102).

For any Lipschitz domain $D$ an operator $\gamma_
D : C^0(\bar D)\to
C^0(\Gamma)$ can be extended to a continuous and surjective
operator $\gamma_D: H^1_{\mbox{\rm\scriptsize loc}}(D)
\to H^{1/2}(\partial D).$

We denote by $\gamma_N$ the Neumann trace operator
$$(\gamma_Nu)(x) := (\mbox{\rm grad}\, u(x) , n(x))_{\mathbb R^n}, \,\,x \in \partial D,\,\, u \in C^1(\bar D),$$
$\gamma_
N : C^1(\bar D)\to
C^0(\partial D)$ and can be extended to a continuous and surjective
operator $\gamma_N: H^1_{\mbox{\rm\scriptsize loc}}(D,\Delta)
\to H^{-1/2}(\partial D).$  This operator will further be denoted by $\partial/\partial\nu.$



$<\cdot,\cdot>_{H^{-1/2}(\partial D)\times H^{1/2}(\partial D)}$ denotes the duality relation between spaces $H^{-1/2}(\partial D)$ and $H^{1/2}(\partial D),$
which is an extension of the inner product in $L^2(\partial D)$ in the following sense:
if $r\in L^2(\partial D),$ then the following relation holds
$$
<r,w>_{H^{-1/2}(\partial D)\times
H^{1/2}(\partial D)}=\int_{\partial D}r\bar w\,d\partial D \quad \forall w \in H^{1/2}(\partial D).\eqno{*}
$$\label{k2}


Let $H$ be a Hilbert space over the set of complex numbers $\mathbb C$ with the inner product $(\cdot,\cdot)_H$ and norm $\|\cdot\|_H.$
By $L^2(\Sigma,H)$ we denote the Bochner space composed of random\footnote {Random element $\xi$ with values in Hilbert space $H$ is considered as a function $\xi:\Sigma\to H$ imaging random events  $E\in\mathcal B$
to Borel sets in $H$ (Borel $\sigma$-algebra in $H$
is generated by open sets in $H$).} elements $\xi=\xi(\omega)$ defined on a certain probability space $(\Sigma, \mathcal B, P)$ with values in
$H$ such that
\begin{equation}\label{g2rand}
\|\xi\|_{L^2(\Sigma,H)}^2
=\int_{\Sigma}\|\xi(\omega)\|_H^2dP(\omega)<\infty.
\end{equation}
In this case there exists the Bochner integral $ \mathbb
E\xi:=\int_{\Sigma}\xi(\omega)\,dP(\omega)\in H$
which is called the mathematical expectation or the mean value of random element
$\xi(\omega)$ and satisfies the condition
\begin{equation}\label{g2rand1}
(h,\mathbb
E\xi)_H=\int_{\Sigma}(h,\xi(\omega))_H\,dP(\omega)\quad\forall
h\in H.
\end{equation}
Being applied to random variable $\xi$ this expression leads to a usual definition (value) of its mathematical expectation because the Bochner integral \eqref{g2rand} reduces to a Lebesgue integral with probability measure $dP(\omega).$

In $L^2(\Sigma,H)$ one can introduce the inner product
\begin{equation}\label{g2rand2}
(\xi,\eta)_{L^2(\Sigma,H)}:=\int_{\Sigma}(\xi(\omega),
\eta(\omega))_H\,dP(\omega)\quad\forall\xi,\eta\in L^2(\Sigma,H).
\end{equation}
Applying the sign of mathematical expectation, one can write
relationships \eqref{g2rand}$-$\eqref{g2rand2} as
\begin{equation}\label{g2rand3}
\|\xi\|_{L^2(\Sigma,H)}^2 =\mathbb E\|\xi(\omega)\|_H^2,
\end{equation}
\begin{equation}\label{g2rand4}
(h,\mathbb E\xi)_H=\mathbb E(h,\xi(\omega))_H\quad\forall h\in H,
\end{equation}
\begin{equation}\label{g2rand5}
(\xi,\eta)_{L^2(\Sigma,H)}:=\mathbb E(\xi(\omega),
\eta(\omega))_H\quad\forall\xi,\eta\in L^2(\Sigma,H).
\end{equation}
$L^2(\Sigma,H)$ equipped with norm \eqref{g2rand3} and inner product \eqref{g2rand5} is a Hilbert space.

 Consider the problem of finding a solution to the exterior Neumann problem for the Helmholtz equation.

Assume that $\Omega\in \mathbb R^2$ is a bounded domain such that $\partial \Omega=\Gamma$ is a Lipschitz contour and $\Omega_0$ is a bounded subdomain of $\mathbb R^2\setminus \bar\Omega,$ $\bar \Omega_0\subset\mathbb R^2\setminus \bar\Omega.$
 Given a function\label{f} $f$ defined in the domain $\mathbb R^2\setminus \bar\Omega$ such as $f=0$ outside $\Omega_0,$ $f\in L^2(\Omega_0)$ and a function $g\in H^{-1/2}(\Gamma),$ find $\psi\in
H^{1}_{\mbox{\rm\scriptsize loc}}((\mathbb R^2\setminus \bar\Omega),\Delta),
$
 such that
\begin{equation}\label{g21y}
-(\Delta +k^2)\psi(x)=f(x)\quad \mbox{in}\quad\mathbb R^2\setminus \bar\Omega,
\end{equation}
\begin{equation}\label{g22y}
\frac{\partial\psi}{\partial\nu}=g\,\,\mbox{on}\,\,\Gamma,
\end{equation}
\begin{equation}\label{g23y}
\frac{\partial\psi}{\partial
r}-ik\psi=o(1/r^{1/2}),\,\,r=|x|,\,\,r\to \infty
\end{equation}
with an equivalent variation formulation: find $\psi\in
H^{1}_{\mbox{\rm\scriptsize loc}}(\mathbb R^2\setminus \bar\Omega)$
 such that
 \begin{equation}\label{g2var3y}
 \int_{\mathbb R^2\setminus\bar\Omega}(\nabla \psi\nabla \overline\theta-k^2\psi\,\overline\theta)\,dx=\int_{\Omega_0}f \overline\theta\,dx+\int_{\Gamma}g
\overline\theta\,d\Gamma
 \end{equation}
for all $\theta\in H^{1}_{\mbox{\scriptsize comp}}(\mathbb
R^2\setminus \bar\Omega)$ such that $\psi$ satisfies the
Sommerfeld radiation condition \eqref{g23y}. Here we suppose that
$k$ is the wave number with $\mbox{\rm Im\,$k$}\geq 0,$ $k\neq 0.$

Consider also the following problem: find
$\psi\in H^{1}
((\Omega_R\setminus \bar\Omega),\Delta)$ such that
\begin{equation}\label{g21yx}
-(\Delta +k^2)\psi(x)=f(x)\,\, \mbox{in}\,\,\Omega_R\setminus \bar\Omega,
\end{equation}
\begin{equation}\label{g22yx}
\frac{\partial\psi}{\partial\nu}=g\,\,\mbox{on}\,\,\Gamma,
\end{equation}
\begin{equation}\label{g23yx}
\frac{\partial\psi}{\partial\nu}=M^{(1)}_k\psi\,\,\mbox{on}\,\,\Gamma_R,
\end{equation}
where $f\in L^2(\Omega_0)$ and $g\in H^{-1/2}(\Gamma),$ $M_k^{(1)}: \,H^{1/2}(\Gamma_R)\to H^{-1/2}(\Gamma_R)$ is
the Dirichlet-to-Neumann map (DtN map) defined by
\begin{equation}\label{g2M1}
(M^{(1)}_k\psi)(R,\theta):=\frac k{2\pi}\sum_{n\in\mathbb Z}\frac {H_{n}^{(1)\prime}(kR)}
{H_{n}^{(1)}(kR)}\int_0^{2\pi}\psi(R,\phi)e^{in(\theta-\phi)}d\,\phi,
\end{equation}
and $\Omega_R$ is a large disk containing $\bar\Omega$ and
$\mbox{\rm supp}\, f.$ It is known that problems
\eqref{g21y}$-$\eqref{g23y} and \eqref{g21yx}$-$\eqref{g23yx} are
equivalent in the following sense (see
\cite{BIBLCak}$-$\cite{BIBLEMasm}). If $\psi\in
H^{1}_{\mbox{\rm\scriptsize loc}}((\mathbb R^2\setminus
\bar\Omega),\Delta), $
is a solution of \eqref{g21y}$-$\eqref{g23y}, then the restriction
of $\psi$ to $\Omega_R\setminus \bar\Omega$ belongs to
$H^{1} ((\Omega_R\setminus \bar\Omega),\Delta)$
and is a solution to \eqref{g21yx}$-$\eqref{g23yx}. Conversely, if
$\psi\in H^{1} ((\Omega_R\setminus \bar\Omega),\Delta)$ is a
solution to \eqref{g21yx}$-$\eqref{g23yx}, then this solution
extended to the domain $\mathbb R^2\setminus \bar\Omega_R$ by
\begin{equation}\label{g2ext}
\psi(r_P,\theta_P)=
\frac k{2\pi}\sum_{n\in\mathbb Z}\frac {H_{n}^{(1)}(kr)}
{H_{n}^{(1)}(kR)}\int_0^{2\pi}
\psi_0(R,\phi)e^{in(\theta_P-\phi)}d\,\phi,\quad r_P\geq R,
\end{equation}
belongs to $H^{1}_{\mbox{\rm\scriptsize loc}}((\mathbb
R^2\setminus \bar\Omega),\Delta)$ and satisfies
\eqref{g21y}$-$\eqref{g23y}. Here $\psi_0(R,\phi):=\psi(R,\phi)$
is the trace of the solution to problem
\eqref{g21yx}$-$\eqref{g23yx} on $\Gamma_R$ and $(r_P,\theta_P)$
are polar coordinates of the point $P\in \mathbb R^2\setminus
\bar\Omega_R.$

To formulate an equivalent  variational setting of problem
\eqref{g21yx}$-$\eqref{g23yx}, we introduce the continuous
sesquilinear form
  $a(\cdot,
\cdot): H^1(\Omega_R\setminus\bar\Omega)\times H^1(\Omega_R\setminus\bar\Omega)\to \mathbb C$ defined as
\begin{equation}\label{g2sesq}
a(\psi,\theta):=\int_{\Omega_R\setminus\bar\Omega}(\nabla \psi\nabla \overline\theta-k^2\psi\,\overline\theta)\,dx
-\int_{\Gamma_R}M^{(1)}_k\psi\,
\overline\theta\,d\Gamma_R.
\end{equation}
Then an equivalent variational formulation of problem
\eqref{g21yx}$-$\eqref{g23yx} can be written as follows: find
$\psi\in H^1(\Omega_R\setminus\bar\Omega)$ such that
\begin{equation}\label{g2var}
a(\psi,\theta)=l(\theta) \quad\forall\theta\in H^1(\Omega_R\setminus\bar\Omega),
\end{equation}
where
\begin{equation}\label{g2semi}
l(\theta):=\int_{\Omega_0}f \overline\theta\,dx+\int_{\Gamma}g
\overline\theta\,d\Gamma
\end{equation}
is a continuous semilinear functional on $H^1(\Omega_R\setminus\bar\Omega).$

In order to set an adjoint problem of
\eqref{g21yx}$-$\eqref{g23yx} which will be used below  we
introduce a sesquilinear form
\begin{equation}\label{g2adj}
 a^*(\psi,\theta):=\int_{\Omega_R\setminus\bar\Omega}(\nabla \psi\nabla \overline\theta-\bar k^2\psi\,\overline\theta)\,dx-\int_{\Gamma_R}M^{(2)}_{\bar k}\psi\,\overline\theta\,d\Gamma_R,
\end{equation}
where $M^{(2)}_{\bar k}: \,H^{1/2}(\Gamma_R)\to H^{-1/2}(\Gamma_R)$ is
the map defined by
\begin{equation}\label{g2M2}
(M^{(2)}_{\bar k}u)(R,\theta):=\frac {\bar k}{2\pi}\sum_{z\in\mathbb Z}\frac {H_{n}^{(2)\prime}(\bar kR)}
{H_{n}^{(2)}(\bar kR)}\int_0^{2\pi}u(R,\phi)e^{in(\theta-\phi)}d\,\phi.
\end{equation}

\begin{predl}
The sesquilinear form $a^*(\psi,\theta)$ is an adjoint of $a(\psi,\theta).$
\end{predl}
\begin{proof}
Defining
$$
a_1(\psi,\theta):=\int_{\Omega_R\setminus\bar\Omega}(\nabla \psi\nabla \overline\theta-k^2\psi\,\overline\theta)\,dx \quad\forall\psi,\theta\in H^1(\Omega_R\setminus\bar\Omega)
$$
and
$$
a_2(\psi,\theta):= \int_{\Gamma_R}M^{(1)}_k\psi\,
\overline\theta\,d\Gamma_R\quad\forall\psi,\theta\in H^1(\Omega_R\setminus\bar\Omega)
$$
we have
\begin{equation}\label{g2lem1}
a(\psi,\theta)=a_1(\psi,\theta)-a_2(\psi,\theta).
\end{equation}
Obviously,
\begin{equation}\label{g2lem2}
a_1^*(\psi,\theta):=\overline{a_1(\theta,\psi)}
=\overline{\int_{\Omega_R\setminus\bar\Omega}(\nabla \theta \nabla\overline\psi-k^2\theta\,\overline\psi)\,dx}
=\int_{\Omega_R\setminus\bar\Omega}(\nabla \psi\nabla \overline\theta-\bar k^2\psi\,\overline\theta)\,dx.
\end{equation}
Taking into account that
$$
\overline{H_{n}^{(1)}(kR)}=H_{n}^{(2)}(\bar kR),\quad
\overline{H_{n}^{(1)\prime}(kR)}=H_{n}^{(2)\prime}(\bar kR),
$$
we find
$$
a_2^*(\psi,\theta):=\overline{a_2(\theta,\psi)}= \overline{\int_{\Gamma_R}M^{(1)}_k\theta\,
\overline\psi\,d\Gamma_R}=\int_{\Gamma_R}\overline{M^{(1)}_k\theta}\,
\psi\,d\Gamma_R
$$
$$
=R\int_0^{2\pi}\overline{\frac k{2\pi}\sum_{n\in\mathbb Z}\frac {H_{n}^{(1)\prime}(kR)}
{H_{n}^{(1)}(kR)}\int_0^{2\pi}
\theta(R,\alpha)e^{in(\chi-\alpha)}\,d\alpha}\,\psi(R,\chi)\,d\chi
$$
$$
=R\int_0^{2\pi}\frac{\bar k}{2\pi}\sum_{n\in\mathbb Z}\frac {\overline{H_{n}^{(1)\prime}(kR)}}
{\overline{H_{n}^{(1)}(kR)}}\int_0^{2\pi}
\overline{\theta(R,\alpha)}e^{-in(\chi-\alpha)}\,d\alpha\,\psi(R,\chi)\,d\chi
$$
$$
=R\int_0^{2\pi}\frac{\bar k}{2\pi}\sum_{n\in\mathbb Z}\frac {H_{n}^{(2)\prime}(\bar kR)}
{H_{n}^{(2)}(\bar kR)}\int_0^{2\pi}
\overline{\theta(R,\alpha)}e^{-in(\chi-\alpha)}\,d\alpha\,\psi(R,\chi)\,d\chi
$$
\begin{equation}\label{g2lem3}
=R\int_0^{2\pi}\frac{\bar k}{2\pi}\sum_{n\in\mathbb Z}\frac {H_{n}^{(2)\prime}(\bar kR)}
{H_{n}^{(2)}(\bar kR)}\int_0^{2\pi}
\psi(R,\chi)e^{in(\alpha-\chi)}\,d\chi\,\overline{\theta(R,\alpha)}
\,d\alpha
=\int_{\Gamma_R}M^{(2)}_{\bar k}\psi\,
\overline\theta\,d\Gamma_R.
\end{equation}
From \eqref{g2lem1}$-$\eqref{g2lem3} and the equality
$$
a^*(\psi,\theta)=\overline{a(\theta,\psi)}
=a_1^*(\psi,\theta)-a_2^*(\psi,\theta),
$$
we obtain the required assertion.

\end{proof}
Now we can state the variational problem adjoint of \eqref{g2var}:

{\it Given functions $f$ and $g$ introduced on page \pageref{f},  find $\psi\in H^1(\Omega_R\setminus\Omega)$ such that
\begin{equation}\label{g2varad}
a^*(\psi,\theta)=l(\theta)\quad\forall\theta\in H^1(\Omega_R\setminus\bar\Omega),
\end{equation}
where functional $l(\theta)$  is defined by \eqref{g2semi}.}

It is easy to see that variational problem \eqref{g2varad} as well as
the following problems:

(i)\,\,  find $\psi\in
H^{1}((\Omega_R\setminus \bar\Omega),\Delta)$ such that
\begin{equation}\label{g21yxa}
-(\Delta +\bar k^2)\psi(x)=f(x)\,\, \mbox{in}\,\,\Omega_R\setminus \bar\Omega,
\end{equation}
\begin{equation}\label{g22yxa}
\frac{\partial\psi}{\partial\nu}=g\,\,\mbox{on}\,\,\Gamma,
\end{equation}
\begin{equation}\label{g23yxa}
\frac{\partial\psi}{\partial\nu}=M_{\bar k}^{(2)}\psi
\,\,\mbox{on}\,\,\Gamma_R,
\end{equation}
and
(ii)\,\, find $\psi\in H^{1}_{\mbox{\rm\scriptsize loc}}((\mathbb R^2\setminus \bar\Omega),\Delta)$ such that
\begin{equation}\label{g21ya}
-(\Delta +\bar k^2)\psi(x)=f(x)\,\, \mbox{in}\,\,\mathbb R^2\setminus \bar\Omega,
\end{equation}
\begin{equation}\label{g22ya}
\frac{\partial\psi}{\partial\nu}=g\,\,\mbox{on}\,\,\Gamma,
\end{equation}
\begin{equation}\label{g23ya}
\frac{\partial\psi}{\partial
r}+i\bar k\psi=o(1/r^{1/2}),\,\,r=|x|,\,\,r\to \infty,\,\,\mbox{если}\,\, \mbox{\rm
Im\,}k\geq 0
\end{equation}
are equivalent and for any $R$ there exists a positive constant
$\alpha>0$ independent of $f$ and $g$ (but dependent on $R$) such that
\begin{equation}\label{g2cont1}
\|\psi\|_{H^1(\Omega_R\setminus\bar\Omega)}\leq \alpha(\|f\|_{H^{-1}(\Omega_0)}+\|g\|_{H^{-1/2}(\Gamma)}).
\end{equation}

If any of the data in BVPs \eqref{g21yx}$-$\eqref{g23yx} or \eqref{g21yxa}$-$\eqref{g23yxa} is random (e.g. forcing function $f$ or Neumann boundary data $g$), then the solution $\psi$ will be a random function and corresponding stochastic BVPs are formulated as follows.

Given $f\in L^2(\Sigma,L^2(\Omega_0)),$ $g\in L^2(\Sigma,H^{-1/2}(\Gamma)),$ find $\psi\in L^2(\Sigma,H^1(\Omega_R))$ such that
\begin{equation}\label{g2stvar}
\mathbb E a(\psi,\theta)=\mathbb E l(\theta) \quad\forall\theta\in L^2(\Sigma,H^1(\Omega_R\setminus\bar\Omega)),
\end{equation}
and
\begin{equation}\label{g2stvarad}
\mathbb E a^*(\psi,\theta)=\mathbb E l(\theta) \quad\forall\theta\in L^2(\Sigma,H^1(\Omega_R\setminus\bar\Omega)),
\end{equation}
where
\begin{equation}\label{g2}
\mathbb E l(\theta):=\mathbb E \left\{\int_{\Omega_0}f(x,\omega) \overline{\theta(x,\omega)}\,dx+\int_{\Gamma}g(\cdot,\omega)
\overline{\theta(\cdot,\omega)}\,d\Gamma\right\}
\end{equation}
is a continuous semilinear functional on $L^2(\Sigma,L^2(H^1(\Omega_R\setminus\bar\Omega)),$
$\mathbb E a(\cdot,
\cdot)$ and $\mathbb E a^*(\cdot,
\cdot): L^2(\Sigma,L^2(H^1(\Omega_R\setminus\bar\Omega))\times L^2(\Sigma,L^2(H^1(\Omega_R\setminus\bar\Omega))\to \mathbb C$
are continuous sesquilinear forms defined as
\begin{multline}\label{g2stsesq}
\mathbb E a(\psi,\theta):=\mathbb E\left\{\int_{\Omega_R\setminus\bar\Omega}(\nabla \psi(x,\omega)\nabla \overline{\theta(x,\omega)}-k^2\psi(x,\omega)\,\overline\theta{(x,\omega)}\,dx
\right.\\ \left.-\int_{\Gamma_R}M^{(1)}_k\psi(\cdot,\omega)\,
\overline{\theta(\cdot,\omega)}\,d\Gamma_R\right\}
\end{multline}
and
\begin{multline}\label{g2stsesqad}
\mathbb E a^*(\psi,\theta):=\mathbb E\left\{\int_{\Omega_R\setminus\bar\Omega}(\nabla \psi(x,\omega)\nabla \overline{\theta(x,\omega)}-\bar k^2\psi(x,\omega)\,\overline{\theta(x,\omega)}\,dx
\right.\\ \left.-\int_{\Gamma_R}M^{(2)}_{\bar k}\psi(\cdot,\omega)\,\overline{\theta(\cdot,\omega)}\,d\Gamma_R\right\}.
\end{multline}

It is known that problems \eqref{g2stvar} and \eqref{g2stvarad} have unique solutions
and there exists a positive constant
$\alpha>0$ independent of $f$ and $g$ such that
\begin{equation}\label{g2alpha1}
\|\psi\|_{L^2(\Sigma,H^1(\Omega_R))}\leq \alpha(\|f\|_{L^2(\Sigma,L^2(\Omega_0))}+\|g\|_{L^2(\Sigma,H^{-1/2}(\Gamma))}).
\end{equation}
Such problems were investigated in \cite{BIBLElm} and \cite{BIBLEGhan}, including the construction of finite element methods of their numerical solution.

Problem \eqref{g2stvar} is equivalent to the following ones: find
$\psi\in L^2(\Sigma,H^{1} ((\Omega_R\setminus
\bar\Omega),\Delta))$ such that
\begin{equation}\label{g21yxl}
-(\Delta +k^2)\psi(x,\omega)=f(x,\omega)\,\, \mbox{in}\,\,\Omega_R\setminus \bar\Omega,
\end{equation}
\begin{equation}\label{g22yxl}
\frac{\partial\psi(\cdot,\omega)}{\partial\nu}=g(\cdot,\omega)\,\,\mbox{on}\,\,\Gamma,
\end{equation}
\begin{equation}\label{g23yxl}
\frac{\partial\psi(\cdot,\omega)}{\partial\nu}=M^{(1)}_k\psi(\cdot,\omega)
\,\,\mbox{on}\,\,\Gamma_R,
\end{equation}
or
find $\psi\in L^2(\Sigma,H^{1}
(\Omega_R\setminus \bar\Omega))$
satisfying
\begin{equation}\label{g2stvar1}
a(\psi(\cdot,\omega),\theta)=l(\theta) \quad\forall\theta\in H^1(\Omega_R\setminus\bar\Omega).
\end{equation}
Analogously, problem \eqref{g2stvarad} is equivalent to the
following problems: find $\psi\in L^2(\Sigma,H^{1}
((\Omega_R\setminus \bar\Omega),\Delta))$ such that
\begin{equation}\label{g21yxlad}
-(\Delta +\bar k^2)\psi(x,\omega)=f(x,\omega)\,\, \mbox{in}\,\,\Omega_R\setminus \bar\Omega,
\end{equation}
\begin{equation}\label{g22yxlad}
\frac{\partial\psi(\cdot,\omega)}{\partial\nu}=g(\cdot,\omega)\,\,\mbox{on}\,\,\Gamma,
\end{equation}
\begin{equation}\label{g23yxlad}
\frac{\partial\psi(\cdot,\omega)}{\partial\nu}=M^{(2)}_{\bar k}\psi(\cdot,\omega)\,\,\mbox{on}\,\,\Gamma_R,
\end{equation}
or
find $\psi\in L^2(\Sigma,H^{1}
(\Omega_R\setminus \bar\Omega))$
satisfying
\begin{equation}\label{g2stvarad1}
a^*(\psi(\cdot,\omega),\theta)=l(\theta) \quad\forall\theta\in H^1(\Omega_R\setminus\bar\Omega).
\end{equation}

The right-hand sides in (\ref{g21yxl})$-$(\ref{g2stvarad1}) are considered for every realization of random fields $f(\cdot,\omega)$ and $g(\cdot,\omega)$ which belong with probability
 $1$ to the spaces $L^2(\Omega_0)$ and $L^2(\Gamma)$, respectively, and the equalities are satisfied almost certainly.

\subsection{Statement of the estimation problem}

Consider the exterior Neumann problem for the Helmholtz equation: find a distribution $\varphi\in\mathcal D'(\mathbb R^2\setminus \bar\Omega)$ such that
\begin{equation}\label{g20}
\varphi\in H^{1}_{\mbox{\rm\scriptsize loc}}((\mathbb R^2\setminus \bar\Omega),\Delta),
\end{equation}
\begin{equation}\label{g21}
-(\Delta+ k^2)\varphi(x)=f(x)\quad \mbox{in}\quad\mathbb R^2\setminus\bar\Omega,
\end{equation}
\begin{equation}\label{g22}
\frac{\partial\varphi}{\partial\nu}=g\,\,\mbox{on}\,\,\Gamma,
\end{equation}
\begin{equation}\label{g23}
\frac{\partial\varphi}{\partial
r}-ik\varphi=o(1/r^{1/2}),\,\,r=|x|,\,\,r\to \infty,\,\, \mbox{\rm
Im\,}k\geq 0,
\end{equation}
where $k$ is the wave number with $\mbox{\rm Im\,$k$}\geq 0,$ $f$
is a source term distributed in bounded subdomain\footnote{This
means that the function $f$ is defined in the domain $\mathbb
R^2\setminus \bar\Omega,$ $f=0$ outside $\Omega_0,$ and $f\in
L^2(\Omega_0).$} $\Omega_0$ in $\mathbb R^2\setminus\bar\Omega$,
$f\in L^2(\Omega_0),$ and
 $g\in L^2(\Gamma).$


 BVP \eqref{g20}$-$\eqref{g23} simulates, in particular,
 acoustic or electromagnetic scattering
 from
 an infinite sound-hard (perfectly conducting) cylinder with cross-section $\Omega.$

Denote by $G_0$ the set of pairs of functions $(\tilde f,\tilde g)$
satisfying the inequality
\begin{equation} \label{g28}
\int_{\Omega_0} Q_1(\tilde f-f_0)(x)\overline{(\tilde f -f_0)}(x)\, dx +\int_{\Gamma} Q_2(\tilde g-g_0)\overline{(\tilde g -g_0)}\,d\Gamma\leq
1,
\end{equation}
 and by $G_1$ the set of random functions
$\tilde \xi(\cdot)=(\tilde \xi_1(\cdot),
\ldots,\tilde \xi_{m}(\cdot))$ defined on $\Omega_1\times\cdots\times\Omega_k$ with integrable second moments $\mathbf E|\tilde\xi_k(x)|^2$ satisfying conditions
\begin{equation} \label{g27}
\mathbf E\tilde \xi_{k}(x)=0,\,\,k= \overline{1,m}.
\end{equation}
\begin{equation} \label{g29}
\sum_{k=1}^{m} \int_{\Omega_k} \mathbf E|\tilde \xi_k(x)|^2
r_k^2(x)\, dx\leq 1,
\end{equation}
where $f_0$ is defined in the domain $\mathbb R^2\setminus \bar\Omega,$ $f_0=0$ outside $\Omega_0,$ $f_0\in L^2(\Omega_0)$ and $g_0\in L^2(\Gamma),$ $f_0$ and $g_0$ are prescribed functions, $Q_1$ and $Q_2$ are
Hermitian operators
in $L^2(\Omega_0)$ and $L^2(\Gamma),$ respectively, for which there exist bounded inverse operators
$Q_1^{-1}$ and $Q_2^{-1},$ and
$\tilde r_k(x), \,\, k= \overline{1,m},$
 are nonvanishing functions continuous on sets $\bar \Omega_k$.

We suppose that functions $f(x)$ and $g(x)$ in equations
(\ref{g21}) and (\ref{g22})
are not known exactly;
it is known only that
$(f,g)\in G_0.$

Assume that in subdomains $\Omega_k,$ $k=\overline{1,m},$ of domain $\Omega$ the following functions are observed
\begin{equation} \label{g25}
y_{k}(x) = \int_{\Omega_k} g_{k}(x,y)\varphi(y)\,dy
+ \xi_{k}(x,\omega), \quad x \in \Omega_k, \quad k=
\overline{1,m},
\end{equation}
where $\varphi$ is a solution of BVP (\ref{g20})$-$(\ref{g23}), $g_k(x,y) \in L^2(\Omega_k \times \Omega_k)$
 are
prescribed functions and $\xi_k(x,\omega)$ are the choice functions
of random fields $\xi_k(x)$ with unknown second moments such that $ \xi(\cdot)=( \xi_1(\cdot),
\ldots,\xi_{m}(\cdot))\in G_1.$

Let $l_0$ be a given function  defined in a bounded subdomain $\omega_0\subset\Omega$ belonging to $L^2(\omega_0).$

The estimation problem consists in the following. From
observations (\ref{g25}) of the state $\varphi(x)$ of the system
described by BVP (\ref{g20})$-$(\ref{g23}) under conditions
(\ref{g27})$-$(\ref{g29}) it is necessary to estimate the value of
the linear functional
\begin{equation} \label{g210}
l(\varphi)= \int_{\omega_0}\overline{ l_0(x)}\varphi(x)\,dx
\end{equation}
in the class of the estimates linear with respect to observations  which have
the form
\begin{equation}\label{g2linest}
\widehat {l(\varphi)}= \sum_{k=1}^{m}\int_{\Omega_k}
\overline{ u_k(x)}y_k(x)\, dx+c,
\end{equation}
where $u_k \in L^2(\Omega_k),\,\, k= \overline{1,m},\,\,c
\in \mathbb C.$

Denote by $u=(u_1, \ldots , u_{m})$ an element belonging to
$H:=L^2(\Omega_1)\times\ldots\times L^2(\Omega_k)$.
\begin{predll}
An estimate $\widehat {\widehat {l(\varphi)}}$ is called a minimax estimate of the $l(\varphi)$ if an element $\hat u=(\hat u_1,\ldots, \hat u_m)\in H$ and a number  $\hat c\in \mathbb C$ are determined
from the condition
$$
\sup_{(\tilde f,\tilde g) \in G_0,\, \tilde \xi \in G_1}
\mathbb E[l(\tilde\varphi)-\widehat {l(\tilde\varphi)}]^2\to \inf_{u\in H,\,c\in \mathbb C}.
$$
Here $\tilde\varphi$ is a solution to problem
(\ref{g20})$-$(\ref{g23}) when $f=\tilde f,$ $g=\tilde g,$ and
$\widehat {l(\tilde\varphi)}=\sum_{k=1}^{m}\int_{\Omega_k}
\overline{ u_k(x)}\tilde y_k(x)\, dx+c,$ $\tilde y_{k}(x) =
\int_{\Omega_k} g_{k}(x,y)\tilde \varphi(y)\,dy + \tilde
\xi_{k}(x),$ $x \in \Omega_k, \,\, k= \overline{1,m},$

 The
quantity
\begin{equation}\label{g210p}
\sigma := \{\sup_{(\tilde f,\tilde g) \in G_0,\, \tilde \xi \in G_1}
\mathbb E[l(\tilde\varphi)-\widehat {\widehat {l(\tilde\varphi)}}]^2\}^{1/2}
\end{equation}
is called the error of the minimax estimation of $l(\varphi).$
\end{predll}

Thus, the minimax estimate is an estimate minimizing the maximal
mean-square estimation error calculated for the ``worst''
implementation of perturbations.


\subsection{Reduction of the estimation problem to an optimal control problem}

\begin{predl}\label{g2lemma2}
The problem of finding the minimax estimate of
$l(\varphi)$ is equivalent to the problem of optimal control of a system described by a BVP
\begin{equation}\label{g20g}
z(\cdot;u)\in H^{1}
((\Omega_R\setminus \bar\Omega),\Delta),
\end{equation}
\begin{equation}\label{g21g}
-(\Delta +\bar
k^2)z(x;u)=\chi_{\omega_0}(x)l_0(x)-\sum_{k=1}^m
\chi_{\Omega_k}(x)\int_{\Omega_k}\overline{ g_{k}(\eta,x)}u_k(\eta)\,d\eta
\,\, \mbox{in}\,\,\Omega_R\setminus\bar \Omega,
\end{equation}
\begin{equation}\label{g22g}
\frac{\partial z(\cdot;u)
}{\partial\nu}=0\,\,\mbox{on}\,\,\Gamma,
\end{equation}
\begin{equation}\label{g23g}
\frac{\partial z(\cdot;u)}{\partial
r}=M^{(2)}_{\bar k}z(\cdot;u)\,\,\mbox{on}\,\,\Gamma_R
\end{equation}
with the quality criterion
$$
I(u):=\int_{\Omega_0} Q_1^{-1}z(x;u)\overline{z(x;u)}\, dx +\int_{\Gamma} Q_2^{-1}z(\cdot;u)\overline{z(\cdot;u)}\,d\Gamma
$$
\begin{equation} \label{g24g}
+\sum_{k=1}^{m} \int_{\Omega_k} r_k^{-2}(x)|u_k(x)|^2\, dx\,\to \inf_{ u\in H},
\end{equation}
where $R$ is chosen so that $\bar\Omega_i\subset\Omega_R\setminus\bar \Omega,$ $i=\overline{0,m},$ $\bar\omega_0\subset\Omega_R\setminus\bar \Omega.$
\end{predl}
\begin{proof}
Taking into account \eqref{g25},  \eqref{g210}, and \eqref{g2linest}, we obtain
$$
l(\tilde\varphi)-\widehat{l(\tilde\varphi)}
=\int_{\omega_0}\overline{l_0(x)}\tilde\varphi(x)\,dx
-\sum_{k=1}^{m}\int_{\Omega_k}
\overline{u_k(x)}\tilde y_k(x)\,dx-c
$$
$$
=\int_{\omega_0}\overline{l_0(x)}\tilde\varphi(x)\,dx
-\sum_{k=1}^{m}\int_{\Omega_k}
\overline{u_k(x)}\int_{\Omega_k} g_{k}(x,y)\tilde \varphi(y)\,dy\,dx  -\sum_{k=1}^{m}\int_{\Omega_k}
\overline{u_k(x)}\tilde\xi_k(x)\,dx-c
$$
$$
=\int_{\omega_0}\overline{l_0(x)}\tilde\varphi(x)\,dx
-\sum_{k=1}^{m}\int_{\Omega_k}
\int_{\Omega_k} g_{k}(y,x)\overline{u_k(y)}\,dy\, \tilde\varphi(x)\,dx  -\sum_{k=1}^{m}\int_{\Omega_k}
\overline{u_k(x)}\tilde\xi_k(x)\,dx-c
$$
$$
=\int_{\Omega_R\setminus\bar\Omega}
\overline{\left(\chi_{\omega_0}(x)l_0(x)
-\sum_{k=1}^{m}
\chi_{\Omega_k}(x)\int_{\Omega_k} \overline{g_{k}(y,x)}u_k(y)\,dy\right)}\tilde\varphi(x)\,dx
$$
\begin{equation}\label{g24gg}
  -\sum_{k=1}^{m}\int_{\Omega_k}
\overline{u_k(x)}\tilde\xi_k(x)\,dx-c.
\end{equation}
For any fixed $u=(u_1, \ldots , u_{m})\in H$ introduce the
function $z(x;u)$ as a unique solution of problem
\eqref{g20g}$-$\eqref{g23g}. According to the equivalent
variational formulation of this problem, it means that $z(x;u)$
satisfies the integral identity
$$
a^*(z(\cdot;u),\theta)=\int_{\Omega_R\setminus\bar\Omega}(\nabla z(x;u)\nabla\overline\theta
(x)-\bar k^2z(x;u)\,\overline\theta(x))\,dx-\int_{\Gamma_R}M_{\bar k}^{(2)}z(\cdot;u)\,\overline\theta\,d\Gamma_R
$$
\begin{equation}\label{g25g}
=\int_{\Omega_R\setminus\bar\Omega}\left(\chi_{\omega_0}(x)l_0(x)
-\sum_{k=1}^{m}
\chi_{\Omega_k}(x)\int_{\Omega_k} \overline{g_{k}(y,x)}u_k(y)\,dy\right)\overline{\theta(x)}\,dx
\quad\forall\theta\in H^1(\Omega_R\setminus\bar\Omega).
\end{equation}
Set $\theta=\tilde\varphi$ in \eqref{g25g}. Then we obtain
$$
a^*(z(\cdot;u),\tilde\varphi)=\int_{\Omega_R\setminus\bar\Omega}
(\nabla z(x;u)\nabla\overline{\tilde\varphi(x)}
-\bar k^2z(x;u)\,\overline{\tilde\varphi(x)}\,dx-\int_{\Gamma_R}M_{\bar k}^{(2)}z(\cdot;u)\,\overline{\tilde\varphi}\,d\Gamma_R
$$
\begin{equation}\label{g26g}
=\int_{\Omega_R\setminus\bar\Omega}\left(\chi_{\omega_0}(x)l_0(x)
-\sum_{k=1}^{m}
\chi_{\Omega_k}(x)\int_{\Omega_k} \overline{g_{k}(y,x)}u_k(y)\,dy\right)\overline{{\tilde\varphi}}\,dx.
\end{equation}
On the other hand, since $\tilde\varphi$ is a solution of problem \eqref{g20}$-$\eqref{g23}
with $f=\tilde f$ and $g=\tilde g,$ setting $\psi=\varphi$ and $\theta=z(\cdot;u)$ in \eqref{g2sesq}, we find
$$
a(\tilde\varphi,z(\cdot;u))=\int_{\Omega_R\setminus\bar\Omega}
(\nabla \tilde\varphi(x)\nabla\overline{z(x;u)}
-k^2\tilde\varphi(x)\,\overline{z(x;u)}\,dx-\int_{\Gamma_R}M_{ k}^{(1)}\tilde\varphi\,\overline{z(\cdot;u)}\,d\Gamma_R
$$
\begin{equation}\label{g27g}
=\int_{\Omega_0}\tilde f(x)\overline{z(x;u)}\,dx
+\int_{\Gamma}\tilde g\overline{z(\cdot;u)}\,d\Gamma.
\end{equation}
 By Lemma 1.1,
   $
\overline{a^*(z,\tilde\varphi)}=a(\tilde\varphi,z).$ This identity and
\eqref{g24g}, \eqref{g26g}, and \eqref{g27g} imply
$$
l(\tilde\varphi)-\widehat{l(\tilde\varphi)}=
\int_{\Omega_R\setminus\bar\Omega}
\overline{\left(\chi_{\omega_0}(x)l_0(x)
-\sum_{k=1}^{m}
\chi_{\Omega_k}(x)\int_{\Omega_k} \overline{g_{k}(y,x)}u_k(y)\,dy\right)}\tilde\varphi(x)\,dx
$$
\begin{equation}\label{g28g}
=\int_{\Omega_0}\tilde f(x)\overline{z(x;u)}\,dx
+\int_{\Gamma}\tilde g\overline{z(\cdot;u)}\,d\Gamma
-\sum_{k=1}^{m}\int_{\Omega_k}
\overline{u_k(x)}\tilde\xi_k(x)\,dx-c.
\end{equation}

Taking into consideration the relationship $\mathbf {D}\eta=\mathbf
E|\eta-\mathbf E\eta|^2=\mathbf E|\eta|^2-|\mathbf E\eta|^2$ that couples dispersion $\mathbf D\eta$
of the complex random variable $\eta=\eta_1+i\eta_2$ and its
expectation $\mathbf E\eta=\mathbf E\eta_1+i\mathbf E\eta_2,$ we obtain from the last formulas
$$
\mathbf E\left|l(\tilde
\varphi)-\widehat {l(\tilde \varphi)}\right|^2
=\left|\int_{\Omega_0}\tilde f(x)\overline{z(x;u)}\,dx
+\int_{\Gamma}\tilde g\overline{z(\cdot;u)}\,d\Gamma \right|^2 +\mathbf
E\left|\sum_{k=1}^{m}\int_{\Omega_k}
\overline{u_k(x)}\,\tilde\xi_k(x)\,dx\right|^2.
$$
Therefore,
$$ \inf_{c \in
\mathbb C}\sup_{(\tilde f,\tilde g)\in G_0, \tilde \xi\in
G_1} \mathbf E|l(\tilde
\varphi)-\widehat {l(\tilde \varphi)}|^2=
$$
$$
=\inf_{c \in \mathbb C}\sup_{(\tilde f,\tilde g)\in G_0}\left|\int_{\Omega_0}\tilde f(x)\overline{z(x;u)}\,dx
+\int_{\Gamma}\tilde g\overline{z(\cdot;u)}\,d\Gamma-c\right|^2
$$
\begin{equation}\label{g2exh}
+ \sup_{ \tilde \xi\in
G_1}\mathbf E\left|\sum_{k=1}^{m}\int_{\Omega_k}
\overline{u_k(x)}\,\tilde\xi_k(x)\,dx\right|^2.
\end{equation}
In order to calculate the first term on the right-hand side of (\ref{g2exh})
make use of the generalized Cauchy$-$Bunyakovsky inequality in
(\ref{g28}). We have
$$
\inf_{c \in \mathbb C}\sup_{(\tilde f,\tilde g)\in G_0}\left|\int_{\Omega_0}\tilde f(x)\overline{z(x;u)}\,dx
+\int_{\Gamma}\tilde g\overline{z(\cdot;u)}\,d\Gamma-c\right|^2
$$
\begin{multline*}
=\inf_{c \in \mathbb C}\sup_{(\tilde f,\tilde g)\in G_0}\left|\int_{\Omega_0}(\tilde f(x)-f_0(x))\overline{z(x;u)}\,dx
+\int_{\Gamma}(\tilde g-g_0)\overline{z(\cdot;u)}\,d\Gamma\right.\\+
\left.\int_{\Omega_0}f_0(x)\overline{z(x)}\,dx
+\int_{\Gamma}g_0\overline{z(\cdot;u)}\,d\Gamma-c\right|^2
\end{multline*}
\begin{multline*}
\leq \left\{\int_{\Omega_0}Q_1^{-1}z(x)\overline{z(x;u)}\,dx
+\int_{\Gamma}Q_2^{-1}z(\cdot;u)\overline{z(\cdot;u)}\,d\Gamma\right\}
\\
\times \left\{\int_{\Omega_0}Q_1(\tilde f-f_0)(x)\overline{(\tilde f(x)-f_0(x))}\,dx
+\int_{\Gamma}Q_2(\tilde g-g_0)\overline{(\tilde g-g_0)}\,d\Gamma\right\}
\end{multline*}
\begin{equation}\label{g2ex2h}
\leq
\int_{\Omega_0}Q_1^{-1}z(x;u)\overline{z(x;u)}\,dx
+\int_{\Gamma}Q_2^{-1}z(\cdot;u)\overline{z(\cdot;u)}\,d\Gamma.
\end{equation}
The direct substitution shows that that inequality
(\ref{g2ex2h}) is transformed to an equality on the element  $(\tilde
f^{(0)}(\cdot),\tilde g^{(0)}),$
where
$$
\tilde f^{(0)}(x):=\frac 1{d}Q_1^{-1}z(x;u)+f_0(x),
$$
$$
\tilde g^{(0)}:=\frac 1{d}Q_2^{-1}z(\cdot;u)+g_0,
$$
$$
d={\Bigl(\int_{\Omega_0}Q_1^{-1}z(x;u)\overline{z(x;u)}\,dx
+\int_{\Gamma}Q_2^{-1}z(\cdot;u)\overline{z(\cdot;u)}
\,d\Gamma\Bigr)^{1/2}}.
$$
Therefore
$$\inf_{c\in \mathbb C} \sup_{(\tilde f,\tilde g) \in G_0,\, \tilde \xi \in G_1}\left|\int_{\Omega_0}(\tilde f(x)-f_0(x))\overline{z(x;u)}\,dx
+\int_{\Gamma}(\tilde g-g_0)\overline{z(\cdot;u)}\,d\Gamma-c\right|^2
$$
\begin{equation}\label{g2lem4}
=\int_{\Omega}Q_1^{-1}z(x;u)\overline{ z(x;u)}\,dx
+\int_{\Gamma}Q_2^{-1}z(\cdot;u)\overline{ z(\cdot;u)}\,d\Gamma
\end{equation}
with
$$
c=\int_{\Omega_0}\overline{z(x;u)}
f_0(x)\,dx+\int_\Gamma\overline{z(\cdot;u)}
g_0\,d\Gamma.
$$
In order to calculate the second term on the right-hand side of  (\ref{g2exh}),
note that the Cauchy$-$Bunyakovsky inequality and (\ref{g29})
yield
$$
\left|\sum_{k=1}^{m}\int_{\Omega_k}
\overline{u_k(x)}\,\tilde\xi_k(x)\,dx\right|^2\leq
$$
$$
\sum_{k=1}^{m}\int_{\Omega_k}
r_k^{-2}(x)|u_k(x)|^2\,dx\sum_{k=1}^{m}
\int_{\Omega_k}r_k^{2}(x)|\tilde\xi_k(x)|^2\,dx,
$$
the latter implies
$$
\sup_{\tilde \xi\in G_1}\mathbf E\left|\sum_{k=1}^{m}\int_{\Omega_k}
\overline{u_k(x)}\,\tilde\xi_k(x)\,dx\right|^2\leq
\sum_{k=1}^{m}\int_{\Omega_k}
r_k^{-2}(x)|u_k(x)|^2.
$$
However\label{g2L1}
$$
\mathbf E\left|\sum_{k=1}^{m}\int_{\Omega_k}
\overline{u_k(x)}\,\tilde\xi^{(0)}_k(x)\,dx\right|^2
=\sum_{k=1}^{m}\int_{\Omega_k}
r_k^{-2}(x)|u_k(x)|^2 ,
$$
where
$$
\tilde\xi^{(0)}_k(x)=\frac{\nu r_k^{-2}(x)u_k(x)}{\{\sum_{k=1}^{m}\int_{\Omega_k}
r_k^{-2}(x)|u_k(x)\}^{1/2}},\,\,x\in\Omega_k,\,\,k=\overline{1,m},
$$
$$
 \tilde\xi^{(0)}=(\tilde\xi^{(0)}_1(\cdot),\ldots,
\tilde\xi^{(0)}_m(\cdot))\in
G_1,
$$
and $\nu$ is a random variable with $\mathbf E\nu=0$ and
$\mathbf E|\nu|^2=1.$
Therefore
\label{g2L2.2}
\begin{equation}\label{g2lieh}
\sup_{\tilde \xi\in G_1}\mathbf E\left|\sum_{k=1}^{m}\int_{\Omega_k}
\overline{u_k(x)}\,\tilde\xi_k(x)\,dx\right|^2
=\sum_{k=1}^{m}\int_{\Omega_k}
r_k^{-2}(x)|u_k(x)|^2,
\end{equation}
which proves the required assertion.
The validity of Lemma 2.2 follows now from relationships (\ref{g2exh}), (\ref{g2lem4}), and (\ref{g2lieh}).
\end{proof}

\subsection{Representation of minimax estimates and estimation errors}

In the course of the proof of Theorem 1.1 below, we show that the
solution to the optimal control problem
\eqref{g2lemma2}$-$\eqref{g24g} (and therefore, the determination
of the minimax estimate, in line with Theorem 1.2) is reduced to
the solution of a certain integro-differential equation system.
Namely, the following statement holds.

\begin{pred}
The minimax estimate of $l(\varphi)$ has the form
\begin{equation}\label{g2rt}
\widehat{\widehat{l(\varphi)}}=\sum_{k=1}^m \int_{\Omega_k}
\overline{\hat u_k(x)}y_k(x)\, dx+\hat c,
\end{equation}
where
\begin{equation}\label{g2rt1}
\hat u_k(x)=r_k^2(x)\int_{\Omega_k}g_{k}(x,y)p(y)\,dy,
\,\,k=\overline{1,m},\,\,\hat
c=\int_{\Omega_0}\overline{z(x)}
f_0(x)\,dx+\int_\Gamma\overline{z}
g_0\,d\Gamma,
\end{equation}
and functions $z$ and $p$ are determined from the solution to the following problem:
\begin{equation}\label{g20gr}
z\in H^{1}
((\Omega_R\setminus \bar\Omega),\Delta),
\end{equation}
\begin{equation}\label{g21gr}
-(\Delta +\bar
k^2)z(x)=\chi_{\omega_0}(x)l_0(x)-\sum_{k=1}^m
\chi_{\Omega_k}(x)\int_{\Omega_k}\overline{ g_{k}(\eta,x)}\hat u_k(\eta)\,d\eta
\,\, \mbox{in}\,\,\Omega_R\setminus\bar \Omega,
\end{equation}
\begin{equation}\label{g22gr}
\frac{\partial z
}{\partial\nu}=0\,\,\mbox{on}\,\,\Gamma,
\end{equation}
\begin{equation}\label{g23gr}
\frac{\partial z}{\partial
r}=M^{(2)}_{\bar k}z\,\,\mbox{on}\,\,\Gamma_R,
\end{equation}
\begin{equation}\label{g20yxr}
p\in H^{1}
((\Omega_R\setminus \bar\Omega),\Delta),
\end{equation}
\begin{equation}\label{g21yxr}
-(\Delta +k^2)p(x)=\chi_{\Omega_0}(x)Q_1^{-1}z(x)\,\, \mbox{in}\,\,\Omega_R\setminus \bar\Omega,
\end{equation}
\begin{equation}\label{g22yxr}
\frac{\partial p}{\partial\nu}=Q_2^{-1}z,\,\,\mbox{on}\,\,\Gamma,
\end{equation}
\begin{equation}\label{g23yxr}
\frac{\partial p}{\partial\nu}=M^{(1)}_k p\,\,\mbox{on}\,\,\Gamma_R.
\end{equation}
Problem \eqref{g20gr}$-$\eqref{g23yxr} is uniquely solvable. The restrictions
of the solutions of this problem corresponding to $R=R_1$ and $R=R_2$ on $\Omega_{R=\mbox{\rm\scriptsize min}\{R_1,R_2\}}$ coincide with the solution corresponding to
$R=\mbox{\rm min}\{R_1,R_2\}.$

 The estimation error $\sigma$ is determined by
 the formula\footnote{In the proof of this theorem, we show that the value of $\int_{\omega_0}\overline{ l_0(x)}p(x)\,dx$ is real.}
\begin{equation}\label{k1k4}
\sigma=l(p)^{1/2}=\left(\int_{\omega_0}\overline{ l_0(x)}p(x)\,dx\right)^{1/2}.
\end{equation}
\end{pred}
\begin{proof}
Let us show first that the optimal control problem
(\ref{g20g})$-$(\ref{g24g}) is uniquely solvable; that is, there
exists one and only one element $\hat u \in H,$ at which
functional (\ref{g24g}) attains the minimum value, $I(\hat
u)=\inf_{u \in H}I(u)$.

It is easy to see that the solution $z(x;u)$ of BVP
(\ref{g20g})$-$(\ref{g23g}) can be represented as
\begin{equation}\label{g2ieq}
z(x;u)=z_0(x)+\tilde z(x;u),
\end{equation}
where
$z_0(x)$ and $\tilde z(x;u)$ are solutions of the following BVPs
\begin{equation}\label{g20g2}
z_0\in H^{1}
((\Omega_R\setminus \bar\Omega),\Delta),
\end{equation}
\begin{equation}\label{g21g2}
-(\Delta +\bar
k^2)z_0(x)=\chi_{\omega}(x)l_0(x)
\,\, \mbox{in}\,\,\Omega_R\setminus\bar \Omega,
\end{equation}
\begin{equation}\label{g22g2}
\frac{\partial z_0
}{\partial\nu}=0\,\,\mbox{on}\,\,\Gamma,
\end{equation}
\begin{equation}\label{g23g2}
\frac{\partial z_0}{\partial
r}=M^{(2)}_{\bar k}z_0\,\,\mbox{on}\,\,\Gamma_R
\end{equation}
and
\begin{equation}\label{g20g3}
\tilde z(\cdot;u)\in H^{1}
((\Omega_R\setminus \bar\Omega),\Delta),
\end{equation}
\begin{equation}\label{g21g3}
-(\Delta +\bar
k^2)\tilde z(x;u)=-\sum_{k=1}^m
\chi_{\Omega_k}(x)\int_{\Omega_k}\overline{ g_{k}(\eta,x)}u_k(\eta)\,d\eta
\,\, \mbox{in}\,\,\Omega_R\setminus\bar \Omega,
\end{equation}
\begin{equation}\label{g22g3}
\frac{\partial \tilde z(\cdot;u)
}{\partial\nu}=0\,\,\mbox{on}\,\,\Gamma,
\end{equation}
\begin{equation}\label{g23g3}
\frac{\partial \tilde z(\cdot;u)}{\partial
r}=M^{(2)}_{\bar k}\tilde z(\cdot;u)\,\,\mbox{on}\,\,\Gamma_R.
\end{equation}
Using representations (\ref{g2ieq}) for $z(t;u),$ write functional $I(u)$ in the form
$$
I(u)=\int_{\Omega_0} Q_1^{-1}z(x;u)\overline{z(x;u)}\, dx +\int_{\Gamma} Q_2^{-1}z(\cdot;u)\overline{z(\cdot;u)}\,d\Gamma+
\sum_{k=1}^{m} \int_{\Omega_k} r_k^{-2}(x)|u_k(x)|^2\, dx
$$
$$
=\int_{\Omega_0} Q_1^{-1}(z_0(x)+\tilde z(x;u))\overline{(z_0(x)+\tilde z(x;u))}\, dx +\int_{\Gamma} Q_2^{-1}(z_0+\tilde z(\cdot;u))\overline{(z_0+\tilde z(\cdot;u))}\,d\Gamma
$$
$$
\sum_{k=1}^{m} \int_{\Omega_k} r_k^{-2}(x)|u_k(x)|^2\, dx
=\tilde I(u)+L(u)+c_0,
$$
where
$$
\tilde I(u):=\int_{\Omega_0} Q_1^{-1}\tilde z(x;u)\overline{\tilde z(x;u)}\, dx +\int_{\Gamma} Q_2^{-1}\tilde z(\cdot;u)\overline{\tilde z(\cdot;u)}\,d\Gamma+
\sum_{k=1}^{m} \int_{\Omega_k} r_k^{-2}(x)|u_k(x)|^2\, dx
$$
is a quadratic functional in the space $H$ which corresponds to a semi-bilinear continuous Hermitian form \label{g2cont2}
\begin{multline}\label{g2rep2a}
\pi(u,v)=\int_{\Omega_0} Q_1^{-1}\tilde z(x;u)\overline{\tilde z(x;v)}\, dx +\int_{\Gamma} Q_2^{-1}\tilde z(\cdot;u)\overline{\tilde z(\cdot;v)}\,d\Gamma\\+
\sum_{k=1}^{m} \int_{\Omega_k} r_k^{-2}(x)u_k(x)\overline{v_k(x)}\, dx
\end{multline}
on $H\times H$ and satisfies
\begin{equation}\label{g2rep3}
\tilde I(u)\geq a \|u\|^2_{H_0}\,\,\forall u\in
H,\,\,a=\mbox{const};
\end{equation}
note that
$$
L(u):=2\mbox{Re}\int_{\Omega_0} Q_1^{-1}\tilde z(x;u)\overline{z_0(x)}\, dx+2\mbox{Re}\int_{\Gamma} Q_2^{-1}\tilde z(\cdot;u)\overline{z_0}\, d\Gamma
$$
is a linear continuous functional in $H$
and
$$
c_0:=\int_{\Omega_0} Q_1^{-1}z_0(x)\overline{z_0(x)}\, dx +\int_{\Gamma} Q_2^{-1}z_0\overline{z_0}\,d\Gamma.
$$
Prove, for example, the continuity of form (\ref{g2rep2a}); namely, the inequality
\begin{equation}\label{g2rep5}
|\pi(v,w)|\leq c  \|v\|_{H_0}  \|w\|_{H_0} \quad\forall v,w\in
V,\,\,c=\mbox{const}
\end{equation}
(the continuity of linear functional $L(v)$ is proved in a similar manner).

Using estimate (\ref{g2cont1}) obtained above and the Cauchy$-$Bunyakovsky inequality
we have
$$
|\pi(u,v)| \leq\left(\int_{\Omega_0}Q_1^{-1}\tilde z(x;u)\overline{\tilde
z(x;u)}\,dx\right)^{1/2}\left(\int_{\Omega_0}Q_1^{-1}\tilde z(x;v)\overline{\tilde
z(x;v)}\,dx\right)^{1/2}
$$
$$
+\left(\int_{\Gamma}Q_2^{-1}\tilde z(\cdot;u)\overline{\tilde
z(\cdot;u)}\,d\Gamma\right)^{1/2}\left(\int_{\Gamma}Q_2^{-1}\tilde z(\cdot;v)\overline{\tilde
z(\cdot;v)}\,d\Gamma\right)^{1/2}
$$
$$
+\left(\sum_{k=1}^{m} \int_{\Omega_k} r_k^{-2}(x)|u_k(x)|^2\, dx\right)^{1/2}
\left(\sum_{k=1}^{m} \int_{\Omega_k} r_k^{-2}(x)|v_k(x)|^2\, dx\right)^{1/2}
$$
$$
\leq \left(\int_{\Omega_0}Q_1^{-1}\tilde z(x;u)\overline{\tilde
z(x;u)}\,dx+\int_{\Gamma}Q_2^{-1}\tilde z(\cdot;u)\overline{\tilde
z(\cdot;u)}\,d\Gamma+\sum_{k=1}^{m} \int_{\Omega_k} r_k^{-2}(x)|u_k(x)|^2\, dx\right)^{1/2}
$$
$$
\times \left(\int_{\Omega_0}Q_1^{-1}\tilde z(x;v)\overline{\tilde
z(x;v)}\,dx+\int_{\Gamma}Q_2^{-1}\tilde z(\cdot;v)\overline{\tilde
z(\cdot;v)}\,d\Gamma+\sum_{k=1}^{m} \int_{\Omega_k} r_k^{-2}(x)|v_k(x)|^2\, dx\right)^{1/2}
$$
\begin{multline*}
\leq \Biggl\{\left(\int_{\Omega_0}|Q_1^{-1}\tilde z(x;u)|^2\,dx\right)^{1/2}
\left(\int_{\Omega_0}|\tilde z(x;u)|^2\,dx\right)^{1/2}+\left(\int_{\Gamma}|Q_2^{-1}\tilde z(\cdot;u)|^2\,d\Gamma\right)^{1/2}\Biggr.\\ \Biggl.
\times\left(\int_{\Gamma}|\tilde z(\cdot;u)|^2\,d\Gamma\right)^{1/2}
+\left(\sum_{k=1}^{m} \int_{\Omega_k} r_k^{-4}(x)|u_k(x)|^2\, dx\right)^{1/2}
\left(\sum_{k=1}^{m} \int_{\Omega_k} |u_k(x)|^2\, dx\right)^{1/2}\Biggr\}^{1/2}
\end{multline*}
\begin{multline*}
\times \Biggl\{\left(\int_{\Omega_0}|Q_1^{-1}\tilde z(x;v)|^2\,dx\right)^{1/2}
\left(\int_{\Omega_0}|\tilde z(x;v)|^2\,dx\right)^{1/2}+\left(\int_{\Gamma}|Q_2^{-1}\tilde z(\cdot;v)|^2\,d\Gamma\right)^{1/2}\Biggr.\\ \Biggl.
\times\left(\int_{\Gamma}|\tilde z(\cdot;v)|^2\,d\Gamma\right)^{1/2}
+\left(\sum_{k=1}^{m} \int_{\Omega_k} r_k^{-4}(x)|v_k(x)|^2\, dx\right)^{1/2}
\left(\sum_{k=1}^{m} \int_{\Omega_k} |v_k(x)|^2\, dx\right)^{1/2}\Biggr\}^{1/2}
\end{multline*}
$$
\leq \max\{\|Q_1^{-1}\|,\|Q_2^{-1}\|,\beta\}
\Biggl\{\int_{\Omega_0}|\tilde z(x;u)|^2\,dx+ \int_{\Gamma}|\tilde z(\cdot;u)|^2\,d\Gamma+\sum_{k=1}^{m} \int_{\Omega_k} |u_k(x)|^2\, dx\Biggr\}^{1/2}
$$
$$
\times
\Biggl\{\int_{\Omega_0}|\tilde z(x;v)|^2\,dx+ \int_{\Gamma}|\tilde z(\cdot;v)|^2\,d\Gamma+\sum_{k=1}^{m} \int_{\Omega_k} |v_k(x)|^2\, dx\Biggr\}^{1/2}
$$
$$
\leq \max\{\|Q_1^{-1}\|,\|Q_2^{-1}\|,\beta\}
\bigl\{\|\tilde z(\cdot;u)\|_{H^1(\Omega_0)}^2+ \|\gamma_0 \tilde z(\cdot;u)|_{L^2(\Gamma)}^2+\|u\|_H^2\bigr\}^{1/2}
$$
\begin{equation}\label{jlz9}
\times \bigl\{\|\tilde z(\cdot;v)\|_{H^1(\Omega_0)}^2+ \|\gamma_0 \tilde z(\cdot;v)|_{L^2(\Gamma)}^2+\|v\|_H^2\bigr\}^{1/2},
\end{equation}
where
$\beta:=\max_{1\leq k\leq m} \max_{x\in\bar\Omega_k}r_k^{-4}(x)>0.$

Setting in \eqref{g2cont1} $\psi=\tilde z(\cdot;u),$ $f=-\sum_{k=1}^m
\chi_{\Omega_k}(x)\int_{\Omega_k}\overline{ g_{k}(\eta,x)}u_k(\eta)\,d\eta,$ and $g=0,$ we find
\begin{multline}\label{g2p8}
\|\tilde z(\cdot;u)\|_{H^1(\Omega_0)}\leq \|\tilde z(\cdot;u)\|_{H^1(\Omega_R\setminus\bar\Omega)}
\\ \leq \alpha \left\|-\sum_{k=1}^m
\chi_{\Omega_k}(x)\int_{\Omega_k}\overline{ g_{k}(\eta,x)}u_k(\eta)\,d\eta\right\|
_{L^2(\Omega_R\setminus\bar\Omega)}
\leq c_1\|u\|_{H},
\end{multline}where $c_1$ is a constant independent of  $R.$
The trace theorem and inequality (\ref{g2p8}) imply
\begin{equation}\label{g2p9}
\|\gamma_0 \tilde z(\cdot;u)\|_{L^2(\Gamma)}\leq c_2\|\tilde z(\cdot;u)\|_{H^1(\Omega_R\setminus\bar\Omega)}\leq c_3\|u\|_{H},\quad c_2,c_3=\mbox{const}.
\end{equation}
From \eqref{jlz9}--\eqref{g2p8} we have
\begin{multline*}
|\pi(u,v)| \leq \max\{\|Q_1^{-1}\|,\|Q_2^{-1}\|,\beta\}(c_1\|u\|^2_{H}
+c_4\|u\|^2_{H}+\|u\|^2_{H})^{1/2}\\
\times (c_1\|v\|^2_{H} +c_4\|v\|^2_{H}+\|v\|^2_{H})^{1/2}
\leq c \|u\|_{H}  \|v\|_{H},
\end{multline*}
where $c$ is a constant independent of $R.$

Thus inequality (\ref{g2rep5}) and, consequently, the continuity
of form (\ref{g2rep2a}) are proved.

In line with Remark 1.4 to Theorem 1.1
proved in
\cite{BIBLlio}, p. 11,  the latter statements \label{g2lh}imply the existence of the unique element
 $\hat u\in H$ such that
\begin{equation*}
I(\hat u)=\inf_{u \in H} I(u).
\end{equation*}
Therefore, for any $\tau\in R$ and $v\in H$, the following relations are valid
\begin{equation}\label{g2z43d}
\frac{d}{d\tau}I(\hat u+\tau
v)\Bigl.\Bigr|_{\tau=0}=0\quad\mbox{and}\quad \frac{d}{d\tau}I(\hat
u+i\tau v)\Bigl.\Bigr|_{\tau=0}=0,
\end{equation}
where $i=\sqrt{-1}.$
Since $z(x;\hat u+\tau v)=z(x;\hat
u)+\tau\tilde z(x;v),$ where $\tilde z(x;v)$ is the unique solution to BVP
 (\ref{g20g3})$-$(\ref{g23g3}) at $u=v$ and $l_0=0,$ the first relation in (\ref{g2z43d})
 yields
$$
0=\frac 1{2}\frac{d}{d\tau}I(\hat u+\tau v)|_{\tau=0}
$$
\begin{equation*} =\lim_{\tau\to 0}\frac
1{2\tau}\Bigl\{\Bigl[(Q_1^{-1}z(\cdot;\hat u+\tau v),z(\cdot;\hat
u+\tau v))_{L^2(\Omega_0)} -(Q_1^{-1}z(\cdot;\hat
u),z(\cdot;\hat u))_{L^2(\Omega_0)}\Bigr]\Bigr.
\end{equation*}
\begin{equation*}
+\Bigl[(Q_2^{-1}z(\cdot;\hat u+\tau v),z(\cdot;\hat
u+\tau v))_{L^2(\Gamma)}\Bigr.\Bigr. \Bigl.-(Q_2^{-1}z(\cdot;\hat
u),z(\cdot;\hat u))_{L^2(\Gamma)}\Bigr]
\end{equation*}
$$
\Bigl.+\Bigl[\sum_{k=1}^{m} \bigl\{\left( r_k^{-2}(\hat u_k+\tau v),
(\hat u_k+\tau v)\right)_{L^2(\Omega_k)}-\left( r_k^{-2}\hat u_k,
\hat u_k\right)_{L^2(\Omega_k)}\bigr\}\Bigr]\Bigr\}
$$
$$
=\mbox{Re}\Bigl\{(Q_1^{-1}z(\cdot;\hat u),\tilde z(\cdot;v
))_{L^2(\Omega_0)}+(Q_2^{-1}z(\cdot;\hat u),\tilde z(\cdot;v))_{\Gamma}+\sum_{i=1}^m\left(r_k^{-2}\hat u_k,
v_k\right)_{L^2(\Omega_k)}\Bigr\}.
$$
Similarly, taking into account that $z(x;\hat u+i\tau v)=z(x;\hat
u)+i\tau\tilde z(x;v),$ we find
$$
0=\frac 1{2}\frac{d}{d\tau}I(\hat u+i\tau v)|_{\tau=0}
$$
$$
=\mbox{Im}\Bigl\{(Q_1^{-1}z(\cdot;\hat u),\tilde z(\cdot;v
))_{L^2(\Omega_0)}+(Q_2^{-1}z(\cdot;\hat u),\tilde z(\cdot;v))_{\Gamma}+\sum_{i=1}^m\left(r_k^{-2}\hat u_k,
v_k\right)_{L^2(\Omega_k)}\Bigr\};
$$
consequently,
\begin{equation}\label{g2vc}
(Q_1^{-1}z(\cdot;\hat u),\tilde z(\cdot;v
))_{L^2(\Omega_0)}+(Q_2^{-1}z(\cdot;\hat u),\tilde z(\cdot;v))_{\Gamma}+\sum_{i=1}^m\left(r_k^{-2}\hat u_k,
v_k\right)_{L^2(\Omega_k)}=0.
\end{equation}

Introduce a function $p(x)\in H^{1}
((\Omega_R\setminus \bar\Omega),\Delta)$ as the unique solution to the BVP
\begin{equation}\label{g21yxr1}
-(\Delta +k^2)p(x)=\chi_{\Omega_0}(x)Q_1^{-1}z(x;\hat u)\,\, \mbox{in}\,\,\Omega_R\setminus \bar\Omega,
\end{equation}
\begin{equation}\label{g22yxr1}
\frac{\partial p}{\partial\nu}=Q_2^{-1}z(\cdot;\hat u)\,\,\mbox{on}\,\,\Gamma,
\end{equation}
\begin{equation}\label{g23yxr1}
\frac{\partial p}{\partial\nu}=M^{(1)}_k p\,\,\mbox{on}\,\,\Gamma_R,
\end{equation}
or to an equivalent variational problem
\begin{multline}\label{g2var2}
a(p,\theta)=\int_{\Omega_R\setminus\bar\Omega}
(\nabla p\nabla \overline\theta-k^2 p\,\overline\theta)\,dx-\int_{\Gamma_R}M^{(1)}_k p\,
\overline\theta\,d\Gamma_R\\
=\int_{\Omega_R\setminus\bar\Omega}
\chi_{\Omega_0}(x)Q_1^{-1}z(x;\hat u) \overline\theta\,dx+\int_{\Gamma_R}Q_2^{-1}z(\cdot;\hat u)
\overline\theta\,d\Gamma_R \quad\forall\theta\in H^1(\Omega_R\setminus\bar\Omega).
\end{multline}
Setting in \eqref{g2var2} $\theta=\tilde z(\cdot;v),$
we obtain
\begin{multline}\label{g2nn}
a(p,\tilde z(\cdot;v))=\int_{\Omega_R\setminus\bar\Omega}
(\nabla p\nabla \overline {\tilde z(x;v)}
-k^2 p\,\overline{\tilde z(x;v)})\,dx
-\int_{\Gamma_R}M^{(1)}_k p\,
\overline{\tilde z(\cdot;v)}\,d\Gamma_R\\
=\int_{\Omega_R\setminus\bar\Omega}
\chi_{\Omega_0}(x)Q_1^{-1}z(x;\hat u) \overline{\tilde z(x;v)}\,dx+\int_{\Gamma_R}Q_2^{-1}z(\cdot;\hat u)
\overline{\tilde z(\cdot;v)}\,d\Gamma_R.
\end{multline}

Taking into account the fact that $\tilde z(\cdot;v)$
satisfies variation equation \eqref{g2varad} with $\psi=\tilde z(\cdot;v)$
equivalent to BVP \eqref{g20g3}--\eqref{g23g3} with $u=v$  and
putting in \eqref{g2varad} $\theta=p,$ we have
$$
a^*(\tilde z(\cdot;v),p)=\int_{\Omega_R\setminus\bar\Omega}
(\nabla \tilde z(x;v)\nabla\overline p
(x)-\bar k^2\tilde z(x;v)\,\overline p(x))\,dx-\int_{\Gamma_R}M_{\bar k}^{(2)}\tilde z(\cdot;v)\,\overline p\,d\Gamma_R
$$
\begin{equation}\label{g25gc}
=-\int_{\Omega_R\setminus\bar\Omega}
\sum_{k=1}^{m}
\chi_{\Omega_k}(x)\int_{\Omega_k} \overline{g_{k}(y,x)} v_k(y)\,dy\overline{p(x)}\,dx.
\end{equation}
Since $a^*(\tilde z(\cdot;v),p)=\overline{a(p,\tilde z(\cdot;v))},$
we have
$$
\int_{\Omega_0}
\overline{Q_1^{-1}z(x;\hat u)} \tilde z(x;v)\,dx+\int_{\Gamma_R}\overline{Q_2^{-1}z(\cdot;\hat u)}
\tilde z(\cdot;v)\,d\Gamma_R
$$
$$
=-\int_{\Omega_R\setminus\bar\Omega}
\sum_{k=1}^{m}
\chi_{\Omega_k}(x)\int_{\Omega_k} \overline{g_{k}(y,x)} v_k(y)\,dy\overline{p(x)}\,dx
$$
\begin{equation}\label{g2vc0}
=-\sum_{k=1}^{m}\int_{\Omega_k}
\left(\int_{\Omega_k} \overline{g_{k}(x,y)p(y)}\,dy\right) v_k(x)\,dx.
\end{equation}
From \eqref{g2vc} it follows that
\begin{equation}\label{g2vc1}
\overline{(Q_1^{-1}z(\cdot;\hat u),\tilde z(\cdot;v
))_{L^2(\Omega_0)}}+\overline{(Q_2^{-1}z(\cdot;\hat u),\tilde z(\cdot;v))_{\Gamma}}=-\sum_{k=1}^{m}
\int_{\Omega_k}r_k^{-2}(x)\overline{\hat u_k(x)} v_k(x)\,dx.
\end{equation}
Relations \eqref{g2vc0} and \eqref{g2vc1} imply
\begin{multline*}
\sum_{k=1}^{m}
\int_{\Omega_k}r_k^{-2}(x)\overline{\hat u_k(x)} v_k(x)\,dx\\=
\sum_{k=1}^{m}\int_{\Omega_k}
\left(\int_{\Omega_k} \overline{g_{k}(x,y)p(y)\,dy}\right) v_k(x)\,dx.
\forall v_k\in L^2(\Omega_k),\,\,k=\overline{1,m}.
\end{multline*}
Hence,
$$
\hat u_k(x)=r_k^{2}(x)\int_{\Omega_k}g_{k}(x,y)p(y)\,dy.
$$

Now let us establish the validity of formula \eqref{k1k4}. We have
$$
\sigma^2 := \sup_{(\tilde f,\tilde g) \in G_0,\, \tilde \xi \in G_1}
M[l(\varphi)-\widehat {\widehat {l(\varphi)}}]^2=
$$
$$
\int_{\Omega_0} Q_1^{-1}z(x)\overline{z(x)}\, dx +\int_{\Gamma} Q_2^{-1}z\overline{z}\,d\Gamma+
\sum_{k=1}^{m} \int_{\Omega_k} \hat r_k^2(x)\left|\int_{\Omega_k}g_{k}(x,y)p(y)\,dy\right|^2\, dx
$$
Transform the sum of first two terms. Make use of equalities
\eqref{g2hat2}$-$\eqref{g2hat5} to obtain
\begin{equation}\label{g2t3}
a(p,z)
=\int_{\Omega_R\setminus\bar\Omega}
\chi_{\Omega_0}(x)Q_1^{-1}z(x) \overline{ z(x)}\,dx+\int_{\Gamma}Q_2^{-1}z
\overline{z}\,d\Gamma,
\end{equation}
hence
\begin{equation}\label{g2lo}
\sigma^2=a(p,z)+
\sum_{k=1}^{m} \int_{\Omega_k} \hat r_k^2(x)\left|\int_{\Omega_k}g_{k}(x,y)p(y)\,dy\right|^2\, dx.
\end{equation}
Note that $z$ satisfies \eqref{g20gr}$-$\eqref{g23gr} which yields an integral identity
\begin{multline}
a^*(z,\theta)
=\int_{\Omega_R\setminus\bar\Omega}\left(\chi_{\omega_0}(x)l_0(x)\right.\\\left.
-\sum_{k=1}^{m}
\chi_{\Omega_k}(x)\int_{\Omega_k} \overline{g_{k}(y,x)}\hat u_k(y)\,dy\right)\overline{\theta(x)}\,dx
\quad\forall\theta\in H^1(\Omega_R\setminus\bar\Omega);
\end{multline}
setting $\theta=p$, we find
$$
a^*(z,p)
=\int_{\Omega_R\setminus\bar\Omega}\left(\chi_{\omega_0}(x)l_0(x)
-\sum_{k=1}^{m}
\chi_{\Omega_k}(x)\int_{\Omega_k} \overline{g_{k}(y,x)}\hat u_k(y)\,dy\right)\overline{p(x)}\,dx.
$$
From the latter relations, the formula
$\overline{a^*(z,p})=a(p,z)$, and \eqref{g2lo} it follows that
$$
\sigma^2=\int_{\omega_0}\overline{l_0(x)}p(x)\,dx-\sum_{k=1}^{m}
\int_{\Omega_k}\int_{\Omega_k} g_{k}(y,x)\overline{\hat u_k(y)}\,dyp(x)\,dx
$$
$$
+
\sum_{k=1}^{m} \int_{\Omega_k} \hat r_k^2(x)\left|\int_{\Omega_k}g_{k}(x,y)p(y)\,dy\right|^2\, dx.
$$
However,
$$
\hat u_k(y)=r_k^2(y)\int_{\Omega_k}g_{k}(y,\eta)p(\eta)\,d\eta,
$$
therefore,
$$
\sum_{k=1}^{m}
\int_{\Omega_k}\int_{\Omega_k} g_{k}(y,x)\overline{\hat u_k(y)}\,dyp(x)\,dx
$$
$$
=\sum_{k=1}^{m}
\int_{\Omega_k}\int_{\Omega_k} g_{k}(y,x)r_k^2(y)\overline{\int_{\Omega_k}
g_{k}(y,\eta)p(\eta)\,}d\eta dyp(x)\,dx
$$
$$
=\sum_{k=1}^{m}
\int_{\Omega_k}r_k^2(x)\int_{\Omega_k} g_{k}(x,y)p(y)\,dy\overline{\int_{\Omega_k}
g_{k}(x,\eta)p(\eta)\,}d\eta \,dx
$$
$$
=\sum_{k=1}^{m} \int_{\Omega_k} \hat r_k^2(x)\left|\int_{\Omega_k}g_{k}(x,y)p(y)\,dy\right|^2\, dx.
$$
Finally, we obtain
$$
\sigma^2=\int_{\omega_0}\overline{l_0(x)}p(x)\,dx.
$$
 \end{proof}

 In the following theorem we obtain an alternative representation
for minimax estimate $ \widehat{\widehat {l(\varphi)}}$ that does
not depend on the form of functional $l.$
 \begin{pred}
The minimax estimate of $l(\varphi)$ has the form
\begin{equation}\label{g2rt2}
\widehat{\widehat{l(\varphi)}}=l(\hat\varphi(\cdot,\omega)),
\end{equation}
where function $\hat\varphi$ is a solution to the following problem:
\begin{equation}\label{g20gr2}
\hat p\in L^2(\Sigma,H^{1}
((\Omega_R\setminus \bar\Omega),\Delta)),
\end{equation}\\[-55pt]
\begin{multline}\label{g21gr2}
-(\Delta +\bar
k^2)\hat p(x,\omega)\\=
\chi_{\Omega_k}(x)\sum_{k=1}^{m}\int_{\Omega_k}
 r_k^2(\tau)\overline{ g_k(\tau,x)} \left[y_k(\tau,\omega)-\int_{\Omega_k}
\hat \varphi(\eta,\omega) g_k(\tau,\eta)d\eta \right] d\tau
\,\,\mbox{in}\,\, \Omega_R\setminus\bar \Omega,
\end{multline}
\begin{equation}\label{g22gr2}
\frac{\partial \hat p(\cdot,\omega)
}{\partial\nu}=0\,\,\mbox{on}\,\,\Gamma,
\end{equation}
\begin{equation}\label{g23gr2}
\frac{\partial \hat p(\cdot,\omega)}{\partial
r}=M^{(2)}_{\bar k}\hat p(\cdot,\omega)\,\,\mbox{on}\,\,\Gamma_R
\end{equation}
\begin{equation}\label{g20yxr2}
\hat\varphi\in L^2(\Sigma,H^{1}
((\Omega_R\setminus \bar\Omega),\Delta)),
\end{equation}
\begin{equation}\label{g21yxr2}
-(\Delta +k^2)\hat \varphi(x,\omega)=\chi_{\Omega_0}(x)Q_1^{-1}\hat p(x,\omega)+f_0(x)\,\, \mbox{in}\,\,\Omega_R\setminus \bar\Omega,
\end{equation}
\begin{equation}\label{g22yxr2}
\frac{\partial \hat \varphi(\cdot,\omega)}{\partial\nu}=Q_2^{-1}\hat p(\cdot,\omega)+g_0,\,\,\mbox{on}\,\,\Gamma,
\end{equation}
\begin{equation}\label{g23yxr2}
\frac{\partial \hat \varphi(\cdot,\omega)}{\partial\nu}=M^{(1)}_k \hat \varphi(\cdot,\omega)\,\,\mbox{on}\,\,\Gamma_R,
\end{equation}
where equalities (\ref{g21gr2})$-$(\ref{g23gr2}) and
(\ref{g21yxr2})$-$(\ref{g23yxr2}) are fulfilled with probability
$1.$ Problem \eqref{g20gr2}$-$\eqref{g23yxr2} is uniquely
solvable. The restrictions of the solutions of this problem on
$\Omega_{R=\mbox{\rm\scriptsize min}\{R_1,R_2\}}$ corresponding to
$R=R_1$ and $R=R_2$ coincide with the solution corresponding to
$R=\mbox{\rm min}\{R_1,R_2\}.$
\end{pred}
\begin{proof}
The proof is similar to the that of Theorem 2.1.
Consider the problem of optimal control
of the equation system
\begin{equation}\label{g20gr20g}
\hat p(\cdot,\cdot;u)\in L^2(\Sigma,H^{1}
((\Omega_R\setminus \bar\Omega),\Delta)),
\end{equation}
\begin{equation}\label{g2gbbk21g}
-(\Delta +\bar
k^2)
\hat p(x,\omega;u)=d_0(x,\omega)-
\sum_{k=1}^{m}\int_{\Omega_k}\chi_{\Omega_k}(x)
\overline{g_k(\tau,x)} u(\tau,\omega) d\tau
\,\,\mbox{in}\,\, \Omega_R\setminus\bar \Omega,
\end{equation}
\begin{equation}\label{g2gbbk22g}
\frac{\partial \hat p(\cdot,\omega;u)
}{\partial\nu}=0\,\,\mbox{on}\,\,\Gamma,
\end{equation}
\begin{equation}\label{g2gbb23g}
\frac{\partial \hat p(\cdot,\omega;u)}{\partial
r}=M^{(2)}_{\bar k}\hat p(\cdot,\omega;u)\,\,\mbox{on}\,\,\Gamma_R
\end{equation}
with the cost function
$$
I(u)=\mathbb
E\left\{\int_{\Omega_0} Q_1^{-1}(\hat p(\cdot,\omega;u)
+Q_1f_0)(x)\overline{(\hat p(\cdot,\omega;u)+Q_1f_0)(x)}\, dx\right.
$$
$$
\left. +\int_{\Gamma} Q_2^{-1}(\hat p(\cdot,\omega;u)
+Q_2g_0)(x)\overline{(\hat p(\cdot,\omega;u)+Q_2g_0)(x)}\,d\Gamma\right.
$$
\begin{equation} \label{g24gst}
\left.+
\sum_{k=1}^{m} \int_{\Omega_k} r_k^{-2}(x)|u_k(x;\omega)|^2\, dx\right\}\,\to \inf_{ u\in L^2(\Sigma,H)}.
\end{equation}
where
$$
d_0(x,\omega):=\sum_{k=1}^{m}\int_{\Omega_k}\chi_{\Omega_k}(x)
  \overline{g_k(\tau,x)} r_k^2(\tau)y_k(\tau;\varphi,\xi_k(\omega))\,d\tau,
$$
The form of functional $I(u)$ and proof of Theorem 2.1 suggest that there is one
and only one element $\hat
u\in L^2(\Sigma,H)$ such that
$$
I(\hat u)= \inf_{u\in L^2(\Sigma,H)} I(u).
$$
Next, denoting by $\hat \varphi(t;\omega)$ the unique solution to
the BVP
$$
\hat\varphi(\cdot,\omega)\in L^2(\Sigma,X_{\Omega_R\setminus \bar\Omega}),
$$
$$
-(\Delta +k^2)\hat \varphi(x,\omega)=\chi_{\Omega_0}(x)Q_1^{-1}\hat p(x,\omega;\hat u)+f_0(x)\,\, \mbox{in}\,\,\Omega_R\setminus \bar\Omega,
$$
$$
\frac{\partial \hat \varphi(\cdot,\omega)}{\partial\nu}=Q_2^{-1}\hat p(\cdot,\omega;\hat u)+g_0,\,\,\mbox{on}\,\,\Gamma,
$$
$$
\frac{\partial \hat \varphi(\cdot,\omega)}{\partial\nu}=M^{(1)}_k \hat \varphi(\cdot,\omega)\,\,\mbox{on}\,\,\Gamma_R,
$$
and making use of virtually the same reasoning that led to the proof of Theorem 2.1
(by applying estimate \eqref{g2alpha1}  instead of \eqref{g2cont1}),
we arrive at the equality $\hat u(\tau,\omega)=\int_{\Omega_k}
\hat \varphi(\eta,\omega) g_k(\tau,\eta)d\eta.$ Denoting $\hat p(x,\omega)=\hat
p(x,\omega;\hat u),$ we deduce from the latter statement the unique solvability of BVP
\eqref{g20gr2}$-$\eqref{g23yxr2}.

Now let us prove the representation $\widehat{\widehat
{l(\varphi)}}=l(\hat \varphi).$
By virtue of (\ref{g25}) and (\ref{g2rt1}),
$$
\widehat{\widehat
{l(\varphi)}}=\sum_{k=1}^m\int_{\Omega_k}\overline{\hat u_k(x)}
y_k(x;\varphi,\xi_k(\omega))+\hat c
$$
$$
=\sum_{k=1}^{m}\int_{\Omega_k}
\int_{\Omega_k}r_k^2(x)\overline{ g_k(x,\eta)p(\eta)}\,d\eta
\,y_k(x;\varphi,\xi_k(\omega))dx+\hat c
$$
\begin{equation}\label{g2hat}
=\sum_{k=1}^{m}\int_{\Omega_k}
\int_{\Omega_k}r_k^2(\tau)\overline{ g_k(\tau,x)}y_k(\tau;\varphi,\xi_k(\omega))
\overline{p(x)}\,d\tau\,dx+\hat c.
\end{equation}
The function $\hat p(x,\omega;\hat u):=\hat p(x,\omega)$ is a solution to BVP (\ref{g20gr20g})$-$(\ref{g2gbb23g}) with $u=\hat u$, therefore $\forall\theta\in
H^1(\Omega_R\setminus\Omega)$ the following identity holds
\begin{multline*}
a^*(\hat p(\cdot,\omega),\theta(\cdot)
=\int_{\Omega_R\setminus\bar\Omega}(\nabla \hat p(x,\omega)\nabla\overline{\theta
(x)})-\bar k^2\hat p(x,\omega)\,\overline{\theta(x)})\,dx\\-\int_{\Gamma_R}M_{\bar k}^{(2)}\hat p(\cdot,\omega)\,\overline{\theta(\cdot)}\,d\Gamma_R
\end{multline*}
\begin{multline}\label{g25gz}
=\int_{\Omega_R\setminus\Omega}\sum_{k=1}^{m}\int_{\Omega_k}\chi_{\Omega_k}(x)
 r_k^2(\tau)\overline{ g_k(\tau,x)} \left[y_k(\tau,\omega,\xi_k(\omega))\right.\\ \left.-\int_{\Omega_k}
\hat \varphi(\eta,\omega) g_k(\tau,\eta)d\eta \right] d\tau
\overline{\theta(x}\,dx.
\end{multline}
Setting in this identity $\theta(x)=p(x)$ we obtain\label{g2mirk}
\begin{multline*}
a^*(\hat p(\cdot,\omega),p(\cdot))
=\int_{\Omega_R\setminus\bar\Omega}(\nabla \hat p(x,\omega)\nabla\overline{p
(x)})-\bar k^2\hat p(x,\omega)\,\overline{p(x)})\,dx\\-\int_{\Gamma_R}M_{\bar k}^{(2)}\hat p(\cdot,\omega)\,\overline{p(\cdot)}\,d\Gamma_R
\end{multline*}
\begin{multline}\label{g25gz1}
=\int_{\Omega_R\setminus\Omega}\sum_{k=1}^{m}
\int_{\Omega_k}\chi_{\Omega_k}(x)
 r_k^2(\tau)\overline{ g_k(\tau,x)} \left(y_k(\tau,\varphi,\xi_k(\omega))\right.\\\left.-\int_{\Omega_k}
\hat \varphi(\eta,\omega) g_k(\tau,\eta)d\eta \right) d\tau
\overline{p(x)}\,dx.
\end{multline}
Since $p(x)$ satisfies \eqref{g20yxr}$-$\eqref{g23yxr} and consequently (\ref{g2var2}),
we have
\begin{multline*}
a(p,\theta)=\int_{\Omega_R\setminus\bar\Omega}
(\nabla p(x)\nabla \overline{\theta(x)}-k^2 p(x)\,\overline{\theta(x)})\,dx-\int_{\Gamma_R}M^{(1)}_k p\,
\overline{\theta}\,d\Gamma_R\\
=\int_{\Omega_R\setminus\bar\Omega}
\chi_{\Omega_0}(x)Q_1^{-1}z(x) \overline{\theta(x)}\,dx+\int_{\Gamma}Q_2^{-1}z
\overline{\theta(\cdot)}\,d\Gamma \quad\forall\theta\in
H^1(\Omega_R\setminus\bar\Omega),
\end{multline*}
which yields
\begin{multline}\label{g2bnm}
a(p,\hat p(\cdot,\omega)=\int_{\Omega_R\setminus\bar\Omega}
(\nabla p(x)\nabla \overline{\hat p(x,\omega)}-k^2 p(x)\,\overline{\hat p(x,\omega)})\,dx-\int_{\Gamma_R}M^{(1)}_k p\,
\overline{\hat p(\cdot,\omega)}\,d\Gamma_R\\
=\int_{\Omega_R\setminus\bar\Omega}
\chi_{\Omega_0}(x)Q_1^{-1}z(x) \overline{\hat p(x,\omega)}\,dx+\int_{\Gamma}Q_2^{-1}z
\overline{\hat p(\cdot,\omega)}\,d\Gamma.
\end{multline}
Since
$a^*(\hat p(\cdot;\omega),p)=\overline{a(p,\hat p(\cdot;\omega))},$
(\ref{g25gz1}) and (\ref{g2bnm}) imply
$$
\sum_{k=1}^{m}\int_{\Omega_k}
\int_{\Omega_k}r_k^2(\tau)\overline{ g_k(\tau,x)}\left(y_k(\tau;\varphi,\xi_k(\omega))
-\int_{\Omega_k}\hat\varphi(\eta,\omega)
g_k(\tau,\eta)\,d\eta\right)d\tau\overline{p(x)}\,dx
$$
\begin{equation}\label{g2hat1}
=\int_{\Omega_R\setminus\bar\Omega}
\chi_{\Omega_0}(x)Q_1^{-1}\hat p(x,\omega) \overline{z(x)}\,dx+\int_{\Gamma}Q_2^{-1}\hat p(\cdot,\omega)
\overline{z}\,d\Gamma.
\end{equation}
Equating \eqref{g2hat} and \eqref{g2hat1}, we obtain
$$
\widehat{\widehat
{l(\varphi)}}-\hat c-\sum_{k=1}^{m}\int_{\Omega_k}
\int_{\Omega_k}r_k^2(\tau)\overline{ g_k(\tau,x)}\int_{\Omega_k}\hat\varphi(\eta,\omega)
g_k(\tau,\eta)\,d\eta\,d\tau\overline{p(x)}\,dx
$$
\begin{equation}\label{g2hat2}
=\int_{\Omega_R\setminus\bar\Omega}
\chi_{\Omega_0}(x)Q_1^{-1} \overline{z(x)}\hat p(x,\omega)\,dx+\int_{\Gamma}Q_2^{-1}
\overline{z}\hat p(\cdot,\omega)\,d\Gamma.
\end{equation}
Next, since $\hat\varphi(\cdot,\omega)$ and $z$ satisfy,
respectively, equalities \eqref{g20yxr2}$-$\eqref{g23yxr2}  and
\eqref{g20gr}$-$\eqref{g23gr}, these functions satisfy also the
identities
\begin{multline}\label{g2hat3}
a(\hat\varphi(\cdot,\omega),\theta)=\int_{\Omega_R\setminus\bar\Omega}
\chi_{\Omega_0}(x)\left(Q_1^{-1}\hat p(x,\omega)+f_0(x)\right) \overline{\theta(x)}\,dx\\+\int_{\Gamma}
\left(Q_2^{-1}\hat p(\cdot,\omega)+g_0\right) \overline{\theta}\,d\Gamma
\end{multline}
and
\begin{multline}\label{g2hat4}
a^*(z,\theta)=\int_{\Omega_R\setminus\bar\Omega}
\left(\chi_{\omega_0}(x)l_0(x)\right.\\\left.-\sum_{k=1}^{m}\chi_{\Omega_k}(x)
\int_{\Omega_k}r_k^2(\eta)\overline{ g_k(\eta,x)}\int_{\Omega_k}
g_k(\eta,\varsigma)p(\varsigma)\,d\varsigma\,d\eta\right)
\overline{\theta(x)}\,dx.
\end{multline}
Setting in \eqref{g2hat3} $\theta(x)=z(x)$ and in \eqref{g2hat4}
$\theta(x)=\hat\varphi(x,\omega)$  and taking into notice that
$\overline{a^*(z,\hat\varphi(\cdot,\omega))}=a(\hat\varphi(\cdot,\omega),z),$ we have
$$
\int_{\Omega_R\setminus\bar\Omega}
\chi_{\Omega_0}(x)Q_1^{-1}\hat p(x,\omega)\overline{z(x)}dx
+\int_{\Gamma}
Q_2^{-1}\hat p(\cdot,\omega)\overline{z}\,d\Gamma
+\int_{\Omega_0}f_0(x) \overline{z(x)}\,dx+\int_{\Gamma}
g_0 \overline{z}\,d\Gamma
$$
\begin{multline}\label{g2hat5}
=\int_{\Omega_R\setminus\bar\Omega}
\chi_{\omega_0}(x)l_0(x)\hat\varphi(x,\omega))\,dx\\
-\int_{\Omega_R\setminus\bar\Omega}
\sum_{k=1}^{m}\chi_{\Omega_k}(\eta)
\int_{\Omega_k}r_k^2(\tau)\overline{ g_k(\tau,x)}\int_{\Omega_k}
g_k(\tau,\eta)\hat\varphi(\eta,\omega))\,d\eta\,d\tau
\overline{p(x)}\,dx.
\end{multline}
Representation \eqref{g2rt2} follows now from \eqref{g2hat2} and \eqref{g2hat5} if we take into account \eqref{g2rt1}.
\end{proof}

\begin{predlllll}\label{remark}
If we define a minimax estimate $\hat\varphi(x,\omega)$ of the
unknown solution $\varphi(x)$ of BVP (\ref{g20})--(\ref{g23}) as
the estimate linear with respect to observations (\ref{g25}),
which is determined from the condition of minimum of the maximal
mean square error of the estimate taken over sets $G_0$ and $G_1$,
then it may be shown that, under certain restrictions on $G_0$ and
$G_1$, this minimax estimate of $\varphi(x)$ coincides with the
function $\hat\varphi(x,\omega)$ obtained from the solution to
problem \eqref{g20gr2}$-$\eqref{g23yxr2}.
\end{predlllll}

\subsection{Minimax estimation of the right-hand sides of equalities that enter the statement of the boundary value problem. Representations for minimax estimates and estimation errors}

The problem is to determine a minimax estimate of the value of the functional
\begin{equation}\label{g2skk8}
l(F)=\int_{\Omega_0}\overline{l_0(x)}f(x)\,dx
+\int_{\Gamma}\overline{l_1}g\,d\Gamma
\end{equation}
from observations (\ref{g25})
in the class of estimates linear with respect to observations
\begin{equation}\label{g2linestrh}
\widehat {l(F)}= \sum_{k=1}^{m}\int_{\Omega_k}
\overline{ u_k(x)}y_k(x)\, dx+c,
\end{equation}
where $u_k \in L^2(\Omega_k)$, $k= \overline{1,m}$, $c \in \mathbb
C,$ and $l_0\in L^2(\Omega_0)$ and $l_1\in L^2(\Gamma)$ are given
functions, under the assumption that $F:=(f,g)\in G_0$ and the
errors $ \xi(\cdot)=(\xi_1(\cdot), \ldots, \xi_{m}(\cdot))$ in
observations (\ref{g25}) belong to $G_1,$ where sets $G_0$ and
$G_1$ are defined by  (\ref{g28}), (\ref{g27}), and (\ref{g29}),
respectively.
\begin{predll}\label{g2oz2}
The  estimate of the form
\begin{equation}\label{g2rtrhs}
\widehat{\widehat{l(F)}}=\sum_{k=1}^m \int_{\Omega_k}
\overline{\hat u_k(x)}y_k(x)\, dx+\hat c,
\end{equation}
will be called the minimax estimate of  $l(F)$ if the element $\hat u=(u_1,\ldots, u_m)\in H$ and number $\hat c\in \mathbb C$
are determined form the condition
$$
\sup_{\tilde F\in G_0,\tilde \eta\in G_1} \mathbf E|l(\tilde F)-\widehat{l(\tilde F)}|^2\to \inf_{u\in H,\,c\in \mathbb C}.
$$
Here
\begin{equation}\label{g2ui}
\widehat {l(\tilde F)}= \sum_{k=1}^{m}\int_{\Omega_k}
\overline{ u_k(x)}\tilde y_k(x)\, dx+c,
\end{equation}
$\tilde y_{k}(x) = \int_{\Omega_k} g_{k}(x,y)\tilde \varphi(y)\,dy
+ \tilde \xi_{k}(x),$ $x \in \Omega_k,$
 $k=
\overline{1,m},$
and $\tilde\varphi$ is a solution to BVP \eqref{g20}$-$\eqref{g23} at $f=\tilde f$
and $g=\tilde g.$  The quantity
$$
\sigma := \{\sup_{(\tilde f,\tilde g) \in G_0,\, \tilde \xi \in G_1}
\mathbb E[l(\tilde F)-\widehat {\widehat {l(\tilde F)}}]^2\}^{1/2}
$$
will be called the error of the minimax estimation of expression
\eqref{g2skk8}.
\end{predll}

\begin{predl}\label{g2lemma3}
Finding the minimax estimate of
$l(F)$ is equivalent to the problem of optimal control of a system described by the  BVP
\begin{equation}\label{g20grh}
z(\cdot;u)\in H^{1}
((\Omega_R\setminus \bar\Omega),\Delta),
\end{equation}
\begin{equation}\label{g21grh}
\Delta z(x;u)+\bar
k^2z(x;u)=\sum_{k=1}^m
\chi_{\Omega_k}(x)\int_{\Omega_k}\overline{ g_{k}(\eta,x)}u_k(\eta)\,d\eta
\,\, \mbox{in}\,\,\Omega_R\setminus\bar \Omega,
\end{equation}
\begin{equation}\label{g22grh}
\frac{\partial z(\cdot;u)
}{\partial\nu}=0\,\,\mbox{on}\,\,\Gamma,
\end{equation}
\begin{equation}\label{g23grh}
\frac{\partial z(\cdot;u)}{\partial
r}=M^{(2)}_{\bar k}z(\cdot;u)\,\,\mbox{on}\,\,\Gamma_R
\end{equation}
with the quality criterion
\begin{multline} \label{g24grh}
I(u):=\int_{\Omega_0} Q_1^{-1}(l_0+z(\cdot;u))(x)\overline{(l_0(x)+z(x;u))})\, dx \\+\int_{\Gamma} (Q_2^{-1}(l_1+z(\cdot;u)))\overline{(l_1+z(\cdot;u))}\,d\Gamma+
\sum_{k=1}^{m} \int_{\Omega_k} r_k^{-2}(x)|u_k(x)|^2\, dx\,\to \inf_{ u\in H}.
\end{multline}
\end{predl}
\begin{proof}
Taking into account \eqref{g2skk8}, \eqref{g2linestrh} and \eqref{g25}, we obtain
$$
l(\tilde F)-\widehat{l(\tilde F)}
=\int_{\Omega_0}\overline{l_0(x)}\tilde f(x)\,dx
+\int_{\Gamma}\overline{l_1}g\,d\Gamma
-\sum_{k=1}^{m}\int_{\Omega_k}
\overline{u_k(x)}y_k(x;\tilde\varphi,\tilde\xi_k)\,dx-c
$$
$$
=\int_{\Omega_0}\overline{l_0(x)}\tilde f\,dx+\int_{\Gamma}\overline{l_1}g\,d\Gamma
-\sum_{k=1}^{m}\int_{\Omega_k}
\overline{u_k(x)}\int_{\Omega_k} g_{k}(x,y)\tilde \varphi(y)\,dy\,dx
$$
$$
-\sum_{k=1}^{m}\int_{\Omega_k}
\overline{u_k(x)}\tilde\xi_k(x)\,dx-c
$$
$$
=\int_{\Omega_0}\overline{l_0(x)}\tilde f(x)\,dx+\int_{\Gamma}\overline{l_1}g\,d\Gamma
-\sum_{k=1}^{m}\int_{\Omega_k}
\int_{\Omega_k} g_{k}(y,x)\overline{u_k(y)}\,dy\, \tilde\varphi(x)\,dx
$$
$$
-\sum_{k=1}^{m}\int_{\Omega_k}
\overline{u_k(x)}\tilde\xi_k(x)\,dx-c
$$
$$
=\int_{\Omega_0}
\overline{\l_0(x)}f(x)\,dx+\int_{\Gamma}\overline{l_1}g\,d\Gamma
-\sum_{k=1}^{m}\int_{\Omega_R\setminus\bar\Omega}
\overline{\left(\chi_{\Omega_k}(x)\int_{\Omega_k} \overline{g_{k}(y,x)}u_k(y)\,dy\right)}\tilde\varphi(x)\,dx
$$
\begin{equation}\label{g24ggrh}
-\sum_{k=1}^{m}\int_{\Omega_k}
\overline{u_k(x)}\tilde\xi_k(x)\,dx-c.
\end{equation}
For any fixed $u=(u_1, \ldots , u_{m})\in H$ introduce the
function $z(x;u)$ as a unique solution of problem
\eqref{g20grh}$-$\eqref{g23grh}. According to the equivalent
variational formulation of this problem it means that $z(x;u)$
satisfies the integral identity
$$
a^*(z(\cdot;u),\theta)=\int_{\Omega_R\setminus\bar\Omega}(\nabla z(x;u)\nabla\overline\theta
(x)-\bar k^2z(x;u)\,\overline\theta(x))\,dx-\int_{\Gamma_R}M_{\bar k}^{(2)}z(\cdot;u)\,\overline\theta\,d\Gamma_R
$$
\begin{equation}\label{g25grh}
=-\int_{\Omega_R\setminus\bar\Omega}
\sum_{k=1}^{m}
\chi_{\Omega_k}(x)\int_{\Omega_k} \overline{g_{k}(y,x)}u_k(y)\,dy\overline{\theta(x)}\,dx
\quad\forall\theta\in H^1(\Omega_R\setminus\bar\Omega).
\end{equation}
Set $\theta=\tilde\varphi$ in \eqref{g25g} to obtain
$$
a^*(z(\cdot;u),\tilde\varphi)=\int_{\Omega_R\setminus\bar\Omega}
(\nabla z(x;u)\nabla\overline{\tilde\varphi(x)}
-\bar k^2z(x;u)\,\overline{\tilde\varphi(x)}\,dx-\int_{\Gamma_R}M_{\bar k}^{(2)}z(\cdot;u)\,\overline{\tilde\varphi}\,d\Gamma_R
$$
\begin{equation}\label{g26grh}
=-\int_{\Omega_R\setminus\bar\Omega}
\sum_{k=1}^{m}
\chi_{\Omega_k}(x)\int_{\Omega_k} \overline{g_{k}(y,x)}u_k(y)\,dy\overline{{\tilde\varphi}}\,dx.
\end{equation}
On the other hand, since $\tilde\varphi$ is a solution of problem \eqref{g20}$-$\eqref{g23}
with $f=\tilde f$ and $g=\tilde g,$ setting $\psi=\varphi$ and $\theta=z(\cdot;u)$ in \eqref{g2sesq}, we find
$$
a(\tilde\varphi,z(\cdot;u))=\int_{\Omega_R\setminus\bar\Omega}
(\nabla \tilde\varphi(x)\nabla\overline{z(x;u)}
-k^2\tilde\varphi(x)\,\overline{z(x;u)}\,dx-\int_{\Gamma_R}M_{ k}^{(1)}\tilde\varphi\,\overline{z(\cdot;u)}\,d\Gamma_R
$$
\begin{equation}\label{g27grh}
=\int_{\Omega_0}\tilde f(x)\overline{z(x;u)}\,dx
+\int_{\Gamma}\tilde g\overline{z(\cdot;u)}\,d\Gamma.
\end{equation}
 By Lemma 1, $\overline{a^*(z,\tilde\varphi)}=a(\tilde\varphi,z).$ This
identity together with \eqref{g24ggrh}, \eqref{g26grh}, and
\eqref{g27grh} imply
$$
l(\tilde F)-\widehat{l(\tilde F)}=
\int_{\Omega_0}
\overline{\l_0(x)}\tilde f(x)\,dx+\int_{\Gamma}\overline{l_1}\tilde g\,d\Gamma
$$
$$
-\int_{\Omega_R\setminus\bar\Omega}
\overline{
\sum_{k=1}^{m}
\chi_{\Omega_k}(x)\int_{\Omega_k} \overline{g_{k}(y,x)}u_k(y)\,dy}\tilde\varphi(x)\,dx
-\sum_{k=1}^{m}\int_{\Omega_k}
\overline{u_k(x)}\tilde\xi_k(x)\,dx-c
$$
$$
=\int_{\Omega_0}
\overline{\l_0(x)}\tilde f(x)\,dx+\int_{\Gamma}\overline{l_1}\tilde g\,d\Gamma
$$
$$
+\int_{\Omega_0}\tilde f(x)\overline{z(x;u)}\,dx
+\int_{\Gamma}\tilde g\overline{z(\cdot;u)}\,d\Gamma
-\sum_{k=1}^{m}\int_{\Omega_k}
\overline{u_k(x)}\tilde\xi_k(x)\,dx-c.
$$
\begin{equation*}\label{g28grh}
=\int_{\Omega_0}\tilde f(x)\overline{(l_0(x)+z(x;u))}\,dx
+\int_{\Gamma}\tilde g\overline{(l_1+z(\cdot;u))}\,d\Gamma
-\sum_{k=1}^{m}\int_{\Omega_k}
\overline{u_k(x)}\tilde\xi_k(x)\,dx-c.
\end{equation*}
The latter yields
$$
\mathbf E\left|l(\tilde
F)-\widehat {l(\tilde F)}\right|^2
=\left|\int_{\Omega_0}\tilde f(x)\overline{(l_0(x)+z(x;u))}\,dx
+\int_{\Gamma}\tilde g\overline{(l_1+z(\cdot;u))} \,d\Gamma \right|^2
$$
$$
+\mathbf
E\left|\sum_{k=1}^{m}\int_{\Omega_k}
\overline{u_k(x)}\,\tilde\xi_k(x)\,dx\right|^2.
$$
Therefore,
$$ \inf_{c \in
\mathbb C}\sup_{(\tilde f,\tilde g)\in G_0, \tilde \xi\in
G_1} \mathbf E|l(\tilde
F)-\widehat {l(\tilde F)}|^2=
$$
$$
=\inf_{c \in \mathbb C}\sup_{(\tilde f,\tilde g)\in G_0}\left|\int_{\Omega_0}\tilde f(x)\overline{(l_0(x)+z(x;u))}\,dx
+\int_{\Gamma}\tilde g\overline{(l_1+z(\cdot;u))}\,d\Gamma
-c\right|^2
$$
\begin{equation*}\label{g2exhrh}
+ \sup_{ \tilde \xi\in
G_1}\mathbf E\left|\sum_{k=1}^{m}\int_{\Omega_k}
\overline{u_k(x)}\,\tilde\xi_k(x)\,dx\right|^2.
\end{equation*}
Beginning from this place, we apply the same reasoning as in the proof of Lemma 1.2 (replacing $z(x,u)$ by $l_0(x)+z(x,u)$) to obtain
$$ \inf_{c \in
\mathbb C}\sup_{(\tilde f,\tilde g)\in G_0, \tilde \xi\in
G_1} \mathbf E|l(\tilde
F)-\widehat {l(\tilde F)}|^2=I(u),
$$
where $I(u)$ is determined by formula \eqref{g24grh} for
$$
c=\int_{\Omega_0}f_0(x)\overline{(l_0(x)+z(x;u))}\,dx
+\int_{\Gamma}g_0\overline{(l_1+z(\cdot;u))}\,d\Gamma .
$$
\end{proof}

The following result follows from this lemma,
\begin{pred}
The minimax estimate of $l(F)$ has the form
\begin{equation}\label{g2rtrh}
\widehat{\widehat{l(F)}}=\sum_{k=1}^m \int_{\Omega_k}
\overline{\hat u_k(x)}y_k(x)\, dx+\hat c,
\end{equation}
where
\begin{equation}\label{g2rt1rh}
\hat u_k(x)=r_k^2(x)\int_{\Omega_k}g_{k}(x,y)p(y)\,dy,
\,\,k=\overline{1,m},
\end{equation}
\begin{equation}\label{g2rt1rh1}
\hat
c=\int_{\Omega_0}\overline{(z(x)+l_0(x))}
f_0(x)\,dx+\int_\Gamma\overline{(z+l_1)}
g_0\,d\Gamma,
\end{equation}
and functions $z$ and $p$ are determined from the solution to the following problem:
\begin{equation}\label{g20grrh}
z\in H^{1}
((\Omega_R\setminus \bar\Omega),\Delta),
\end{equation}
\begin{equation}\label{g21grrh}
\Delta z(x) +\bar
k^2z(x)=\sum_{k=1}^m
\chi_{\Omega_k}(x)\int_{\Omega_k}\overline{ g_{k}(\eta,x)}\hat u_k(\eta)\,d\eta
\,\, \mbox{in}\,\,\Omega_R\setminus\bar \Omega,
\end{equation}
\begin{equation}\label{g22grrh}
\frac{\partial z
}{\partial\nu}=0\,\,\mbox{on}\,\,\Gamma,
\end{equation}
\begin{equation}\label{g23grrh}
\frac{\partial z}{\partial
r}=M^{(2)}_{\bar k}z\,\,\mbox{on}\,\,\Gamma_R,
\end{equation}
\begin{equation}\label{g20yxrrh}
p\in H^{1}
((\Omega_R\setminus \bar\Omega),\Delta),
\end{equation}
\begin{equation}\label{g21yxrrh}
-(\Delta +k^2)p=\chi_{\Omega_0}Q_1^{-1}(z+l_0)\,\, \mbox{in}\,\,\Omega_R\setminus \bar\Omega,
\end{equation}
\begin{equation}\label{g22yxrrh}
\frac{\partial p}{\partial\nu}=Q_2^{-1}(z+l_1),\,\,\mbox{on}\,\,\Gamma,
\end{equation}
\begin{equation}\label{g23yxrrh}
\frac{\partial p}{\partial\nu}=M^{(1)}_k p\,\,\mbox{on}\,\,\Gamma_R.
\end{equation}
Problem \eqref{g20grrh}$-$\eqref{g23yxrrh} is uniquely solvable.
The restrictions of the solutions of this problem on
$\Omega_{R=\mbox{\rm min}\{R_1,R_2\}}$ corresponding to
$R=R_1$ and $R=R_2$  coincide with the solution corresponding to
$R=\mbox{\rm min}\{R_1,R_2\}.$

Estimation error $\sigma$ is determined by the formula
$\sigma=l(P)^{1/2},$
where $P=(Q_1^{-1}(l_0+z),Q_2^{-1}(l_1+\gamma_Dz)).$
\end{pred}
\begin{proof}
Similarly to the proof of Theorem 2.1, we will show that the
solution to the optimal control problem
(\ref{g20grh})$-$(\ref{g24grh}) can be reduced to the solution of
the equation system \eqref{g20grrh}$-$\eqref{g23yxrrh}. To this
end, note that there exists the unique element $\hat u \in H$ at
which functional (\ref{g24grh}) attains its minimum, namely,
$I(\hat u)=\inf_{u \in H}I(u).$ In order to prove this statement,
represent $I(u)$  as
\begin{equation}\label{g2frh}
I(u)=\tilde I(u)+L(u)+\tilde c,
\end{equation}
where
$$
\tilde I(u)=\int_{\Omega_0} Q_1^{-1}z(x;u)\overline{z(x;u)}\, dx +\int_{\Gamma} Q_2^{-1}z(\cdot;u)\overline{z(\cdot;u)}\,d\Gamma+
\sum_{k=1}^{m} \int_{\Omega_k} r_k^{-2}(x)|u_k(x)|^2\, dx
$$
is a quadratic functional in space $H$ corresponding to a
semi-bilinear continuous Hermitian form \footnote{Its continuity
is shown in the course of the proof of Theorem 1.1.}
\label{g2cont222}
$$
\pi(u,v)=\int_{\Omega_0} Q_1^{-1}\tilde z(x;u)\overline{\tilde z(x;v)}\, dx +\int_{\Gamma} Q_2^{-1}\tilde z(\cdot;u)\overline{\tilde z(\cdot;v)}\,d\Gamma
$$
\begin{equation}\label{g2rep2}
+\sum_{k=1}^{m} \int_{\Omega_k} r_k^{-2}(x)u_k(x)\overline{v_k(x)}\, dx
\end{equation}
on $H\times H$ which satisfies the inequality
\begin{equation}\label{g2rep3333}
\tilde I(u)\geq a \|u\|^2_{H_0}\,\,\forall u\in
H,\,\,a=\mbox{const},
\end{equation}
$$
L(u):=2\mbox{Re}\int_{\Omega_0} Q_1^{-1}\tilde z(x;u)\overline{l_0(x)}\, dx+2\mbox{Re}\int_{\Gamma} Q_2^{-1}\tilde z(\cdot;u)\overline{l_1}\, d\Gamma,
$$
is a linear continuous functional in $H$, and
$$
\tilde c=\int_{\Omega_0} Q_1^{-1}l_0(x)\overline{l_0(x)}\, dx +
\int_{\Gamma} Q_2^{-1}l_1\overline{l_1}\,d\Gamma.
$$
This statement yields
(see page \pageref{g2lh}) the existence of the unique element $\hat u\in H$ such that
\begin{equation*}
I(\hat u)=\inf_{u \in H} I(u).
\end{equation*}
Therefore, for any $\tau\in R$ and $v\in H$ the relations
\begin{equation}\label{g2z43dd}
\frac{d}{d\tau}I(\hat u+\tau
v)\Bigl.\Bigr|_{\tau=0}=0\quad\mbox{and}\quad \frac{d}{d\tau}I(\hat
u+i\tau v)\Bigl.\Bigr|_{\tau=0}=0
\end{equation}
hold. Taking into account that functions $z(x;\hat u+\tau v)$ and
$z(x;\hat u+i\tau v)$ can be written, respectively, as $z(x;\hat
u+\tau v)=z(x;\hat u)+\tau\tilde z(x;v)$ and $z(x;\hat u+i\tau
v)=z(x;\hat u)+i\tau\tilde z(x;v)$, where  $\tilde z(x;v)$ is the
unique solution to problem (\ref{g20g3})$-$(\ref{g23g3}) at $u=v,$
we deduce from \eqref{g2z43dd} that
\begin{multline}\label{g2vcnn}
(Q_1^{-1}(l_0+z(\cdot;\hat u)),\tilde z(\cdot;v
))_{L^2(\Omega_0)}+(Q_2^{-1}(l_1+z(\cdot;\hat u)),\tilde z(\cdot;v))_{\Gamma}\\+\sum_{k=1}^{m}\left(r_k^{-2}\hat u_k,
v_k\right)_{L^2(\Omega_k)}=0.
\end{multline}

Introduce the function $p(x)\in H^{1}
((\Omega_R\setminus \bar\Omega),\Delta)$ as the unique solution to the BVP
\begin{equation}\label{g21yxr1rh}
-(\Delta +k^2)p(x)=\chi_{\Omega_0}(x)Q_1^{-1}(z(\cdot;\hat u)+l_0)(x)\,\, \mbox{in}\,\,\Omega_R\setminus \bar\Omega,
\end{equation}
\begin{equation}\label{g22yxr1rh}
\frac{\partial p}{\partial\nu}=Q_2^{-1}(z(\cdot;\hat u)+l_1)\,\,\mbox{on}\,\,\Gamma,
\end{equation}
\begin{equation}\label{g23yxr1rh}
\frac{\partial p}{\partial\nu}=M^{(1)}_k p\,\,\mbox{on}\,\,\Gamma_R,
\end{equation}
or to the equivalent variational problem
\begin{multline}\label{g2var2rh}
a(p,\theta)=\int_{\Omega_R\setminus\bar\Omega}
(\nabla p\nabla \overline\theta-k^2 p\,\overline\theta)\,dx-\int_{\Gamma_R}M^{(1)}_k p\,
\overline\theta\,d\Gamma_R\\
=\int_{\Omega_R\setminus\bar\Omega}
\chi_{\Omega_0}(x)Q_1^{-1}(z(\cdot;\hat u)+l_0)(x) \overline\theta\,dx+\int_{\Gamma_R}Q_2^{-1}(z(\cdot;\hat u)+l_1)
\overline\theta\,d\Gamma_R \quad\forall\theta\in H^1(\Omega_R\setminus\bar\Omega).
\end{multline}
Setting in \eqref{g2var2rh} $\theta=\tilde z(\cdot;v),$
we obtain
\begin{multline}\label{g2n}
a(p,\tilde z(\cdot;v))=\int_{\Omega_R\setminus\bar\Omega}
(\nabla p\nabla \overline {\tilde z(x;v)}
-k^2 p\,\overline{\tilde z(x;v)})\,dx
-\int_{\Gamma_R}M^{(1)}_k p\,
\overline{\tilde z(\cdot;v)}\,d\Gamma_R\\
=\int_{\Omega_R\setminus\bar\Omega}
\chi_{\Omega_0}(x)Q_1^{-1}(z(\cdot;\hat u)+l_0)(x) \overline{\tilde z(x;v)}\,dx+\int_{\Gamma_R}Q_2^{-1}(z(\cdot;\hat u)+l_1)
\overline{\tilde z(\cdot;v)}\,d\Gamma_R.
\end{multline}

Taking into account the fact that $\tilde z(\cdot;v)$
satisfies variation equation \eqref{g2varad} with $\psi=\tilde z(\cdot:v)$
and $f(x)=-
\sum_{k=1}^{m}
\chi_{\Omega_k}(x)\int_{\Omega_k} \overline{g_{k}(y,x)} v_k(y)$ which is equivalent to BVP \eqref{g20grh}$-$\eqref{g23grh} with $u=v$  and
setting in \eqref{g2varad} $\theta=p,$ we have
$$
a^*(\tilde z(\cdot;v),p)=\int_{\Omega_R\setminus\bar\Omega}
(\nabla \tilde z(x;v)\nabla\overline p
(x)-\bar k^2\tilde z(x;v)\,\overline p(x))\,dx-\int_{\Gamma_R}M_{\bar k}^{(2)}\tilde z(\cdot;v)\,\overline p\,d\Gamma_R
$$
\begin{equation}\label{g25gcrh}
=-\int_{\Omega_R\setminus\bar\Omega}
\sum_{k=1}^{m}
\chi_{\Omega_k}(x)\int_{\Omega_k} \overline{g_{k}(y,x)} v_k(y)\,dy\overline{p(x)}\,dx.
\end{equation}
Since $a^*(\tilde z(\cdot;v),p)=\overline{a(p,\tilde
z(\cdot;v))},$ we obtain
$$
\int_{\Omega_0}
\overline{Q_1^{-1}(l_0+z(\cdot;\hat u))(x)} \tilde z(x;v)\,dx+\int_{\Gamma_R}\overline{Q_2^{-1}(l_1+z(\cdot;\hat u))}
\tilde z(\cdot;v)\,d\Gamma_R
$$
$$=-\int_{\Omega_R\setminus\bar\Omega}
\sum_{k=1}^{m}
\chi_{\Omega_k}(x)\int_{\Omega_k} \overline{g_{k}(y,x)} v_k(y)\,dy\overline{p(x)}\,dx
$$
\begin{equation}\label{g2vc0rh}
=-\sum_{k=1}^{m}\int_{\Omega_k}
\left(\int_{\Omega_k} \overline{g_{k}(x,y)p(y)}\,dy\right) v_k(x)\,dx.
\end{equation}
From \eqref{g2vcnn} it follows that
$$
\overline{(Q_1^{-1}(l_0+z(\cdot;\hat u)),l_0+\tilde z(\cdot;v
))_{L^2(\Omega_0)}}+\overline{(Q_2^{-1}(l_1+z(\cdot;\hat u)),\tilde z(\cdot;v))_{\Gamma}}
$$
\begin{equation}\label{g2vc1rh}
=-\sum_{k=1}^{m}
\int_{\Omega_k}r_k^{-2}(x)\overline{\hat u_k(x)} v_k(x)\,dx.
\end{equation}
Relations \eqref{g2vc0rh} and \eqref{g2vc1rh} imply
\begin{multline*}
\sum_{k=1}^{m}
\int_{\Omega_k}r_k^{-2}(x)\overline{\hat u_k(x)} v_k(x)\,dx\\=
\sum_{k=1}^{m}\int_{\Omega_k}
\left(\int_{\Omega_k} \overline{g_{k}(x,y)p(y)\,dy}\right) v_k(x)\,dx.
\forall v_k\in L^2(\Omega_k),\,\,k=\overline{1,m}.
\end{multline*}
Hence,
$$
\hat u_k(x)=r_k^{2}(x)\int_{\Omega_k}g_{k}(x,y)p(y)\,dy.
$$

Now we can determine estimation error $\sigma$. Substituting
$$\hat u_k(x)=r_k^2(x)\int_{\Omega_k}g_{k}(x,y)p(y)\,dy,
\,\,k=\overline{1,m},$$ to the formula $I(\hat u)=\sigma^2,$ in
which $I(u)$ is calculated according to \eqref{g24grh} we obtain,
taking into notice that $z(t)=z(t;\hat u),$
$$
\sigma^2=
\int_{\Omega_0} Q_1^{-1}(l_0+z)(x)\overline{(l_0(x)+z(x))}\, dx +\int_{\Gamma} Q_2^{-1}(l_1+z)\overline{(l_1+z)}\,d\Gamma
$$
$$
+\sum_{k=1}^{m} \int_{\Omega_k} \hat r_k^2(x)\left|\int_{\Omega_k}g_{k}(x,y)p(y)\,dy\right|^2\, dx.
$$
Transform the sum of the first two terms. To do this, make use of equalities
\eqref{g20yxrrh}$-$\eqref{g23yxrrh} yielding
\begin{equation}\label{g2t3rh}
a(p,z)
=\int_{\Omega_R\setminus\bar\Omega}
\chi_{\Omega_0}(x)Q_1^{-1}(l_0+z)(x) \overline{ z(x)}\,dx+\int_{\Gamma}Q_2^{-1}(l_1+z)
\overline{z}\,d\Gamma,
\end{equation}
so that
$$
\sigma^2=a(p,z)+\int_{\Omega_R\setminus\bar\Omega}
\chi_{\Omega_0}(x)Q_1^{-1}(l_0+z)(x) \overline{l_0(x) }\,dx+\int_{\Gamma}Q_2^{-1}(l_1+z)
\overline{l_1}\,d\Gamma
$$
\begin{equation}\label{g2lorh}
+\sum_{k=1}^{m} \int_{\Omega_k} \hat r_k^2(x)\left|\int_{\Omega_k}g_{k}(x,y)p(y)\,dy\right|^2\, dx.
\end{equation}
However, because $z$ satisfies \eqref{g20grrh}$-$\eqref{g23grrh}, the following integral identity holds
$$
a^*(z,\theta)
=-\int_{\Omega_R\setminus\bar\Omega}\sum_{k=1}^{m}
\chi_{\Omega_k}(x)\int_{\Omega_k} \overline{g_{k}(y,x)}\hat u_k(y)\,dy\overline{\theta(x)}\,dx
\quad\forall\theta\in H^1(\Omega_R\setminus\bar\Omega);
$$
setting in this identity $\theta=p$, we have
$$
a^*(z,p)
=-\int_{\Omega_R\setminus\bar\Omega}\sum_{k=1}^{m}
\chi_{\Omega_k}(x)\int_{\Omega_k} \overline{g_{k}(y,x)}\hat u_k(y)\,dy\overline{p(x)}\,dx.
$$
The last formula, the relation
$
\overline{a^*(z,p})=a(p,z)$, and \eqref{g2lorh} give us
$$
\sigma^2=\int_{\Omega_0}\overline{l_0(x)}Q_1^{-1}(l_0+z)(x)\,dx
\int_{\Gamma}\overline{l_1}Q_2^{-1}(l_1+z)\,d\Gamma
$$
$$
-\sum_{k=1}^{m}
\int_{\Omega_k}\int_{\Omega_k} g_{k}(y,x)\overline{\hat u_k(y)}\,dyp(x)\,dx
+\sum_{k=1}^{m} \int_{\Omega_k} \hat r_k^2(x)\left|\int_{\Omega_k}g_{k}(x,y)p(y)\,dy\right|^2\, dx.
$$

Repeating literally the end of the proof of Theorem 2.3 we obtain
$$
\sum_{k=1}^{m}
\int_{\Omega_k}\int_{\Omega_k} g_{k}(y,x)\overline{\hat u_k(y)}\,dyp(x)\,dx
$$
$$
=\sum_{k=1}^{m} \int_{\Omega_k} \hat r_k^2(x)\left|\int_{\Omega_k}g_{k}(x,y)p(y)\,dy\right|^2\, dx,
$$
and finally
$$
\sigma^2=\int_{\Omega_0}\overline{l_0(x)}Q_1^{-1}(l_0+z)(x)\,dx
\int_{\Gamma}\overline{l_1}Q_2^{-1}(l_1+z)\,d\Gamma
=l(P).
$$
 \end{proof}

In Theorem \ref{g2th2} stated below we obtain another
representation for minimax estimate $ \widehat{\widehat {l(F)}}$,
not depending on the form of functional $l.$
\begin{pred}\label{g2th2}
The minimax estimate of $l(F)$ has the form
\begin{equation}\label{g2reprk}
\widehat{\widehat {l(F)}}=l(\hat F),
\end{equation}
where $\hat F=(\hat f,\hat g)$, $\hat f(x)=\hat f(x,\omega)=Q_1^{-1}\hat p(x,\omega)+f_0(x)$,
$\hat g=\hat g(\cdot,\omega)=Q_2^{-1}\gamma_D\hat p(\cdot,\omega)+g_0$, and function $\hat p=\hat p(\cdot,\omega)$ is determined from the solution to problem
 \eqref{g20gr2}$-$\eqref{g23yxr2}.
\end{pred}
\begin{proof}
Taking into notice \eqref{g2rtrh}$-$\eqref{g2rt1rh1}, we find
$$
\widehat{\widehat{l(F)}}=\sum_{k=1}^m \int_{\Omega_k}
\overline{\hat u_k(x)}y_k(x)\, dx+\hat c
$$
\begin{equation}\label{g2hatrh}
=\sum_{k=1}^{m}\int_{\Omega_k}
\int_{\Omega_k}r_k^2(\tau)\overline{ g_k(\tau,x)}y_k(\tau)
\overline{p(x)}\,d\tau\,dx+\hat c,
\end{equation}
where $\hat c$ is determined from \eqref{g2rt1rh1}.

Next, repeating literally the proof of Theorem 1.2 on page \pageref{g2mirk}, we arrive at the relatioship
\begin{multline*}
a^*(\hat p(\cdot,\omega),p(\cdot))
=\int_{\Omega_R\setminus\bar\Omega}(\nabla \hat p(x,\omega)\nabla\overline{p
(x)})-\bar k^2\hat p(x,\omega)\,\overline{p(x)})\,dx\\-\int_{\Gamma_R}M_{\bar k}^{(2)}\hat p(\cdot,\omega)\,\overline{p(\cdot)}\,d\Gamma_R
\end{multline*}
\begin{multline}\label{g25gz1rh}
=\int_{\Omega_R\setminus\Omega}\sum_{k=1}^{m}
\int_{\Omega_k}\chi_{\Omega_k}(x)
 r_k^2(\tau)\overline{ g_k(\tau,x)} \Bigl(y_k(\tau)\Bigr.\\\Bigl.-\int_{\Omega_k}
\hat \varphi(\eta,\omega) g_k(\tau,\eta)d\eta \Bigr) d\tau
\overline{p(x)}\,dx.
\end{multline}
Taking into account that $p(x)$ satisfies \eqref{g20yxrrh}$-$\eqref{g23yxrrh}, we obtain
\begin{multline*}
a(p,\theta)=\int_{\Omega_R\setminus\bar\Omega}
(\nabla p(x)\nabla \overline{\theta(x)}-k^2 p(x)\,\overline{\theta(x)})\,dx-\int_{\Gamma_R}M^{(1)}_k p\,
\overline{\theta}\,d\Gamma_R\\
=\int_{\Omega_R\setminus\bar\Omega}
\chi_{\Omega_0}(x)Q_1^{-1}(z+l_0)(x) \overline{\theta(x)}\,dx+\int_{\Gamma}Q_2^{-1}(z+l_1)
\overline{\theta(\cdot)}\,d\Gamma \quad\forall\theta\in
H^1(\Omega_R\setminus\bar\Omega);
\end{multline*}
consequently,
\begin{multline}\label{g2bnmrh}
a(p,\hat p(\cdot,\omega)=\int_{\Omega_R\setminus\bar\Omega}
(\nabla p(x)\nabla \overline{\hat p(x,\omega)}-k^2 p(x)\,\overline{\hat p(x,\omega)})\,dx-\int_{\Gamma_R}M^{(1)}_k p\,
\overline{\hat p(\cdot,\omega)}\,d\Gamma_R\\
=\int_{\Omega_R\setminus\bar\Omega}
\chi_{\Omega_0}(x)Q_1^{-1}(z+l_0)(x) \overline{\hat p(x,\omega)}\,dx+\int_{\Gamma}Q_2^{-1}(z+l_1)
\overline{\hat p(\cdot,\omega)}\,d\Gamma.
\end{multline}
From \eqref{g25gz1rh} and the latter equality, it follows that
$$
\sum_{k=1}^{m}\int_{\Omega_k}
\int_{\Omega_k}r_k^2(\tau)\overline{ g_k(\tau,x)}\left(y_k(\tau)
-\int_{\Omega_k}\hat\varphi(\eta,\omega)
g_k(\tau,\eta)\,d\eta\right)d\tau\overline{p(x)}\,dx
$$
\begin{equation}\label{g2hat1rh}
=\int_{\Omega_R\setminus\bar\Omega}
\chi_{\Omega_0}(x)Q_1^{-1}\hat p(x,\omega) \overline{(z(x)+l_0(x))}\,dx+\int_{\Gamma}Q_2^{-1}\hat p(\cdot,\omega)
\overline{(z+l_1)}\,d\Gamma.
\end{equation}
Equating \eqref{g2hatrh} and \eqref{g2hat1rh}, we find
$$
\widehat{\widehat{l(F)}}-\hat c-\sum_{k=1}^{m}\int_{\Omega_k}
\int_{\Omega_k}r_k^2(\tau)\overline{ g_k(\tau,x)}\int_{\Omega_k}\hat\varphi(\eta,\omega)
g_k(\tau,\eta)\,d\eta\,d\tau\overline{p(x)}\,dx
$$
$$
+\int_{\Omega_R\setminus\bar\Omega}
\chi_{\Omega_0}(x)Q_1^{-1}\overline{l_0}(x)\hat p(x,\omega) \,dx+\int_{\Gamma}Q_2^{-1}\overline{l_1}\hat p(\cdot,\omega)
\,d\Gamma
$$
\begin{equation}\label{g2100rh}
=\int_{\Omega_R\setminus\bar\Omega}
\chi_{\Omega_0}(x)Q_1^{-1}\overline{z}(x)\hat p(x,\omega) \,dx+\int_{\Gamma}Q_2^{-1}\overline{z}\hat p(\cdot,\omega)
\,d\Gamma.
\end{equation}

Next, since $\hat\varphi(\cdot,\omega)$ and $z$  satisfy, respectively, equalities \eqref{g20yxr2}$-$\eqref{g23yxr2} and \eqref{g20grrh}$-$\eqref{g23grrh},
these functions also satisfy the identities
\begin{equation}\label{g2hat3rh}
a(\hat\varphi(\cdot,\omega),\theta)=\int_{\Omega_R\setminus\bar\Omega}
\chi_{\Omega_0}(x)\left(Q_1^{-1}\hat p(x,\omega)+f_0(x)\right) \overline{\theta(x)}\,dx+\int_{\Gamma}
\left(Q_2^{-1}\hat p(\cdot,\omega)+g_0\right) \overline{\theta}\,d\Gamma
\end{equation}
and
\begin{equation}\label{g2hat4rh}
a^*(z,\theta)=-\int_{\Omega_R\setminus\bar\Omega}
\sum_{k=1}^{m}\chi_{\Omega_k}(x)
\int_{\Omega_k}r_k^2(\eta)\overline{ g_k(\eta,x)}\int_{\Omega_k}
g_k(\eta,\varsigma)p(\varsigma)\,d\varsigma\,d\eta
\overline{\theta(x)}\,dx.
\end{equation}
Setting in \eqref{g2hat3rh} $\theta(x)=z(x)$ and in \eqref{g2hat4rh}
$\theta(x)=\hat\varphi(x,\omega)$  and taking into account that
$\overline{a^*(z,\hat\varphi(\cdot,\omega))}=a(\hat\varphi(\cdot,\omega),z),$ we obtain
$$
\int_{\Omega_R\setminus\bar\Omega}
\chi_{\Omega_0}(x)Q_1^{-1}\hat p(x,\omega)\overline{z(x)}dx
+\int_{\Gamma}
Q_2^{-1}\hat p(\cdot,\omega)\overline{z}\,d\Gamma
+\int_{\Omega_0}f_0(x) \overline{z(x)}\,dx+\int_{\Gamma}
g_0 \overline{z}\,d\Gamma
$$
\begin{equation}\label{g2hat5rh}
=
-\int_{\Omega_R\setminus\bar\Omega}
\sum_{k=1}^{m}\chi_{\Omega_k}(\eta)
\int_{\Omega_k}r_k^2(\tau)\overline{ g_k(\tau,x)}\int_{\Omega_k}
g_k(\tau,\eta)\hat\varphi(\eta,\omega))\,d\eta\,d\tau
\overline{p(x)}\,dx.
\end{equation}

From \eqref{g2100rh} and \eqref{g2hat5rh}, it follows, in view of \eqref{g2rt1rh1}, that
$$
\widehat{\widehat{l(F)}}-\int_{\Omega_0}\overline{z(x)}
f_0(x)\,dx-\int_\Gamma\overline{z}
g_0\,d\Gamma-\int_{\Omega_0}\overline{l_0(x)}
f_0(x)\,dx-\int_\Gamma\overline{l_1}
g_0\,d\Gamma
$$
$$
+\int_{\Omega_0}\overline{z(x)}
f_0(x)\,dx+\int_\Gamma\overline{z}g_0\,d\Gamma
$$
$$
=\int_{\Omega_0}
\overline{l_0(x)}Q_1^{-1}\hat p(x,\omega)\,dx
+\int_{\Gamma}
\overline{l_1}Q_2^{-1}\hat p(\cdot,\omega)\,d\Gamma,
$$
thus
$$
 \widehat{\widehat{l(F)}}
 =\int_{\Omega_0}\overline{l_0(x)}(f_0(x)+Q_1^{-1}\hat p(x,\omega))\,dx+\int_{\Gamma}
\overline{l_1}(l_1+Q_2^{-1}\hat p(\cdot,\omega))\,d\Gamma
$$
$$
=\int_{\Omega_0}\overline{l_0(x)}\hat f(x)\,dx+\int_{\Gamma}
\overline{l_1}\hat g\,d\Gamma=l(\hat F).
$$
\end{proof}

\begin{predlllll} \label{remark2}
If we define a minimax estimate $\hat F(x,\omega)$ of the
element\footnote{Here $f$ and $g$ are the functions  entering the
statement of BVP (\ref{g20})--(\ref{g23}) and $F=(f,g)\in G_0$.}
$F=(f,g)$ as an estimate linear with respect to observations
(\ref{g25}), which is determined from the condition of minimum of
the maximal mean square error of the estimate taken over sets
$G_0$ and $G_1,$ then it may be established that, under certain
restrictions on $G_0$ and $G_1$, this minimax estimate of $F$
coincides with the element $\hat F=(\hat f,\hat g),$ where $\hat
f=Q_1^{-1}\hat p(x,\omega)+f_0(x)$ and $\hat
g=Q_2^{-1}\gamma_D\hat p(\cdot,\omega)+g_0$, and the function
$\hat p=\hat p(\cdot,\omega)$ is determined from the solution to
problem \eqref{g20gr2}$-$\eqref{g23yxr2}.
\end{predlllll}


Using Theorems 1.1$-$1.4 together with the solution techniques
employing the so-called DtN finite-element methods elaborated for
problems \eqref{g21yx}$-$\eqref{g23yx} and
\eqref{g21yxl}$-$\eqref{g23yxl}, one can construct algorithms of
numerical solution to problems \eqref{g20gr}$-$\eqref{g23yxr},
\eqref{g20gr2}$-$\eqref{g23yxr2}, and
\eqref{g20grrh}$-$\eqref{g23yxrrh} and obtain the required minimax
estimates.

\begin{predlllll}
All results of this section remain valid in the three-dimensional case.
For example, finding minimax estimates of the solutions to the BVPs for the Helmholtz equation that describe  diffraction of acoustic waves
by an obstacle $\Omega\in \mathbb R^3$ can be reduced to the solution of integro-differential equation systems \eqref{g20gr}$-$\eqref{g23yxr} and
\eqref{g20gr}$-$\eqref{g23yxr}; the domain $\mathbb R^2\setminus\bar\Omega,$ should be replaced by  $\mathbb R^3\setminus\bar\Omega,$  and plane domains $\Omega_i,$ $i=1,k,$ on which observations are made should be considered, as well as supports $\Omega_0$ and $\omega_0$ of functions $f$ and $l_0,$ as spatial domains. The Sommerfeld condition
\begin{equation}\label{*}
\frac{\partial\varphi}{\partial
r}-ik\varphi=o(1/r^{1/2}),\,\,r=|x|=\sqrt{x_1^2+x_2^2},\,\,r\to \infty
\end{equation}
should be replaced by
\begin{equation}\label{**}
\frac{\partial\varphi}{\partial
r}-ik\varphi=o(1/r),\,\,r=|x|=\sqrt{x_1^2+x_2^2+x_3^2},\,\,r\to \infty;
\end{equation}
 operators $M^{(j)}_k\psi$ defined on a circle $\Gamma_R$ should be replaced by the following operators defined on  a sphere $\Gamma_R$ of radius $R$
\begin{equation}\label{g2M22}
(M^{(j)}_k\psi)(R,\theta,\phi):=\frac k{4\pi}\sum_{n=0}^\infty\frac {h_{n}^{(j)\prime}(kR)}
{h_{n}^{(j)}(kR)}\sum_{n=-m}^mu_{mn}Y_{mn}(\theta,\phi),
\end{equation}
where $(R,\theta,\phi)$ are the spherical coordinates of the point $x=(x_1,x_2,x_3)\in \Gamma_R,$
$u_{mn}=\int_{0}^\pi\int_{0}^{2\pi}\psi(R,\theta',\phi')
\overline{Y_{mn}(\theta',\phi')}\sin\theta'\,d\,\theta'd\,\phi',$
$h_{n}^{(j)}(x)$ are the spherical Hankel functions $(j=1,2)$,
$Y_{mn}(\theta,\phi):=\frac 1{2}\sqrt{\frac{(2n+1)(n-|m|)!}
{\pi(n+|m|)!}}P_n^{|m|}(\cos\theta)e^{im\phi}$ are the normalized spherical functions, and
$P_n^{m}(t)$ are the associated Legendre functions ($-n\leq m\leq n,$ $n=0,1,\ldots,\infty$);
and formula \eqref{g2ext} should be replaced by
\begin{equation}\label{g2ext'}
\varphi(r_P,\theta_P,\phi_P)=\frac 1{4\pi}\sum_{n=0}^\infty
\frac {h_{n}^{(1)}(kr_P)}
{h_{n}^{(1)}(kR)}\sum_{m=-n}^nu_{mn}Y_{mn}(\theta_P,\phi_P),\quad r_P\geq R.
\end{equation}
\end{predlllll}

\begin{predlllll} The analysis performed in Section 2 enables us to state that
 all the results of this section remain valid, for example, when in BVP \eqref{g20}$-$\eqref{g23} the Helmholtz equation \eqref{g21} is replaced by
\begin{equation}\label{g21gn}
-(\Delta+ k^2n(x))\varphi(x)=f(x)\quad \mbox{in}\quad\mathbb R^n\setminus\bar\Omega,
\end{equation}
where $n=2,3,$ $k=\mbox{const}>0,$ the function $n\in C(\mathbb R^n\setminus\Omega)$ is positive in
$\mathbb R^n\setminus\Omega$,
and $n(x)=1$ in the domain
$\Omega_{R_0}\setminus\bar\Omega$ for a certain $R_0>0.$

In this case one has to choose in the equation systems \eqref{g20gr}$-$\eqref{g23yxr} and
\eqref{g20gr2}$-$\eqref{g23yxr2} which specify minimax mean-square estimates the value of $R$
greater than $R_0$, set $k=\bar k,$ and replace in equations \eqref{g21gr}, \eqref{g21yxr},
\eqref{g21gr2}, and \eqref{g21yxr2} $k^2$ by  $k^2n(x).$
\end{predlllll}

Note also that applying DtN finite-element methods to  problems
 \eqref{g20gr}$-$\eqref{g23yxr} and
\eqref{g20gr2}$-$\eqref{g23yxr2} one can construct approximate methods of their solution.

\newpage
\makeatletter
\renewcommand{\section}{\@startsection{section}{1}%
{\parindent}{3.5ex plus 1ex minus .2ex}%
{2.3ex plus.2ex}{\normalfont\Large{\bf PART\ \ }}} \makeatother

\begin{center}
\section[
Minimax estimation of the solutions to the boundary value  problems
from observations distributed on a system of surfaces.
Reduction to a surface integral equation systems
]{}
\end{center}

\begin{quote}
{\bf Minimax estimation of the solutions to the boundary value  problems
from observations distributed on a system of surfaces.
Reduction to a surface integral equation systems
}
\end{quote}

\subsection{Statement of the problem}
Before to formulate the estimation problem which is the subject of
analysis of the present chapter, let us introduce the necessary
notations and functional spaces

Let $\Lambda$ be closed or unclosed $(n-1)$-dimensional Lipschitz
surface in $\mathbb R^n$. By $d\Lambda$ we will denote the element
of measure on surface $\Lambda$ and by $L^2(\Lambda)$ the space of
square integrable
functions on
 $\Lambda.$

Let $\gamma$ be an unclosed $(n-1)$-dimensional smooth $C^{\infty}$-surface
in $\mathbb R^n,$ $\partial \gamma$ its boundary
whose points do not belong to $\gamma,$ $\partial
\gamma\cap\gamma=\emptyset,$ and $\hat\gamma$ a closed smooth $C^{\infty}$-surface ($(n-1)$-dimensional manifold without edge) that contains $\gamma,$ $\gamma\subset\hat\gamma,$ and divides
$\mathbb R^n$ into two domains, bounded and unbounded. Set
$$
H^{1/2}(\gamma):=\{v|_{\gamma}:\,\,v\in H^{1/2}(\hat\gamma)\}.
$$
The norm in space $H^{1/2}(\gamma)$ is determined according to the formula
$$
v\in H^{1/2}(\gamma)\Longrightarrow
\|v\|_{H^{1/2}(\gamma)}=\inf_{V\in
H^{1/2}(\hat\gamma),\,V|_{\gamma}=v }\|V\|_{H^{1/2}(\hat\gamma)}.
$$
Denote by  $H^{-1/2}(\gamma)=(H^{1/2}(\gamma))'$ the space dual to
 $H^{1/2}(\gamma).$ Below, the duality relation
$<r,w>_{\gamma}$ on $H^{-1/2}(\gamma)\times H^{1/2}(\gamma)$ will be also denoted by  $\int_{\gamma}r\bar w\,d\gamma$ because for this relation the condition (*) on page
\pageref{k2} is valid in which $\Gamma$ should be replaced by $\gamma.$
Note that the elements of $H^{-1/2}(\gamma)$ extended to
$\hat\gamma\setminus\gamma$ by zero values belong to $H^{-1/2}(\hat\gamma).$

Let $\rho(x)$ be a function regular on surface $\gamma$ which is equivalent to the distance from a point
 $x$ to the boundary $\partial \gamma$
of $\gamma$ (this distance will be denoted by
$d(x,\partial\gamma)$) in a vicinity\footnote{It means that $\lim_{x\to
x_0}\frac{\rho(x)}{d(x,\partial\gamma)}=c=\mbox{const}\neq 0\quad
\forall x_0\in \partial\gamma.$ Such functions exist because
$\partial \gamma$ is an infinitely differentiable manifold
\cite{BIBLlima}.} of $\partial \gamma.$ Following
\cite{BIBLlima} and \cite{BIBLces} introduce the space
$$
H_{0\,0}^{1/2}(\gamma):=\{u\in H^{1/2}(\gamma),\, \rho^{-1/2}u\in
L^2(\gamma)\}=\{u\in H^{1/2}(\gamma),\,\tilde u \in H^{1/2}(\hat
\gamma)\}
$$
with the norm
$$
||u||_{H_{0\,0}^{1/2}(\gamma)}=\left(||u||^2_{H^{1/2}
(\gamma)}+||\rho^{-1/2}u||^2_{L^2(\gamma)}\right)^{1/2},
$$
where function $\tilde u$ denotes the extension of $u$ by $0$ on
$\hat\gamma\setminus\gamma.$

By $\left(H_{0\,0}^{1/2}(\gamma)\right)'$ we will denote a
space conjugate
to $H_{0\,0}^{1/2}(\gamma).$ Then, in line with
\cite{BIBLces}, p. 43, we have
$$
\left(H_{0\,0}^{1/2}(\gamma)\right)'=\{f=f_0+f_1,\,f_0 \in
H^{-1/2}(\gamma),\, \rho^{1/2}f_1\in L^2(\gamma)\}.
$$

Now let us formulate the estimation problem. Assume that the state
$\varphi (x)$ of a system is determined as a solution
 to the Neumann problem\footnote{In this chapter we will restrict ourselves to the case $n=3$. The results obtained
 for $n=3$ remain valid for $n=2$ after corresponding replacement of Sommerfeld radiation condition and
 fundamental solution.}
\begin{equation} \label{1}
\varphi \in H^1_{\mbox{\rm loc}}((\mathbb R^3\setminus \bar\Omega),\Delta),
\end{equation}
\begin{equation} \label{2}
(\Delta+k^2)\varphi(x)=0\quad \mbox {in}\quad \mathbb R^3\setminus
\bar\Omega,
\end{equation}
\begin{equation} \label{3}
\frac{\partial\varphi}{\partial\nu_{A}}=h \quad \mbox {on} \quad
\Gamma,
\end{equation}
\begin{equation}\label{3'}
\frac{\partial\varphi}{\partial
r}-ik\varphi=o(1/r),\,\,r=|x|,\,\,r\to \infty,
\end{equation}
where\footnote{It is known (see \cite{BIBLob}, page 221) that
there exists a uniquely determined continuous operator which we
denote by $\frac{\partial }{\partial \nu}$ and which maps space
$H^1(\Omega,\Delta)$ into space $H^{-1/2}(\Gamma)$ and is such
that $\forall u \in H^1(\Omega,\Delta),\,\forall v\in H^1(\Omega)$
the following representation  (Green's formula) holds: $
\sum_{i=1}^n\int_{\Omega}\frac{\partial u}{\partial
x_i}\frac{\partial v}{\partial x_i}dx=-\int_{\Omega}\Delta
u\,v\,dx+\int_{\Gamma}\frac{\partial u }{\partial \nu}v\,d\Gamma,
$ where the integrals over $\Gamma$ should be understood as the
duality relations $<\frac{\partial u }{\partial \nu},\mu v>$ on
$H^{-1/2}(\Gamma)\times H^{1/2}(\Gamma).$ This operator is called
the normal derivative in relation to $-\Delta$;
the operator $\frac{\partial}{\partial\nu}$ is defined by $
\frac{\partial u}{\partial \nu}= \sum_{i,j=1}^{n}\,
\,\frac{\partial u}{\partial x_j}\, \cos (\nu,x_i) $ when $u\in
C^{\infty}(\bar\Omega)$, where $\nu$ is the unit normal vector of
$\Gamma$ external with respect to domain $\Omega$ and $\cos
(\nu,x_i)$ is the $i$th directional cosine of $\nu.$} $\Omega\in
\mathbb R^3$ is a bounded domain with a connected complement such
that $\partial \Omega=\Gamma$ is a surface of class $C^2,$ $h\in
L^2(\Gamma).$

It is known that problem
 (\ref{1})$-$(\ref{3'}) is uniquely solvable.

Let $\gamma_i,$ $i=\overline{1,N},$ be smooth simply-connected oriented surfaces in $\mathbb R^3$ with smooth boundaries
$\partial \gamma_i,$ $\partial \gamma_i\cap\gamma_i=\emptyset,$ contained in the domain \label{wzaaq}$\mathbb R^3\setminus\bar\Omega$ that have no intersections pairwise, $\bar\gamma_i\cap\bar\gamma_j=\emptyset,$ $i\neq
j,$ $\bar\gamma_i\subset\mathbb R^3\setminus\bar\Omega.$ Let the orientation of
$\gamma_i$ be determined by a continuous family of normals
$\nu(x),\,x\in \gamma_i.$

Assume that on surfaces $\gamma_i$
the following functions are observed
\begin{equation} \label{4}
y_{i}\n(x)=\int_{\gamma_i}
K_{i}\p(x,y)\varphi(y)\,d\gamma_{i_y}+\int_{\gamma_i}
K_{i}\q(x,y)\frac{\partial\varphi(y)}{\partial\nu}\,d\gamma_{i_y}
+ \eta_{i}\n(x),
\end{equation}
\begin{equation} \label{5}
y_{i}\m(x)=\int_{\gamma_i}
K_{i}\s(x,y)\varphi(y)\,d\gamma_{i_y}+\int_{\gamma_i}
K_{i}\w(x,y)\frac{\partial\varphi(y)}{\partial\nu}\,d\gamma_{i_y}
+\eta_{i}\m(x),
\end{equation}
$$
x \in \gamma_i, \quad i= \overline{1,N},
$$
where $\varphi(x)$ is the solution\footnote{Note that
for any subdomain
$\omega$ such that $\bar\omega\in \Omega,$ the solution $\varphi$
to problem (\ref{1})$-$(\ref{3}) belongs to $H^2(\omega)$; therefore, according to the trace theorem,
$\varphi|_{\gamma_i},\frac{\partial\varphi(y)}{\partial\nu}
\left.\right|_{\gamma_i}\in L^2(\gamma_i),$ $i= \overline{1,N},$
and integrals (\ref{4})$-$(\ref{5}) make sense.\label{page3}} to BVP (\ref{1})$-$(\ref{3}),
 $\eta_i\n(x)$ and $\eta_i\m(x)$ are the observation errors
that are choice functions of random fields defined on surfaces $\gamma_i;$
$K_i^{(r,j)}\in
L^2(\gamma_i\times\gamma_i),\,\,r,j=1,2,$ are functions defined on
$\gamma_i\times\gamma_i$; and
integral operators $G_i^{(j)}$ with kernels
$K_{i}^{(j,2)}(\xi,x)$ defined according to\\[-20pt]
\begin{equation}\label{ya1}
G_i^{(j)}\psi(x)=\int_{\gamma_i}\!\!
K_{i}^{(j,2)}(\xi,x)\psi(\xi)d\gamma_{i_{\xi}},\quad j=1,2,
\end{equation}\\[-27pt]
are linear bounded operators acting from $L^2(\gamma_i)$ to $H_{0\,0}^{1/2}(\gamma_i),\,\,
i=\overline{1,N}$ (as an example of such kernels, one may take degenerated kernels
$K_{i}^{(j,2)}(\xi,x)=\sum_{r=1}^la_r^{(ij)}(\xi)b_r^{(ij)}(x),$
where $a_r^{(ij)}\in L^2(\gamma_i),$ $b_r^{(ij)} \in
H_{0\,0}^{1/2}(\gamma_i)$).

From the physical viewpoint, observations of the form (\ref{4}),
(\ref{5}) enable one, e.g. in stationary problems of hydro
acoustics, to observe independently both the pressure and the
normal velocity component as well as their linear combinations on
a system of surfaces $\gamma_i,$ $i=\overline{1,N}$.

Denote by $G_0$ the set of functions $\tilde h,$ $\tilde h \in
L^2(\Gamma)$ that satisfy the condition
\begin{equation} \label{7}
\int_{\Gamma}|\tilde h-h_0|^2
q_1^2 d\Gamma\leq 1,
\end{equation}
where $h_0\in L^2(\Gamma)$ is a given function. By $G_1$ we denote the set of random vector-functions $\tilde
\eta(\cdot)=(\tilde\eta^{(1)}_1(\cdot),\ldots,
\tilde\eta^{(1)}_{N}(\cdot),$ $\tilde\eta^{(2)}_1(\cdot), \ldots,
\tilde\eta^{(2)}_{N}(\cdot));$ their components
$\tilde\eta^{(1)}_{i}(x)$ and $\tilde\eta^{(2)}_i(x)$ are random fields
defined on surfaces $\gamma_i,$ $i=\overline{1,N}$ having square integrable second moments and satisfying the conditions
\begin{gather} \label{6}
\mathbf E\tilde\eta_{i}^{(1)}(x)=0, \quad \mathbf E\tilde\eta_{i}^{(2)}(x) =0,\quad i= \overline{1,N},\\
\label{8}  \!\!\sum_{i=1}^{N}\! \int_{\gamma_i}\!\! \mathbf E|\tilde\eta\n_i(x)|^2 \left(r\n_i(x)\right)^2 d\gamma_i+
\sum_{i=1}^{N}\! \int_{\gamma_i}\!\! \mathbf E|\tilde\eta\m_i(x)|^2 \left(r\m_i(x)\right)^2 d\gamma_i\leq 1,
\end{gather}
where $q_1(x),r_i \n(x), r_i \m(x),\,\, i= \overline{1,N},$ are functions continuous on $\Gamma$ and $\bar
\gamma_i,$ respectively, that do not vanish on these sets.

Assume also that in equalities (\ref{1})$-$(\ref{3})
function $h(x)$ and the second moments  $\mathbf E|\eta\n_i(x)|^2$ and $\mathbf E|\eta\m_i(x)|^2$ of random fields
$\eta\n_i(x)$ and $\eta\m_i(x)$ in observations (\ref{4}) and
(\ref{5}) are not known exactly, and it is known only that
\begin{equation}\label{lllxxx}
h\in G_0,\quad
\eta(\cdot)=(\eta^{(1)}_{1},\ldots,\eta^{(1)}_{N}
(\cdot),\eta^{(2)}_1(\cdot), \ldots, \eta^{(2)}_{N}(\cdot))\in
G_1.
\end{equation}

Suppose that a function $l_0 \in L^2(\omega_0)$ is defined in a domain
$\omega_0,$ $\bar\omega_0\subset\mathbb R^3\setminus\bar\Omega$.
The problem is as follows: given observations (\ref{4}),
(\ref{5}) of the state $\varphi(x)$ of a system  described by the Neumann BVP (\ref{1})$-$(\ref{3}) under the conditions that $h\in
G_0$ and $\eta\in G_1,$ estimate the value of the linear functional
\begin{equation} \label{9}
l(\varphi)= \int_{\omega_0}\overline{l_0(x)}\varphi(x)\, dx
\end{equation}
in the class of estimates linear with respect to observations that have the form
\begin{equation} \label{10}
\widehat{l(\varphi)}=\sum_{i=1}^N\int_{\gamma_i}\left(\overline{u_i\n(x)}
y_i\n(x) +
\overline{u_i\m(x)}y_i\m(x)\right)\, d\gamma_i+c,
\end{equation}
where $u_i\n,u_i\m\in L^2(\gamma_i),\,\,i= \overline{1,N},\,\,c \in
\mathbb C.$

Put $u:=(u_1\n, \ldots ,
u_{N}\n,$ $u_1\m, \ldots ,u_{N}\m)$ $\in H:=\left(L^2(\gamma_1)
\times\ldots \times L^2(\gamma_N)\right)^2.$

\begin{predll} {\it An estimate
$$
\widehat{\widehat
{l(\varphi)}}=\sum_{i=1}^N\int_{\gamma_i}\left(\overline{\hat u_i\n(x)}
y_i\n(x) + \overline{\hat
u_i\m(x)}y_i\m(x)\right)\, d\gamma_i+\hat c,
$$
in which functions $\hat u_i\n,$ $\hat u_i\m$ and number $\hat c$
are determined from the condition
\begin{equation} \label{11}
\sup_{\tilde h \in G_0,\, \tilde \eta
\in G_1}\mathbf E|l(\tilde\varphi)-\widehat
{l(\tilde\varphi)}|^2 \to \inf_{u\in H,\,c\in \mathbb C} ,
\end{equation}
where
\begin{equation} \label{llx}
\widehat
{l(\tilde\varphi)}=\sum_{i=1}^N\int_{\gamma_i}\left(\overline{u_i\n(x)}
\tilde y_i\n(x) +
\overline{u_i\m(x)}\tilde y_i\m(x)\right)\,
d\gamma_i+c,
\end{equation}
\begin{equation} \label{4iu}
\tilde y_{i}\n(x)=\int_{\gamma_i}
K_{i}\p(x,y)\tilde \varphi(y)\,d\gamma_{i_y}+\int_{\gamma_i}
K_{i}\q(x,y)\frac{\partial\tilde \varphi(y)}{\partial\nu}\,d\gamma_{i_y}
+ \tilde \eta_{i}\n(x),
\end{equation}
\begin{equation} \label{5iu}
\tilde y_{i}\m(x)= \int_{\gamma_i}
K_{i}\s(x,y)\tilde \varphi(y)\,d\gamma_{i_y}+\int_{\gamma_i}
K_{i}\w(x,y)\frac{\partial\tilde \varphi(y)}{\partial\nu}\,d\gamma_{i_y}
+\tilde \eta_{i}\m(x),
\end{equation}
$$
x \in \gamma_i, \quad i= \overline{1,N},
$$\\[-30pt]
and $\tilde\varphi(x)$ is the solution to the Neumann BVP at  $h=\tilde h,$
will be called a minimax estimate of expression (\ref{9}).

The quantity
\begin{equation} \label{12d}
\sigma:=\{\sup_{\tilde h \in G_0,\,
\tilde \eta \in G_1}\mathbf E|l(\tilde \varphi)-\widehat{\widehat
{l(\tilde\varphi)}}|^2\}^{1/2}
\end{equation}
will be called the error of the minimax estimation of $l(\varphi).$}
\end{predll}



\subsection{Auxiliary statements}

In this section we will prove that finding the minimax estimate is equivalent to a certain problem of optimal control of a system described by elliptic equations with conjugation conditions on surfaces
 $\gamma_i,$ $i=\overline{1,N}.$

In order to state the conjugation problems under study
and prove the existence of their solutions it is necessary to introduce the corresponding  Sobolev spaces
and trace theorems for surfaces with edges. First, let us formulate several definitions.

For $f \in \mathcal D(\mathbb R^n)$ the function $u(x) :=
\int_{\mathbb R^n}
\Phi_k(x,y)f(y)dy$ solves $-\Delta u-k^2u = f$
(and complies with Sommerfeld radiation conditions). Here
$$
\Phi_k(x,y) = \left\{
    \begin{array}{ccc}
    \frac{i}{4}H_0^{(1)} (k|x-y|) &\mbox{in}& \mathbb R^2,  \\
    \frac 1{4\pi}\frac{e^{ik|x-y|}}{|x-y|}&\mbox{in}& \mathbb R^3
    \end{array}
 \right.
$$
is the fundamental solution to the Helmholtz operator.

Introduce the Newton potential operator $(N_kf)(x) :=
\int_{\mathbb R^n} \Phi_k(x,y)f(y)dy.$
 \begin{predllll}(see, for example, \cite{BIBLCak})
$N_k$ can be extended to a bounded operator $N_k :
H^{-1}_{\mbox{\rm\small
comp}}(\mathbb R^n) \to H^1_{\mbox{\rm\small
loc}}(\mathbb R^n).$
\end{predllll}
For $\psi,\chi \in C(\Gamma)$ introduce the functions
\begin{equation}\label{simple}
(\mathcal V^k_{\Gamma}\psi)(x) := \int_\Gamma\Phi_k(x,y)
\psi(y)\,d\Gamma_y,
\end{equation}
\begin{equation}\label{double}
(\mathcal W^k_{\Gamma} \chi)(x):= \int_\Gamma\frac{\partial
\Phi_k(x,y)}{\partial \nu_y} \chi(y)\,d\Gamma_y, \,\,\quad x\notin \Gamma,
\end{equation}\label{poten}
called the single and double layer potentials. Let us formulate the results
contained in \cite{BIBLces} in the following form.
\begin{predllll}
$\mathcal V^k_{\Gamma}$ and $\mathcal W^k_{\Gamma}$ can be
extended to bounded operators $\mathcal V^k_{\Gamma}:
\,H^{-1/2}(\Gamma)\to H^1_{\mbox{\rm\small loc}}(\mathbb R^3)\cap
H^1_{\mbox{\rm\small loc}}\left(\Omega\cup(\mathbb R^3\setminus
\bar\Omega),\Delta\right)$ and $\mathcal
W^k_{\Gamma}:\,H^{1/2}(\Gamma)\to H^1_{\mbox{\rm\small
loc}}\left(\Omega\cup(\mathbb R^3\setminus
\bar\Omega),\Delta\right).$\footnote{ Here, the boundedness of,
e.g., the operator $\mathcal V_\Gamma^k: \,H^{-1/2}(\Gamma)\to
H^1_{\mbox{\rm loc}}(\mathbb R^3)$ means that for any
cutoff function  $\alpha\in C_{\mbox{\rm comp}}^\infty(\mathbb R^3)$ the operator
$\alpha\mathcal
V_\Gamma^k:\,H^{-1/2}(\Gamma)\to H^1(\mathbb R^3)$ is bounded.}
\end{predllll}

Let $\gamma_0$ be a smooth bounded oriented simply-connected unclosed smooth surface
in $\mathbb R^3$ with a smooth boundary $\partial \gamma_0,$ $\partial \gamma_0\cap\gamma_0=\emptyset.$
Let its orientation be specified by a continuous family of unit normals
$\nu(x),\,x\in \gamma_0.$ Denote by $\Omega_0$
such a bounded open set in $\mathbb R^3$ with a smooth boundary
$\partial\Omega_0=:\hat\gamma_0$ containing surface $\gamma_0$ that the normal vector $\nu$ to $\gamma_0$ is directed outside
$\Omega_0$; by $\gamma_{0-}$ we denote the side of $\gamma_0,$
whose orientation coincides with that of the external side of surface $\partial\Omega_0$ and by $\gamma_{0+}$ the opposite side of $\partial\Omega_0.$
\label{vvk}

Let $\gamma_1,\ldots,\gamma_N$ be the surfaces of the type introduced above
(see  p.~\pageref{wzaaq}) with the sides $\gamma_{1+},\gamma_{1-},\ldots,
\gamma_{N+},\gamma_{N-}$ respectively.

 Set\label{Omega'}
$\Omega'=(\mathbb R^3\setminus\bar\Omega)\setminus
\cup_{i=1}^{N}\bar\gamma_i$ and for any function $v$ defined in
$\Omega'$ denote by $v|_{\gamma_{i+}}$ and $v|_{\gamma_{i-}}$ the
restriction of $v$ to $\gamma_{i+} $ and to $\gamma_{i-},$
$i=\overline{1,N},$ respectively.


Assume that the functions $g,\tilde g\in H^{0}_{\mbox{\rm\small
comp}}(\Omega')=:L^2_{\rm \small comp}(\Omega'),$
$\alpha, \tilde \alpha\in H^{-1/2}(\Gamma),$ $\omega_i^{(1)},
\varrho_i^{(1)}\in
H_{0\,0}^{1/2}(\gamma_i),$ and $\omega_i^{(2)},\varrho_i^{(2)}\in H^{-1/2}(\gamma_i),$
$i=\overline{1,N}$ are defined in domain $\Omega'$.

Consider two problems.

1. Find function $u$ that satisfies the conditions\label{page1}
\begin{equation}\label{wzv1}
u\in H^1_{\mbox{\rm loc}}(\Omega',\Delta),
\end{equation}
\begin{equation}\label{wzv2}
-(\Delta+k^2)u(x)=g(x) \quad\mbox{in}\quad \Omega',\end{equation}
\begin{equation}\label{wzv3}
\frac{\partial u}{\partial \nu}=\alpha \quad\mbox{on}\quad
\Gamma,
\end{equation}
\begin{equation}\label{wzv4}
[u]_{\gamma_{i}}=\omega_{i}\n \quad \mbox{on}\quad \gamma_i,\quad
i=\overline{1,N},
\end{equation}
\begin{equation}\label{wzv5}
\left[\frac{\partial
u}{\partial\nu}\right]_{\gamma_{i}}=\omega_{i}\m \quad
\mbox{on}\quad \gamma_i,\quad i=\overline{1,N}.
\end{equation}
\begin{equation}\label{wzv5'}
\frac{\partial u}{\partial
r}-iku=o(1/r),\,\,r=|x|,\,\,r\to \infty,\,\,\mbox{if}\,\, \mbox{\rm
Im\,}k\geq 0.
\end{equation}
Here
$[u]_{\gamma_0}=u|_{\gamma_{i+}}-u|_{\gamma_{i-}}\in
H_{0\,0}^{1/2}(\gamma_i),$ $u|_{\gamma_{i+}},u|_{\gamma_{i-}}\in
H^{1/2}(\gamma_i),$ $\left[\frac{\partial
u}{\partial\nu}\right]_{\gamma_i}=\left.\frac{\partial
u}{\partial\nu}\right|_{\gamma_{i+}}-\left.\frac{\partial
u}{\partial\nu}\right|_{\gamma_{i-}}\in H^{-1/2}(\gamma_i),$
$\left.\frac{\partial
u}{\partial\nu}\right|_{\gamma_{i+}},\left.\frac{\partial
u}{\partial\nu}\right|_{\gamma_{i-}}\in
\left(H_{0\,0}^{1/2}(\gamma_i)\right)',$ $i=\overline{1,N},$
$\left.\frac{\partial u}{\partial\nu}\right|_{\Gamma}\in
H^{-1/2}(\Gamma),$ and (\ref{wzv2})$-$(\ref{wzv5}) should be understood as equalities of elements from spaces, respectively, $L^2(\Omega'),$
$H^{-1/2}(\Gamma),$ $H_{0\,0}^{1/2}(\gamma_i),$ and
$H^{-1/2}(\gamma_i),$ $i=\overline{1,N}$.

2. Find function $v$ that satisfies the conditions
\begin{equation}\label{wzuv35}
v\in H^1_{\mbox{\rm loc}}(\Omega',\Delta),
\end{equation}
\begin{equation} \label{wzuv36} -(\Delta+\bar k^2)v(x)=\tilde g(x)
\quad\mbox{in}\quad \Omega',
\end{equation}
\begin{equation}
\label{wzuv37} \frac{\partial v}{\partial \nu}=\tilde \alpha
\quad\mbox{on}\quad \Gamma,
\end{equation}
\begin{equation}
 \label{wzuv38} [v]_{\gamma_{i}}
 =\varrho_{i}\n \quad \mbox{on}\quad
\gamma_i,\quad i=\overline{1,N},
\end{equation}
\begin{equation}
\label{wzuv39}\left[\frac{\partial
v}{\partial\nu}\right]_{\gamma_{i}}=\varrho_{i}\m
\quad \mbox{on}\quad \gamma_i,\quad i=\overline{1,N}.
\end{equation}
\begin{equation}\label{wzuv39'}
\frac{\partial v}{\partial
r}+i\bar ku=o(1/r),\,\,r=|x|,\,\,r\to \infty,
\end{equation}
where
 $[v]_{\gamma_0}=v|_{\gamma_{i+}}-v|_{\gamma_{i-}}\in
H_{0\,0}^{1/2}(\gamma_i),$ $v|_{\gamma_{i+}},v|_{\gamma_{i-}}\in
H^{1/2}(\gamma_i),$ $\left[\frac{\partial
v}{\partial\nu}\right]_{\gamma_i}=\left.\frac{\partial
v}{\partial\nu}\right|_{\gamma_{i+}}-\left.\frac{\partial
v}{\partial\nu}\right|_{\gamma_{i-}}\in H^{-1/2}(\gamma_i),$
$\left.\frac{\partial
v}{\partial\nu}\right|_{\gamma_{i+}},\left.\frac{\partial
v}{\partial\nu}\right|_{\gamma_{i-}}\in
\left(H_{0\,0}^{1/2}(\gamma_i)\right)',$ $i=\overline{1,N},$
$\left.\frac{\partial v}{\partial\nu}\right|_{\Gamma}\in
H^{-1/2}(\Gamma),$ and (\ref{wzuv36})$-$(\ref{wzuv39})
should be understood as equalities of elements from the corresponding spaces.

\begin{predlllll} The choice of spaces for problems 1 and 2 is governed\label{wzt1} by the trace theorems (see \cite{BIBLces}, pp. 44, 45 and \cite{BIBLob}, pp. 180, 181).
\end{predlllll}

In order to prove the existence and uniqueness of solutions to problems 1 and 2 we formulate one more known result.
Namely, let $\gamma_0$ and $\hat\gamma_0$ be the surfaces introduced on p. \pageref{vvk}.
By virtue of the definition of spaces $
H_{0\,0}^{1/2}(\gamma_0)$ and $H^{-1/2}(\gamma_0)$, the elements of these spaces extended by zero on $\hat\gamma_0\setminus\gamma_0$ are the elements of $H^{1/2}(\hat\gamma_0)$ and $H^{-1/2}(\hat\gamma_0)$, respectively.

Set
\begin{equation}\label{}
X^1_{\gamma_0}:=\{u\in D'(\mathbb R^3),\,\, u\in
H^1(\Omega_R\setminus\bar\gamma_0),\,\,\Delta u\in H^1(\Omega_R \setminus\bar\gamma_0)\,\,\forall\,\,\mbox{sufficiently large}\,\,R\}.
\end{equation}
If the tilde sign marks the zero extension on $\hat\gamma_0\setminus\bar\gamma_0$ of an element defined on $\gamma_0,$ the following statement holds.
\begin{pred}
Let\label{ol} $\rho\in H_{0\,0}^{1/2}(\gamma_0)$ and $
\rho'\in H^{-1/2}(\gamma_0).$ Then $\mathcal W^{k}_{\gamma_0} \rho:=\mathcal W^{k}_{\hat\gamma_0}\tilde \rho$ and $\mathcal V^{k}_{\gamma_0} \rho':=\mathcal V^{k}_{\hat\gamma_0}\tilde \rho'$ belong to
$X^1_{\gamma_0},$ and $\mathcal V \rho'\in H^1(\mathbb R^3).$
\end{pred}

Formulate the BVP: find $u\in H^1_{\rm loc}(\mathbb R^3\setminus\bar\gamma_0)$
satisfying
\begin{equation}\label{z0}
u\in H^1_{\rm loc}(\mathbb R^3\setminus\bar\gamma_0),
\end{equation}
\begin{equation}\label{z1}
\Delta u(x)+k^2u(x)=0\,\, \mbox{in}\,\,\mathbb R^3\setminus
\bar\gamma_0,
\end{equation}
\begin{equation}\label{z2}
[u]_{\gamma_0}=\rho,\quad\left[\frac{\partial u
}{\partial\nu}\right]_{\gamma_0}=\rho',
\end{equation}
\begin{equation}\label{z3}
\frac{\partial u}{\partial
r}-iku=o(1/r),\,\,r=|x|,\,\,r\to \infty.
\end{equation}
The next statement is an immediate corollary of the last theorem.
\begin{pred}\label{ska}
Let $\rho\in H_{0\,0}^{1/2}(\gamma_0)$ and $
\rho'\in H^{-1/2}(\gamma_0).$ Then BVP (\ref{z0})$-$(\ref{z3})
has a unique solution $u\in X^1_{\gamma_0}$ which, for $x$ outside $\gamma_0,$ can be given by the formula
\begin{equation}\label{z4}
u=\mathcal V^{k}_{\gamma_0}\rho'-\mathcal W^{k}_{\gamma_0}\rho.
\end{equation}
Here $[u]_{\gamma_0}=u|_{\gamma_{0+}}-u|_{\gamma_{0-}}\in
H_{0\,0}^{1/2}(\gamma_0),$ $u|_{\gamma_{0+}},u|_{\gamma_{0-}}\in
H^{1/2}(\gamma_0),$ $\left[\frac{\partial
u}{\partial\nu_{A}}\right]_{\gamma_0}=\left.\frac{\partial
u}{\partial\nu_{A}}\right|_{\gamma_{0+}}-\left.\frac{\partial
u}{\partial\nu_{A}}\right|_{\gamma_{0-}}\in H^{-1/2}(\gamma_0),$
and $\left.\frac{\partial
u}{\partial\nu_{A}}\right|_{\gamma_{0+}},\,\left.\frac{\partial
u}{\partial\nu_{A}}\right|_{\gamma_{0-}}\in
\left(H_{0\,0}^{1/2}(\gamma_0)\right)'.$
\end{pred}

Now let us prove, e.g. for problem 2, that the following statement is valid.
\begin{pred}
BVP (\ref{wzuv35})$-$(\ref{wzuv39'}) is uniquely solvable and
the estimate
\begin{multline}\label{yr3}
\|v\|_{H^1(\Omega'\cap\Omega_R)}\leq C_0\Bigl(\|\tilde g\|_{H^{-1}(\Omega_0)}+\|\tilde\alpha\|_{H^{-1/2}(\Gamma)}\Bigr.\\ \Bigl.+
\sum_{i=1}^N\| \varrho_{i}\n\|_{H_{0\,0}^{1/2}(\gamma_i)}+\sum_{i=1}^N\| \varrho_{i}\m\|_{H^{-1/2}(\gamma_i)}\Bigr),
\end{multline}
holds, where  $\Omega_0$ is the support of function $\tilde g$ and $C_0$ is a constant which does not depend
on $\tilde g,$ $\tilde\alpha,$ $\varrho_{i}\n$, and $\varrho_{i}\m,$
$i=\overline{1,N}.$
\end{pred}
\begin{proof}
Set
\begin{equation*}
 v_1(x)=\sum_{i=1}^N\mathcal V^{-\bar k}_{\gamma_i} \varrho_{i}\m(x) -
 \sum_{i=1}^N\mathcal W^{-\bar k}_{\gamma_i} \varrho_{i}\n(x),
\quad x\in\mathbb R^3\setminus(\cup_{i=1}^N\bar\gamma_i).
\end{equation*}
where, in accordance with definition on page \pageref{ol},
$\mathcal V^{-\bar k}_{\gamma_i}$ and $\mathcal W^{-\bar
k}_{\gamma_i}$ are single and double layer potentials determined
on unclosed surfaces $\gamma_i,$ $i=1,N,$ corresponding to the
wave number $-\bar k.$ Then by theorem \ref{ska},
function $v_1(x)$ is the unique solution to the problem
\begin{equation}\label{wzuv35i}
v_1\in X^1_{\cup_{i=1}^N\gamma_i},
\end{equation}
\begin{equation} \label{wzuv36i}
 -(\Delta+\bar k^2)v_1(x)=0
\quad\mbox{in}\quad \Omega',
\end{equation}
\begin{equation}
\label{wzuv37i} \frac{\partial v_1}{\partial \nu}=0
\quad\mbox{on}\quad \Gamma,
\end{equation}
\begin{equation}
 \label{wzuv38a}  [v_1]_{\gamma_{i}}
 =\varrho_{i}\n \quad \mbox{on}\quad
\gamma_i,\quad i=\overline{1,N},
\end{equation}
\begin{equation}
\label{wzuv39i}\left[\frac{\partial
v_1}{\partial\nu}\right]_{\gamma_{i}}=\varrho_{i}\m
\quad \mbox{on}\quad \gamma_i,\quad i=\overline{1,N}.
\end{equation}
\begin{equation}\label{wzuv39i'}
\frac{\partial v_1}{\partial
r}+i\bar kv_1=o(1/r),\,\,r=|x|,\,\,r\to \infty,
\end{equation}
Here $[v_1]_{\gamma_0}=v_1|_{\gamma_{i+}}-v_1|_{\gamma_{i-}}\in
H_{0\,0}^{1/2}(\gamma_i),$ $v_1|_{\gamma_{i+}},v_1|_{\gamma_{i-}}\in
H^{1/2}(\gamma_i),$ $\left[\frac{\partial
v_1}{\partial\nu}\right]_{\gamma_i}=\left.\frac{\partial
v_1}{\partial\nu}\right|_{\gamma_{i+}}-\left.\frac{\partial
v_1}{\partial\nu}\right|_{\gamma_{i-}}\in H^{-1/2}(\gamma_i),$
$\left.\frac{\partial
v_1}{\partial\nu}\right|_{\gamma_{i+}},\left.\frac{\partial
v_1}{\partial\nu}\right|_{\gamma_{i-}}\in
\left(H_{0\,0}^{1/2}(\gamma_i)\right)',$ $i=\overline{1,N},$
and $\left.\frac{\partial v_1}{\partial\nu}\right|_{\Gamma}\in
H^{-1/2}(\Gamma).$
Denote by $v_2$ the unique solution to the problem
\begin{equation}\label{wzuv37a0}
v_2\in H^1_{\mbox{\rm\small loc}}\left((\mathbb R^3\setminus \bar\Omega),\Delta\right)
\end{equation}
\begin{equation}\label{wzuv37a1}
 -(\Delta +\bar k^2)v_2(x)=0,\quad x\in\mathbb R^3\setminus\bar\Omega,
\end{equation}
\begin{equation}\label{wzuv37a}
\frac{\partial v_2}{\partial \nu}=-\frac{\partial v_1}{\partial
\nu}
\quad\mbox{on}\quad \Gamma,
\end{equation}
\begin{equation}\label{wzuv37'}
\frac{\partial v_2}{\partial
r}+i\bar kv_2=o(1/r),\,\,r=|x|,\,\,r\to \infty,
\end{equation}
and by $v_3$ the unique solution to the problem
\begin{equation}\label{wzuv37a2}
v_3\in H^1_{\mbox{\rm\small loc}}\left((\mathbb R^3\setminus \bar\Omega),\Delta\right)
\end{equation}
\begin{equation}\label{wzuv37a3}
 -(\Delta +\bar k^2)v_3(x)=\tilde g,\quad x\in\mathbb R^3\setminus\bar\Omega,
\end{equation}
\begin{equation}\label{wzuv37aa}
\frac{\partial v_3}{\partial \nu}=\tilde \alpha
\quad\mbox{on}\quad \Gamma,
\end{equation}
\begin{equation}\label{wzuv37''}
\frac{\partial v_3}{\partial
r}+i\bar kv_3=o(1/r),\,\,r=|x|,\,\,r\to \infty.
\end{equation}
Then the function
\begin{equation}\label{yr-3}
v(x)=v_1(x)+v_2(x)+v_3(x)
\end{equation}
will be the unique solution to problem
(\ref{wzuv36})$-$(\ref{wzuv39'}). From equalities $\mathcal
V^{-\bar k}_{\gamma_i} \varrho_{i}\m(x)=\mathcal V^{-\bar
k}_{\hat\gamma_i} \tilde\varrho_{i}\m(x)$ and $\mathcal W^{-\bar
k}_{\gamma_i} \varrho_{i}\m(x)=\mathcal W^{-\bar k}_{\hat\gamma_i}
\tilde\varrho_{i}\m(x),$ and boundedness of operstors $\mathcal
V^{-\bar k}_{\hat\gamma_i}$ and $\mathcal W^{-\bar
k}_{\hat\gamma_i}$ in the corresponding spaces, we obtain the
following estimates for $v_1$ \footnote{Here and below $C_i$ are
constants that do not depend on the data of the problems in
question.}
\begin{equation*}
\|v_1\|_{H^1(\Omega'\cap\Omega_R,\Delta)}\leq C_1\sum_{i=1}^N\|\tilde \varrho_{i}\n\|_{H^{1/2}(\hat\gamma_i)},\quad
\|v_1\|_{H^1(\Omega'\cap\Omega_R,\Delta)}\leq  C_2\sum_{i=1}^N\|\tilde \varrho_{i}\m\|_{H^{-1/2}(\hat\gamma_i)}
\end{equation*}
and hence
\begin{equation}\label{yr-2}
\|v_1\|_{H^1(\Omega'\cap\Omega_R,\Delta)}\leq  C_3\left(\sum_{i=1}^N\| \varrho_{i}\n\|_{H_{0\,0}^{1/2}(\gamma_i)}+2\sum_{i=1}^N\| \varrho_{i}\m\|_{H^{-1/2}(\gamma_i)}\right).
\end{equation}
For functions $v_2$ and $v_3$ we apply estimate \eqref{g2cont1} to obtain
\begin{equation}\label{yr}
\|v_2\|_{H^1(\Omega_R\setminus\bar\Omega)}\leq C_4\|\gamma_Nv_1\|_{H^{-1/2}(\Gamma)},
\end{equation}
\begin{equation}\label{yr1}
\|v_3\|_{H^1(\Omega_R\setminus\bar\Omega)}\leq C_5\left(\|\tilde g\|_{H^{-1}(\Omega_R\setminus\bar\Omega)}+\|\tilde\alpha\|_{H^{-1/2}(\Gamma)}\right).
\end{equation}
Using the trace theorem (\cite{BIBLob}, pp. 180, 181) and (\ref{yr-2}), we prove the estimates
\begin{equation}\label{yr2}
\|\gamma_N v_1\|_{H^{-1/2}(\Gamma)}\leq C_5\|v_1\|_{H^1(\Omega'\cap\Omega_R,\Delta)}\leq C_6\left(
\sum_{i=1}^N\| \varrho_{i}\n\|_{H_{0\,0}^{1/2}(\gamma_i)}+\sum_{i=1}^N\| \varrho_{i}\m\|_{H^{-1/2}(\gamma_i)}\right).
\end{equation}
From (\ref{yr-3})$-$(\ref{yr2}) it follows that
\begin{multline*}
\|v\|_{H^1(\Omega'\cap\Omega_R)}\leq C_0\Bigl(\|\tilde g\|_{H^{-1}(\Omega'\cap\Omega_R)}+\|\tilde\alpha\|_{H^{-1/2}(\Gamma)}\Bigr.\\ \Bigl.+
\sum_{i=1}^N\| \varrho_{i}\n\|_{H_{0\,0}^{1/2}(\gamma_i)}+\sum_{i=1}^N\| \varrho_{i}\m\|_{H^{-1/2}(\gamma_i)}\Bigr).
\end{multline*}
\end{proof}

For problem 1 the corresponding theorem is obtained in a similar manner.

\subsection{General form of the guaranteed estimates and expression for the estimation error}

Introduce, for every fixed $u\in H=\left(L^2(\gamma_1)
\times\ldots \times L^2(\gamma_N)\right)^2$
the function $z(x;u)$ as a solution to the problem
\begin{equation} \label{13}
z \in H^1_{\mbox{\rm loc}}(\Omega',\Delta), \end{equation}
\begin{equation}\label{14}
-(\Delta+\bar k^2)z(x;u)=\chi_{\omega_0}(x)l_0(x)\quad \mbox {in} \quad
\Omega',
\end{equation}
\begin{equation}\label{15} \frac{\partial
z}{\partial\nu_{A^*}}=0 \quad \mbox {on} \quad \Gamma,
\end{equation}
\begin{equation*}
\!\!\!\!\!\!\!\!\!\![z(x;u)]_{\gamma_i}=
\int_{\gamma_i}\left[\overline{K_i\q(\xi,x)}u_i\n(\xi)+
\overline{K_i\w(\xi,x)}u_i\m(\xi)\right]\,d\gamma_{i_{\xi}},
\end{equation*}
\begin{equation}\label{16}
\left[\frac{\partial z(x;u)}{\partial\nu_{A^*}}\right]_{\gamma_i}=
-\int_{\gamma_i}\left[\overline{K_i\p(\xi,x)}u_i\n(\xi)+
\overline{K_i\s(\xi,x)}u_i\m(\xi)\right]\,d\gamma_{i_{\xi}}
\end{equation}
\begin{equation*}\quad \mbox {on}\quad \gamma_i,\quad i=
\overline{1,N}.
\end{equation*}
\begin{equation}\label{16'}
\frac{\partial z(x;u)}{\partial
r}+i\bar kz(x;u)=o(1/r),\,\,r=|x|,\,\,r\to \infty.
\end{equation}

\begin{predl}\label{lem2.1}
Finding the minimax estimate of
$l(\varphi)$ is equivalent to the problem of optimal control of the system described by BVP (\ref{13})$-$(\ref{16'}) with the cost function
\begin{equation*} I(u_1\n, \ldots ,u_{N}\n,u_1\m, \ldots ,u_{N}\m)
=\int_{\Gamma}q_1^{-2}(x)
z^2(x;u)\,d\Gamma
\end{equation*}
\begin{equation}\label{17}
+\sum_{i=1}^{N}\int_{\gamma_i}
(r_i\n(x))^{-2}(u_i\n(x))^2d\gamma_i+\sum_{i=1}^{N}\int_{\gamma_i}
(r_i\m(x))^{-2}(u_i\m(x))^2d\gamma_i \to \min_{u\in H}.
\end{equation}
\end{predl}
\begin{proof}
Denote by $\Omega_i$ an open
\label{ax} subdomain in $\mathbb R^3\setminus\bar\Omega$
($\bar\Omega_i\subset\mathbb R^3\setminus\bar\Omega$) such that $\partial\Omega_i$ contains $\gamma_i,$ its boundary
$\partial\Omega_i$ is simply-connected and smooth, and the normal vector $\nu$ to  $\gamma_i$ is directed
outside $\Omega_i$. We also assume that $\Omega_i\cap\Omega_j=\emptyset$
for all $i\neq j,$
$i,j=\overline{1,N}.$

Set $\tilde
\Omega_R:=(\Omega_R\setminus\bar\Omega)\setminus\cup_{i=1}^N\bar\Omega_i$ ($R$ is assumed to be sufficiently large),
$\hat\gamma_i=\partial\Omega_i$ and denote by $\hat\gamma_{i-}$ and
$\hat\gamma_{i+}$ the external and internal sides of surface
$\hat\gamma_i.$ Next, simplifying the notation in the surface integrals, denote by $z_+$, $\left(\frac{\partial
z}{\partial \nu_{A^*}}\right)_+$ and $z_-,$ $\left(\frac{\partial
z}{\partial \nu_{A^*}}\right)_-$ the traces
$\left.z\right|_{\gamma_{k+}},$ $\left.\frac{\partial z}{\partial
\nu_{A^*}}\right|_{\gamma_{k+}}$ (or
$\left.z\right|_{\hat\gamma_{k+}},$ $\left.\frac{\partial
z}{\partial \nu_{A^*}}\right|_{\hat\gamma_{k+}}$) and
$\left.z\right|_{\gamma_{k-}}$, $\left.\frac{\partial z}{\partial
\nu_{A^*}}\right|_{\gamma_{k-}}$ (or
$\left.z\right|_{\hat\gamma_{k-}},$ $\left.\frac{\partial
u}{\partial \nu_{A^*}}\right|_{\hat\gamma_{k-}}$) of functions $z(x;u)$
or $\frac{\partial z(x;u)}{\partial \nu_{A^*}}$ on sides
$\gamma_{k+}$ (or $\hat\gamma_{k+}$) and $\gamma_{k-}$ (or
$\hat\gamma_{k-}$) of surface
 $\gamma_k$ (or $\hat\gamma_k$).

Taking into consideration relationships (\ref{llx})$-$(\ref{5iu}),
(\ref{13})$-$(\ref{16'}), and applying to $\tilde\varphi(x)$ and
$z(x;u)$ in domains $\Omega_i,$ $i=1,N,$ and $\tilde \Omega_R$ the
second Green formula\footnote{One can apply the second Green
formula because $z(\cdot;u),\,\varphi \in H^1(\tilde \Omega_R),$
$z(\cdot;u),\,\varphi \in H^1( \Omega_i),$ $i=\overline{1,N},$ and
$\Delta z,\, \Delta\varphi \in L^2(\tilde \Omega_R),$ $\Delta z,\,
\Delta\varphi \in L^2(\Omega_i),$ $i=\overline{1,N},$}, we obtain,
using the equalities $[ z(\cdot;u)]_{\hat\gamma_i\setminus
\gamma_i} =\left[\frac{\partial z(\cdot;u)}{\partial
\nu}\right]_{\hat\gamma_i\setminus \bar \gamma_i}=0,$
$$
l(\tilde\varphi)-\widehat{l(\tilde\varphi)}=
$$
$$
=\int_{\omega_0}\overline{l_0(x)} \tilde\varphi(x)
\,dx-\sum_{i=1}^N \int_{\gamma_{i_x}}\left(\overline{u_i\n(x)}
\tilde y_i\n(x) +
\overline{u_i\m(x)}\tilde y_i\m(x)\right)\, d\gamma_i-c
$$
$$
=\int_{\tilde\Omega_R}\chi_{\omega_0}(x)\overline{l_0(x)} \tilde\varphi(x)
\,dx+\sum_{i=1}^N \int_{\Omega_i}\chi_{\omega_0}(x)\overline{l_0(x)} \tilde\varphi(x)
\,dx
$$
$$
-\sum_{i=1}^N \int_{\gamma_{i_x}}\left(\overline{u_i\n(x)}
\tilde y_i\n(x) +
\overline{u_i\m(x)}\tilde y_i\m(x)\right)\, d\gamma_i-c
$$
$$
=-\int_{\tilde\Omega_R}\tilde\varphi(x)\overline{(\Delta+\bar
k^2)z(x;u)}\,dx-\sum_{i=1}^N \int_{\Omega_i}\tilde\varphi(x)\overline{(\Delta+\bar
k^2)z(x;u)}\,dx
$$
$$
-\sum_{i=1}^N \int_{\gamma_i}
\overline{u_i\n(x)}\left[\int_{\gamma_i}K_i\p(x,\xi)
\tilde\varphi(\xi)\,d\gamma_{i_{\xi}}+
\int_{\gamma_i}K_i\q(x,\xi)\frac{\partial
\tilde\varphi(\xi)}{\partial \nu_{A^*}}
\,d\gamma_{i_{\xi}}\right]\,d\gamma_{i_x}
$$
$$
-\sum_{i=1}^N \int_{\gamma_i}
\overline{u_i\m(x)}\left[\int_{\gamma_i}K_i\s(x,\xi)
\tilde\varphi(\xi)\,d\gamma_{i_{\xi}}+
\int_{\gamma_i}K_i\w(x,\xi)\frac{\partial
\tilde\varphi(\xi)}{\partial \nu_{A^*}}
\,d\gamma_{i_{\xi}}\right]\,d\gamma_{i_x}
$$
$$
-\sum_{i=1}^N\int_{\gamma_i}
\overline{u_i\n(x)}\tilde\eta_i\n(x)\,d\gamma_i
-\sum_{i=1}^N\int_{\gamma_i}
\overline{u_i\m(x)}\tilde\eta_i\m(x)\,d\gamma_i-c
$$
$$
=-\int_{\tilde\Omega_R}(\Delta+
k^2)\tilde\varphi(x)\overline{z(x;u)}\,dx+
\int_{\Gamma}\overline{z}\frac{\partial\tilde\varphi}{\partial\nu}
\,d\Gamma -\sum_{i=1}^N\int_{\Omega_i}(\Delta+
k^2)\tilde\varphi(x)\overline{z(x;u)}\,dx
$$
$$
-\sum_{i=1}^N\int_{\hat\gamma_i}\left(\overline{z_-(x;u)}
\frac{\partial\tilde\varphi(x)}{\partial\nu}
-\tilde\varphi(x)\overline{\left(\frac{\partial z(x;u)}
{\partial\nu}\right)_{-}}\right)\,d\hat\gamma_i
$$
$$
+\sum_{i=1}^N\int_{\hat\gamma_i}\left(\overline{z_+(x;u)}
\frac{\partial\tilde\varphi(x)}{\partial\nu}
-\tilde\varphi(x)\overline{\left(\frac{\partial z(x;u)}
{\partial\nu}\right)_{+}}\right)\,d\hat\gamma_i
+\Sigma_R(z(\cdot;u),\tilde\varphi)
$$
$$
-\sum_{i=1}^N\int_{\gamma_i}\tilde\varphi(\xi)\int_{\gamma_i}
\left[K_i\p(x,\xi)\overline{u_i\n(x)}+K_i\s(x,\xi)\overline{u_i\m(x)}
\right]\,d\gamma_{i_x}
\,d\gamma_{i_{\xi}}
$$
$$
-\sum_{i=1}^N\int_{\gamma_i}\frac{\partial
\tilde\varphi(\xi)}{\partial\nu} \int_{\gamma_i}
\left[K_i\q(x,\xi)\overline{u_i\n(x)}+K_i\w(x,\xi)\overline{u_i\m(x)}
\right]\,d\gamma_{i_x}
\,d\gamma_{i_{\xi}}
$$
$$
-\sum_{i=1}^N\int_{\gamma_i}\overline{u_i\n(x)}\tilde\eta_i\n(x)\,d\gamma_i
-\sum_{i=1}^N\int_{\gamma_i}\overline{u_i\m(x)}\tilde\eta_i\m(x)\,d\gamma_i-c
$$
$$
=\int_{\Gamma}\overline{z} \tilde
h\,d\Gamma +\sum_{i=1}^N\int_{\gamma_i}(\overline{z_+(x,u)}-\overline{z_-(x,u)})
\frac{\partial\tilde\varphi(x)}{\partial\nu}\,d\hat\gamma_i
$$
$$
-\sum_{i=1}^N\int_{\gamma_i}\left[\overline{\left(\frac{\partial
z(x,u)}{\partial\nu}\right)_+} -\overline{\left(\frac{\partial
z(x,u)}{\partial\nu}\right)_-}\right]
\tilde\varphi(x)\,d\hat\gamma_i
$$
$$
-\sum_{i=1}^N\int_{\gamma_i}\tilde\varphi(x)\int_{\gamma_i}
\left[K_i\p(\xi,x)\overline{u_i\n(\xi)}+K_i\s(\xi,x)\overline{u_i\m(\xi)}\right]
\,d\gamma_{i_{\xi}}\,d\gamma_{i_x}+\Sigma_R(z(\cdot;u),\tilde\varphi)
$$
$$
-\sum_{i=1}^N\int_{\gamma_i}\frac{\partial
\tilde\varphi(x)}{\partial\nu} \int_{\gamma_i}
\left[K_i\q(\xi,x)\overline{u_i\n(\xi)}+K_i\w(\xi,x)\overline{u_i\m(\xi)}\right]
\,d\gamma_{i_{\xi}}\,d\gamma_{i_x}
$$
$$
-\sum_{i=1}^N\int_{\gamma_i}\overline{u_i\n(x)}\tilde\eta_i\n(x)\,d\gamma_i
-\sum_{i=1}^N\int_{\gamma_i}\overline{u_i\m(x)}\tilde\eta_i\m(x)\,d\gamma_i-c
$$
$$
=\int_{\Gamma} \tilde h
\overline{z}\,d\Gamma+\Sigma_R(z(\cdot;u),\tilde\varphi)
$$
\begin{equation}\label{kjj1}
-\sum_{i=1}^N\int_{\gamma_i}
\overline{u_i\n(x)}\tilde\eta_i\n(x)\,d\gamma_i
-\sum_{i=1}^N\int_{\gamma_i}
\overline{u_i\m(x)}\tilde\eta_i\m(x)\,d\gamma_i-c,
\end{equation}
where by $
\Sigma_R(z(\cdot;u),\tilde\varphi)$ we denote
$$
\Sigma_R(z(\cdot;u),\tilde\varphi):=\int_{\Gamma_R}\left(\overline{z(x;u)}
\frac{\partial\tilde\varphi(x)}{\partial\nu}
-\tilde\varphi(x)\overline{\left(\frac{\partial z(x;u)}
{\partial\nu}\right)}\right)\,d\Gamma_R
$$
Since $z(\cdot;u)$ and $\tilde\varphi$ satisfy, respectively, the Sommerfeld radiation conditions (\ref{16'}) and (\ref{3'}) we obtain an estimate for $\Sigma_R(z(\cdot;u),\tilde\varphi)$,
\begin{multline*}
\Sigma_R(z(\cdot;u),\tilde\varphi)=\int_{\Gamma_R}\overline{z(x;u)}
\left(\frac{\partial\tilde\varphi(x)}{\partial R}-ik\varphi(x)\right)\,d\Gamma_R
\\
-\int_{\Gamma_R}\tilde\varphi(x)\overline{\left(\frac{\partial z(x;u)}
{\partial R}+i\bar kz(x;u)\right)}\,d\Gamma_R
\end{multline*}
\begin{equation*}\label{kjj}
=\int_{\Gamma_R}O(1/R)o(1/R)d\Gamma_R
-\int_{\Gamma_R}O(1/R)o(1/R)d\Gamma_R=o(1)\quad \mbox{при}\quad R\to\infty.
\end{equation*}
From here, passing to the limit as $R\to\infty$ in \eqref{kjj1}, we obtain
$$
l(\tilde\varphi)-\widehat{l(\tilde\varphi)}=
\int_{\Gamma} \tilde h
\overline{z(\cdot;u)}\,d\Gamma
$$
$$
-\sum_{i=1}^N\int_{\gamma_i}
\overline{u_i\n(x)}\tilde\eta_i\n(x)\,d\gamma_i
-\sum_{i=1}^N\int_{\gamma_i}
\overline{u_i\m(x)}\tilde\eta_i\m(x)\,d\gamma_i-c,
$$

The latter equalities together with conditions (\ref{6})$-$(\ref{8}) and the known relation $\mathbf {D}\xi=\mathbf E\xi^2-|\mathbf E\xi|^2$
that couples dispersion $\mathbf D\xi$ of random variable $\xi$ and its expectation $\mathbf E\xi,$ yield
$$
\inf_{c \in \mathbb C} \sup_{\tilde h \in G_0,\, \tilde \eta \in
G_1} \mathbf E|l(\tilde\varphi)-\widehat{l(\tilde\varphi)}|^2=
$$
$$
=\inf_{c \in \mathbb C}\sup_{\tilde h \in G_0} \Bigl|
\int_{\Omega}\int_{\Gamma} \tilde h \overline{z(\cdot;u)}
\, d\Gamma-c \Bigr|^2
$$
\begin{equation}\label{18}
+\sup_{\tilde \eta \in G_1} \mathbf E\Bigl|\sum_{i=1}^{N}\int_{\gamma_i}
\overline{u_i\n(x)}\tilde\eta_i\n(x)\,d\gamma_i+
\sum_{i=1}^{N}\int_{\gamma_i}
\overline{u_i\m(x)}\tilde\eta_i\m(x)\,d\gamma_i\Bigr|^2.
\end{equation}
In order to calculate the supremum in the right-hand side of (\ref{18})
make use of the Cauchy$-$Bunyakovsky inequality.
Introducing the notation
$$
y=\int_{\Gamma} (\tilde h
-h_0)\overline{z(\cdot;u)}\, d\Gamma,
$$
we prove, using relation (\ref{7}), the inequality
$$ |y|\leq\Bigl\{
\int_{\Gamma}q_1^{-2}(x)
|z(x;u)|^2\,d\Gamma\Bigr\}^{\frac 1{2}}
\times\Bigl\{\int_{\Gamma}|\tilde h-h_0|^2 q_1^2 d\Gamma \Bigl.\Bigr\}^{\frac
1{2}}
$$
$$ \leq \Bigl\{
\int_{\Gamma}q_1^{-2}(x)
|z(x;u)|^2\,d\Gamma \Bigr\}^{\frac 1{2}}:=a,
$$
in which the equality holds at  $\tilde
h\in G_0$ and
\begin{equation}\label{lll}
\tilde h(\cdot)=\pm\frac{q_1^{-2}(\cdot) z(\cdot;u)
\left.\right|_{\Gamma}}{a}+h_0.
\end{equation}
Therefore,
$$
\inf_{c \in \mathbb C}\sup_{\tilde h \in G_0} \Bigl|
\int_{\Gamma}
\tilde h \overline{z(\cdot;u)}\, d\Gamma -c\Bigr|^2
$$
$$
=\inf_{c \in \mathbb C}\sup_{|y|\leq a}|y+\int_{\Gamma}\overline{z(\cdot;u)}h_0\,d\Gamma-c|^2=a^2
= \int_{\Gamma}q_1^{-2}(x)
|z(x;u)|^2\,d\Gamma
$$
at
\begin{equation}\label{sas}
c=\int_{\Gamma}\overline{z(\cdot;u)}h_0\,d\Gamma.
\end{equation}
Similarly,
\begin{gather*}
\sup_{\tilde \eta \in G_1} \mathbf E\Biggl|\sum_{i=1}^{N}\int_{\gamma_i}
\overline{u_i\n(x)}\tilde\eta_i\n(x)\,d\gamma_i+\sum_{i=1}^{N}\int_{\gamma_i}
\overline{u_i\m(x)}\tilde\eta_i\m(x)\,d\gamma_i\Biggr|^2\\
\leq\sum_{i=1}^{N} \int_{\gamma_i}(
r_i\n(x))^{-2}|u_i\n(x)|^2d\gamma_i +\sum_{i=1}^{N}
\int_{\gamma_i} ( r_i\m(x))^{-2}|u_i\m(x)|^2d\gamma_i.
\end{gather*}
It is easy to see that in this inequality, the equality holds when $\tilde \eta$ is a random vector-function
with the components
\begin{equation}\label{lllb}
\begin{split}
\tilde\eta_{i}\n(x)=\frac{\xi
\,(r\n_i(x))^{-2}u_i\n(x)}{\left[\sum_{i=1}^{N} \int_{\gamma_i}
(r_i\n(x))^{-2}|u_i\n(x)|^2d\gamma_i +\sum_{i=1}^{N}
\int_{\gamma_i} (r_i\m(x))^{-2}|u_i\m(x)|^2d\gamma_i
\right]^{1/2}},\\ \tilde\eta_{i}\m(x)=\frac{\xi
\,(r\m_i(x))^{-2}u_i\m(x)}{\left[\sum_{i=1}^{N} \int_{\gamma_i}(
r_i\n(x))^{-2}|u_i\n(x)|^2d\gamma_i +\sum_{i=1}^{N}
\int_{\gamma_i} (r_i\m(x))^{-2}|u_i\m(x)|^2d\gamma_i
\right]^{1/2}},
\end{split}
\end{equation}
where $\xi$ is a random value such that $\mathbf E\xi=0$ and $\mathbf E\xi^2=1.$
The latter facts yield
$$
\inf_{c \in \mathbb C} \sup_{\tilde h \in G_0,\, \tilde \eta \in G_1}
\mathbf E|l(\tilde\varphi)- \widehat {l(\tilde\varphi)}|^2 =
I(u_1\n, \ldots ,u_{m_1}\n,u_1\m, \ldots ,u_{m_2}\m),
$$
where functional $I$ is determined according to (\ref{17}) and the infimum with respect
to $c$ is attained at $c=\int_{\Gamma}\overline{z(\cdot;u)}h_0\,d\Gamma$. The lemma is proved.
\end{proof}

Solving the optimal control problem  (\ref{13})$-$(\ref{17}), we arrive at the following statement.
\begin{pred} \label{t4}
The minimax estimate of the value of functional $l(\varphi)$ has the form
\begin{equation} \label{2000}
\widehat{\widehat{l(\varphi)}}
=\sum_{i=1}^N\int_{\gamma_i}\left(\overline{\hat u_i\n(x)}
y_i\n(x)+\overline{\hat
u_i\m(x)}y_i\m(x)\right)\, d\gamma_i+\hat c,
\end{equation}
where
\begin{equation}\label{sas1}
\hat c=\int_{\Gamma}\overline{z}h_0\,d\Gamma,
\end{equation}
\begin{align}   \label{20}
\hat u_i\n(x) &=(r_i\n(x))^2
\int_{\gamma_i}\left[K_i\p(x,\xi)p(\xi)+
K_i\q(x,\xi)\frac{\partial p(\xi)} {\partial
\nu}\right]\,d\gamma_{i_\xi}, \\ \hat u_i\m(x) &=(r_i\m(x))^2
\int_{\gamma_i}\left[K_i\s(x,\xi)p(\xi)+
K_i\w(x,\xi)\frac{\partial p(\xi)} {\partial
\nu}\right]\,d\gamma_{i_\xi},\quad i= \overline{1,N},\nonumber
\end{align}
and function $p(x)$ is determined from the solution to the problem
\begin{equation} \label{r13}
z \in H^1_{\rm loc}(\Omega',\Delta), \end{equation}
\begin{equation}\label{r14}
-(\Delta+\bar k^2)z(x)=\chi_{\omega_0}(x)l_0(x)\quad \mbox {in} \quad \Omega',
\end{equation}
\begin{equation}\label{r15}
\frac{\partial z}{\partial\nu}=0 \quad  \mbox {on} \quad
\Gamma,
\end{equation}
\begin{equation}\label{r15''}
[z(x)]_{\gamma_i}=\int_{\gamma_i}\left[\overline{K_i\q(\xi,x)}\hat
u_i\n(\xi)+ \overline{K_i\w(\xi,x)}\hat u_i\m(\xi)\right]\,d\gamma_{i_{\xi}},
\end{equation}
\begin{equation*}
\left[\frac{\partial z(x)}{\partial\nu}\right]_{\gamma_i}=
\end{equation*}
\begin{equation}\label{r16''''}
= -\int_{\gamma_i}\left[ \overline{K_i\p(\xi,x)}\hat u_i\n(\xi)+
\overline{K_i\s(\xi,x)}\hat u_i\m(\xi)\right]\,d\gamma_{i_{\xi}} \,\, \mbox
{on}\,\, \gamma_i,\,\,i=\overline{1,N},
\end{equation}
\begin{equation}\label{r16'}
\frac{\partial z(x)}{\partial
r}+i\bar kz(x)=o(1/r),\,\,r=|x|,\,\,r\to \infty,
\end{equation}
\begin{equation}\label{re25}
p \in H^1_{\rm loc}(\Omega',\Delta),
\end{equation}
\begin{equation}\label{re26}
(\Delta+k^2)p(x)=0\quad\mbox{in}\quad \Omega',
\end{equation}
\begin{equation}\label{re27}
\frac{\partial p}{\partial\nu}=q_1^{-2}z \quad \mbox {on} \quad
\Gamma,
\end{equation}
\begin{equation}\label{re28'}
[p]_{\gamma_i}=0,\quad \left[\frac{\partial p}{\partial
\nu}\right]_{\gamma_i}=0,\quad i= \overline{1,N}.
\end{equation}
\begin{equation}\label{re28''}
\frac{\partial p(x)}{\partial
r}-ikp(x)=o(1/r),\,\,r=|x|,\,\,r\to \infty.
\end{equation}
where in (\ref{20}) $p(y)$ and $\frac{\partial
p(y)}{\partial \nu}$ denote the
values of the traces of functions
$p$ and boundary values of its conormal derivatives on different sides of surface $\gamma_i.$
Also, $p|_{\gamma_i}\,
\frac{\partial p(y)}{\partial \nu}|_{\gamma_i}\in
L^2(\gamma_i),$ $i= \overline{1,N}.$ Problem (\ref{20})$-$(\ref{re28''}) is uniquely solvable.

The  error $\sigma$ of the minimax estimation of $l(\varphi)$ is
given by the formula
\begin{equation} \label{34}
\sigma=[l(p)]^{1/2}= \left(\int_{\omega_0}
\overline{l_0(x)}p(x)\,dx \right) ^{1/2}.
\end{equation}
\end{pred}

Note that if we replace in (\ref{r15''}) and (\ref{r16''''})
 $\hat u_i\n(x)$ and $\hat u_i\m(x),$
$i=\overline{1,N},$ by their expressions in the right-hand sides of (\ref{20}), then these functions may be excluded from the equality system (\ref{r13})$-$(\ref{re28''}).
\begin{proof}
Let us show that $I(u)$ is a quadratic function on $H.$ Indeed,
since the solution $z(x;u)$ to problem (\ref{13})$-$(\ref{16})
can be represented as $z(x;u)=\tilde z(x;u)+z_0(x),$
\label{p.36} where $\tilde z(x;u)$ is the solution to this problem at $l_0(x)\equiv 0$ and $\tilde z_0(x)$ is the solution to the problem
\begin{equation} \label{13x}
z_0 \in H^1_{\rm loc}(\Omega',\Delta),
\end{equation}
\begin{equation}\label{14x}
-(\Delta+\bar k^2)z_0(x)=\chi_{\omega_0}(x)l_0(x)\quad \mbox {in} \quad
\Omega',
\end{equation}
\begin{equation}\label{15x} \frac{\partial
z_0}{\partial\nu}=0 \quad \mbox {on} \quad \Gamma,
\end{equation}
\begin{equation*}
[z_0(x)]_{\gamma_i}=0 \quad \mbox {on}\quad \gamma_i,\quad i=
\overline{1,N},
\end{equation*}
\begin{equation}\label{16x}
\left[\frac{\partial z_0(x)}{\partial\nu}\right]_{\gamma_i}=0
\quad \mbox {on}\quad \gamma_i,\quad i=
\overline{1,N},
\end{equation}
\begin{equation}\label{16x'}
\frac{\partial z_0(x)}{\partial
r}+i\bar kz_0(x)=o(1/r),\,\,r=|x|,\,\,r\to \infty,
\end{equation}
functional $I(u)$ can be represented as \label{D}
\begin{equation*}
I(u)=\tilde I(u)+L(u)
+\int_{\Gamma}q_1^{-2} |z_0|^2 \,d\Gamma,
\end{equation*}
where
\begin{equation*} \tilde I(u)=
\int_{\Gamma}q_1^{-2}(x)|\tilde z^2(x;u)| \,d\Gamma
\end{equation*}
\begin{equation*} +\sum_{i=1}^{N}\int_{\gamma_i}
 (r_i\n(x))^{-2}|u_i\n(x)|^2d\gamma_i
 +\sum_{i=1}^{N}\int_{\gamma_i}
 (r_i\m(x))^{-2}|u_i\m(x)|^2d\gamma_i,
\end{equation*}
$$
L(u)=
2\mbox{\rm\, Re}\int_{\Gamma}q_1^{-2}(x) \tilde z(x;u)\overline{z_0(x)} \,d\Gamma.
$$
From inequality (\ref{yr3}) and our assumptions concerning operators of the form (\ref{ya1}), we deduce, taking into account that
$\|\cdot\|_{H^{-1/2}(\gamma_i)}\leq c\|\cdot\|_{L^2(\gamma_i)}$
($c=\mbox{const}>0$),
the inequality
\label{pi3}
\begin{equation*}
\|\tilde z(\cdot;u)\|_{H^1(\Omega'_R)}\leq \Bigl(\sum_{i=1}^N
\Bigl\|\int_{\gamma_i}\left[\overline{K_i\q(\xi,x)}u_i\n(\xi)+
\overline{K_i\w(\xi,x)}u_i\m(\xi)\right]\,d\gamma_{i_{\xi}}
\Bigr\|^2_{H_{0\,0}^{1/2}(\gamma_i)}\Bigr)^{1/2}
\end{equation*}
\begin{equation}\label{ba18}
\end{equation}
\begin{equation*}
+\Bigl(\sum_{i=1}^N\Bigl\|\int_{\gamma_i}\left[\overline{K_i\p(\xi,x)}u_i\n(\xi)+
\overline{K_i\s(\xi,x)}u_i\m(\xi)\right]\,
d\gamma_{i_{\xi}}\Bigr\|^2_{H^{-1/2}(\gamma_i)}\Bigr)^{1/2} \leq
c_1\|u\|_H,
\end{equation*}
where $c_1=\mbox{const}$ that does not depend on $u.$
Taking into account (\ref{ba18}) and the trace theorem from
\cite{BIBLlio}, we see that
$u \to \gamma_D\tilde z(\cdot;u)$ is a bounded linear operator that maps Hilbert space $H$ in
$H^{1/2)}(\Gamma).$ From the latter statement, it follows that $\tilde I(u)$ is a quadratic form
which corresponds to a semi-linear continuous Hermitian form
$$
\pi(u,v):=
 \int_{\Gamma}q_1^{-2}(x) \tilde
z(x;u)\overline{\tilde z(x;v)} \,d\Gamma
$$
$$
+\sum_{i=1}^{N}\int_{\gamma_i}
 (r_i\n(x))^{-2}u_i\n(x)\overline{v_i\n(x)}d\gamma_i
 +\sum_{i=1}^{N}\int_{\gamma_i}
 (r_i\m(x))^{-2}u_i\m(x)\overline{v_i\m(x)}d\gamma_i,
$$
and $L(u)$ a linear continuous functional defined on $H.$
Moreover, since
\begin{equation*}
\tilde I(u)\geq \sum_{i=1}^{N}\!\int_{\gamma_i}
(r_i\n(x))^{-2}|u_i\n(x)|^2d\gamma_i
\end{equation*}
\begin{equation*}
+\sum_{i=1}^{N}\int_{\gamma_i}
(r_i\m(x))^{-2}|u_i\m(x)|^2d\gamma_i \geq c\|u\|_H \quad \forall
u\in H,\,\, \mbox{c=const},
\end{equation*}
we obtain, using Remark 1.1 to Theorem 1.1 from \cite{BIBLlio},
that there exists \label{d2} one and only one element $\hat
u=(\hat u_1\n, \ldots ,\hat u_N\n,\hat u_1\m, \ldots ,\hat
u_N\m)\in H$ such that
що
\begin{equation*}
I(\hat u)=\inf_{u \in H}I(u).
\end{equation*}
Therefore, for any fixed $v \in H$ and $\tau \in \mathbb R^1$, the function
$s(\tau):=I(\hat u+\tau v)$ has only one minimum point $\tau
=0,$ so that
$$
\frac {d}{d\tau} I(\hat u+\tau v)\left. \right|_{\tau=0}=0.
$$
This yields
$$
0=\frac {1}{2} \frac {d}{d\tau} I(\hat u+\tau v)
\left.\right|_{\tau=0}
=\lim_{\tau\to 0}\frac 1{2\tau}(I(\hat u+\tau v)-I(\hat u))
$$
$$
=\lim_{\tau\to 0}\frac 1{2\tau}\int_{\Gamma}
q^{-2}(x)\left(z(x;\hat u+\tau v)\overline{z(x;\hat u+\tau v)}-z(x;\hat u)\overline{z(x;\hat u)}\right)\,d\Gamma
$$
\begin{multline*}
+\lim_{\tau\to 0}\frac 1{2\tau}
\sum_{i=1}^{N}\int_{\gamma_i}(r_i\n(x))^{-2}\left[\left(\hat
u_i\n(x)+\tau v_i\n(x)\right)\overline{\left(\hat
u_i\n(x)+\tau v_i\n(x)\right)}\right.\\\left.-\hat u_i\n(x)\overline{\hat u_i\n(x)}\right]d\gamma_i
\end{multline*}
\begin{multline*}\label{jj}
+\lim_{\tau\to 0}\frac 1{2\tau}
\sum_{i=1}^{N}\int_{\gamma_i}(r_i\m(x))^{-2}\left[\left(\hat
u_i\m(x)+\tau v_i\m(x)\right)\overline{\left(\hat
u_i\m(x)+\tau v_i\m(x)\right)}\right.\\\left.-\hat u_i\m(x)\overline{\hat u_i\m(x)}\right]d\gamma_i.
\end{multline*}
Calculate the first limit in the right-hand side of the last relationship.
Taking into notice the notation for $\tilde z(x;v)$
and the equality $z(x;\hat u+\tau v)=z(x;\hat
u)+\tau \tilde z(x;v)$, we have
$$
\lim_{\tau\to 0}\frac 1{2\tau}\int_{\Gamma}
q^{-2}(x)\left(z(x;\hat u+\tau v)\overline{z(x;\hat u+\tau v)}-z(x;\hat u)\overline{z(x;\hat u)}\right)\,d\Gamma
$$
$$
=\lim_{\tau\to 0}\frac 1{2\tau}\int_{\Gamma} q^{-2}(x)
\left[2\tau\mbox{\rm Re\,}(
 z(\cdot;\hat u)\overline{\tilde z(\cdot;v)})+O(\tau^2)\right]\,d\Gamma
$$
$$
=\mbox{\rm Re\,}\int_{\Gamma}q_1^{-2}(x) z(x;\hat u)\overline{\tilde z(x;v)}d\Gamma.
$$
Performing similar calculations for the remaining limits we obtain
$$
0=\mbox{\rm Re\,}\int_{\Gamma}q_1^{-2}(x) z(x;\hat u)\overline{\tilde z(x;v)}d\Gamma
+\sum_{i=1}^{N}\mbox{\rm Re\,}\int_{\gamma_i} (r_i\n(x))^{-2}\hat
u_i\n(x)\overline{v_i\n(x)}d\gamma_i
$$
\begin{equation} \label{29}
+\sum_{i=1}^{N}\mbox{\rm Re\,}\int_{\gamma_i}
(r_i\m(x))^{-2}\hat u_i\m(x)\overline{v_i\m(x)}d\gamma_i.
\end{equation}\\[-25pt]
On the other side, for any fixed $v \in H$ and $\tau \in \mathbb R^1$, the function
$s_1(\tau_1):=I(\hat u+i\tau_1 v)$ has the unique minimum point at $\tau_1
=0,$ so that
$$
\frac {d}{d\tau_1} I(\hat u+i\tau_1 v)\left. \right|_{\tau_1=0}=0,
$$
which yields
$$
0=\mbox{\rm Im\,}\int_{\Gamma}q_1^{-2}(x) z(x;\hat u)\overline{\tilde z(x;v)}d\Gamma
+\sum_{i=1}^{N}\mbox{\rm Im\,}\int_{\gamma_i} (r_i\n(x))^{-2}\hat
u_i\n(x)\overline{v_i\n(x)}d\gamma_i
$$
\begin{equation} \label{29'}
+\sum_{i=1}^{N}\mbox{\rm Im\,}\int_{\gamma_i}
(r_i\m(x))^{-2}\hat u_i\m(x)\overline{v_i\m(x)}d\gamma_i.
\end{equation}\\[-25pt]
From the latter, in line with \eqref{29} and \eqref{29'}, it follows
$$
0=\int_{\Gamma}q_1^{-2}(x) z(x;\hat u)\overline{\tilde z(x;v)}d\Gamma
+\sum_{i=1}^{N}\int_{\gamma_i} (r_i\n(x))^{-2}\hat
u_i\n(x)\overline{v_i\n(x)}d\gamma_i
$$
\begin{equation} \label{29''}
+\sum_{i=1}^{N}\int_{\gamma_i}
(r_i\m(x))^{-2}\hat u_i\m(x)\overline{v_i\m(x)}d\gamma_i.
\end{equation}\\[-25pt]

Introduce function $p(x)$ as the unique solution to the problem\\[-25pt]
\begin{equation}
\label{e25} p \in H^1_{\rm loc}(\Omega',\Delta),
\end{equation}
\begin{equation}\label{e26}
(\Delta+k^2)p(x)=0\quad\mbox{in}\quad \Omega',
\end{equation}
\begin{equation}\label{e27}
\frac{\partial p}{\partial\nu}=q_1^{-2} z(\cdot;\hat u) \quad \mbox {on} \quad
\Gamma,
\end{equation}
\begin{equation} \label{e28}
[p]_{\gamma_i}=0,\quad \left[\frac{\partial p}{\partial
\nu}\right]_{\gamma_i}=0,\quad i= \overline{1,N},
\end{equation}
\begin{equation}\label{e28'}
\frac{\partial p(x)}{\partial
r}-ikp(x)=o(1/r),\,\,r=|x|,\,\,r\to \infty.
\end{equation}
Transform the first term in the right-hand side of (\ref{29})
using equalities (\ref{e25})$-$(\ref{e28'}) and applying the second Green formula in domains $\Omega_i,$ $i=1,N,$ and
$\tilde \Omega_R$ to functions $p$ and $\overline{\tilde z(\cdot;v)}$.  We have
$$
0=\int_{\tilde\Omega_R}-(\Delta+k^2)p(x) \overline{\tilde
z(x;v)}dx+\sum_{i=1}^{N}\int_{\Omega_i}-(\Delta+k^2)p(x)\overline{ \tilde z(x;v)}dx
$$
$$
=\int_{\tilde\Omega_R} -\overline{(\Delta+\bar k^2)\tilde z(x;v)}p(x)\,dx
-\int_{\Gamma} \overline{\tilde z(\cdot;v)}\frac {\partial p}{\partial
\nu_{A}}\,d\Gamma -\Sigma_R(\tilde z(\cdot;v),p)
$$
$$
+\sum_{i=1}^{N}\int_{\Omega_i} -\overline{(\Delta+\bar k^2)\tilde
z(x;v)}p(x)\,dx -\sum_{i=1}^{N}\int_{\hat\gamma_i}\left(\left(\overline{\frac
{\partial \tilde z(\cdot;v)}{\partial \nu_{A^*}}}\right)_-p(x) - \overline{\tilde
z_-(\cdot;v)} \frac {\partial p}{\partial \nu_{A}}\,\right)
d\hat\gamma_i
$$
$$
+\sum_{i=1}^{N}\int_{\hat\gamma_i}\left(\left(\overline{ \frac {\partial
\tilde z(\cdot;v)}{\partial \nu_{A^*}}}\right)_+ p -\overline{\tilde z_+(\cdot;v)}
\frac {\partial p(x)}{\partial \nu_{A}}\right) d\hat\gamma_i
$$
$$
=\int_{\Omega_R'} -\overline{(\Delta+\bar k^2)\tilde z(x;v)}p(x)\,dx
-\int_{\Gamma} \overline{\tilde z(\cdot;v)} \frac {\partial p}{\partial
\nu}\,d\Gamma-\Sigma_R(\tilde z(\cdot;v),p)
$$
$$
-\sum_{i=1}^{N}\int_{\gamma_i} \left(\overline{\tilde z_+(\cdot;v)} -\overline{\tilde
z_-(\cdot;v)}\right) \frac {\partial p}{\partial
\nu}\,d\hat\gamma_i
$$
$$
+\sum_{i=1}^{N}\int_{\gamma_i} \left(\left(\overline{\frac {\partial \tilde
z(\cdot;v) }{\partial \nu_{A^*}}}\right)_+ -\left(\overline{\frac{\partial \tilde
z(\cdot;v)}{\partial \nu_{A^*}}}\right)_- \right)p\,d\gamma_i
$$
$$
=-\int_{\Gamma} \overline{\tilde z(\cdot;v)} q_1^{-2}(x) z(x;\hat u)\,d\Gamma
-\Sigma_R(\tilde z(\cdot;v),p)
$$
$$
-\sum_{i=1}^{N}\int_{\gamma_i}p(x) \int_{\gamma_i}\left[\right.
K_i\p(\xi,x)\overline{v_i\n(\xi)} +K_i\s(\xi,x)\overline{v_i\m(\xi)}\left.\right]
\,d\gamma_{i_{\xi}}\,d\gamma_{i_x}
$$
$$
-\sum_{i=1}^{N}\int_{\gamma_i}\frac {\partial p(x)}{\partial
\nu} \int_{\gamma_i}\left[\right. K_i\q(\xi,x)\overline{ v_i\n(\xi)}
+K_i\w(\xi,x)\overline{ v_i\m(\xi)}\left.\right]
\,d\gamma_{i_{\xi}}\,d\gamma_{i_x}.
$$
Calculating the limit as $R\to\infty$ and taking into account that
$\Sigma_R(\tilde z(\cdot;v),p)=o(1),$
we obtain
$$
0=\int_{\Gamma} \overline{\tilde z(\cdot;v)} q_1^{-2}(x) z(x;\hat u)\,d\Gamma
$$
$$
\sum_{i=1}^{N}\int_{\gamma_i}p(x) \int_{\gamma_i}\left[\right.
K_i\p(\xi,x)\overline{v_i\n(\xi)} +K_i\s(\xi,x)\overline{v_i\m(\xi)}\left.\right]
\,d\gamma_{i_{\xi}}\,d\gamma_{i_x}
$$
$$
\sum_{i=1}^{N}\int_{\gamma_i}\frac {\partial p(x)}{\partial
\nu} \int_{\gamma_i}\left[\right. K_i\q(\xi,x)\overline{ v_i\n(\xi)}
+K_i\w(\xi,x)\overline{ v_i\m(\xi)}\left.\right]
\,d\gamma_{i_{\xi}}\,d\gamma_{i_x}.
$$
Next, by virtue of (\ref{29''}),
$$
\sum_{i=1}^{N}\int_{\gamma_i} (r_i\n(x))^{-2}\hat
u_i\n(x)\overline{v_i\n(x)}d\gamma_i +\sum_{i=1}^{N}\int_{\gamma_i}
(r_i\m(x))^{-2}\hat u_i\m(x)\overline{v_i\m(x)}d\gamma_i
$$
$$
=\sum_{i=1}^{N}\int_{\gamma_i}\overline{v_i\n(x)} \left[K_i\p(x,\xi)p(\xi)
+K_i\q(x,\xi)\frac {\partial p(\xi)}{\partial \nu}\right]
\,d\gamma_{i_{\xi}}\,d\gamma_{i_x}
$$
$$
+\sum_{i=1}^{N}\int_{\gamma_i}\overline{v_i\m(x)} \left[ K_i\s(x,\xi)p(\xi)
+K_i\w(x,\xi)\frac {\partial p(\xi)}{\partial \nu}\right]
\,d\gamma_{i_{\xi}}\,d\gamma_{i_x}.
$$
Rewrite the last equality in the form
\begin{multline*}
\sum_{i=1}^{N}\int_{\gamma_i} (r_i\n(x))^{-2}\!\left[\hat
u_i\n(x)-(r_i\n(x))^2\int_{\gamma_i}
\left(K_i\p(x,\xi)p(\xi)\right.\right. +\\
\phantom{(r_i\n(x))^2}\left.\left. +K_i\q(x,\xi)\frac {\partial
p(\xi)}{\partial \nu}\right) \,d\gamma_{i_{\xi}}\right]
\overline{v_i\n(x)}d\gamma_{i_x}
\end{multline*}\\[-80pt]\\
\begin{equation} \label{hh}\end{equation}\\[-80pt]
\begin{multline*}
\!\!\!\!\!\!\!+\sum_{i=1}^{N}\int_{\gamma_i}(r_i\m(x))^{-2}\left[\hat
u_i\m(x)-(r_i\m(x))^2\int_{\gamma_i}
\left(K_i\s(x,\xi)p(\xi)\right.\right. +\\
\left.\left.+K_i\w(x,\xi)\frac {\partial p(\xi)}{\partial
\nu}\right) \,d\gamma_{i_{\xi}}\right] \overline{v_i\m(x)}d\gamma_{i_x}=0.
\end{multline*}
Setting in (\ref{hh})
$$
v_i\n(x)=\hat u_i\n(x)-(r_i\n(x))^2\int_{\gamma_i}
\left(K_i\p(x,\xi)p(\xi) +K_i\q(x,\xi)\frac {\partial
p(\xi)}{\partial \nu}\right) \,d\gamma_{i_{\xi}},
$$
$$
v_i\m(x)=\hat u_i\m(x)-(r_i\m(x))^2\int_{\gamma_i}
\left(K_i\s(x,\xi)p(\xi) +
K_i\w(x,\xi)\frac {\partial p(\xi)}{\partial
\nu}\right) \,d\gamma_{i_{\xi}},
$$
$i=\overline{1,N},$ we find
\begin{multline*}
\phantom{+}\sum_{i=1}^{N}\int_{\gamma_i}
(r_i\n(x))^{-2}\!\left|\hat u_i\n(x)-(r_i\n(x))^2\int_{\gamma_i}
\left(K_i\p(x,\xi)p(\xi)\right.\right. +\\
\left.\left.+K_i\q(x,\xi)\frac {\partial p(\xi)}{\partial
\nu}\right) \,d\gamma_{i_{\xi}}\right|^2d\gamma_{i_x}
\end{multline*}\\[-40pt]
\begin{multline*}
+\sum_{i=1}^{N}\int_{0\gamma_i}(r_i\m(x))^{-2}\!\left|\hat
u_i\m(x)-(r_i\m(x))^2\int_{\gamma_i}
\left(K_i\s(x,\xi)p(\xi)\right.\right. +\\ \left.\left.
+K_i\w(x,\xi)\frac {\partial p(\xi)}{\partial \nu}\right)
\,d\gamma_{i_{\xi}}\right|^2 d\gamma_{i_x}=0,
\end{multline*}
and consequently,
\begin{equation*}
\begin{split}
\hat u_i\n(x) &=(r_i\n(x))^2
\int_{\gamma_i}\left[K_i\p(x,\xi)p(\xi)+
K_i\q(x,\xi)\frac{\partial p(\xi)} {\partial
\nu}\right]\,d\gamma_{i_\xi}, \\ \hat u_i\m(x) &=(r_i\m(x))^2
\int_{\gamma_i}\left[K_i\s(x,\xi)p(\xi)+
K_i\w(x,\xi)\frac{\partial p(\xi)} {\partial
\nu}\right]\,d\gamma_{i_\xi},\quad i= \overline{1,N}.
\end{split}
\end{equation*}
Substituting these quantities to (\ref{13}) and (\ref{16}), setting $z(x)=z(x;\hat u)$,
and taking into account (\ref{e25})$-$(\ref{e28'}), we arrive at problem
(\ref{r13})$-$(\ref{re28'}); the unique solvability of this problem follows from the fact that functional (\ref{17})
has the unique minimum point $\hat u$.

Now let us establish the validity of
(\ref{34}).
 Substituting expressions (\ref{20}) to (\ref{17}),
we obtain
$$
\sigma^2=I(\hat u_1\n, \ldots ,\hat u_{N}\n,\hat u_1\m, \ldots
,\hat u_{N}\m)=
$$
$$
=\int_{\Gamma} q_1^{-2}(x)
|z(x)|^2\,d\Gamma
$$
$$
+\sum_{i=1}^{N}\int_{\gamma_i}
 (r_i\n(x))^{-2}|\hat u_i\n(x)|^2d\gamma_i+\sum_{i=1}^{N}\int_{\gamma_i}
(r_i\m(x))^{-2}|\hat u_i\m(x)|^2d\gamma_i
$$
$$
=\int_{\Gamma}
q_1^{-2}(x) |z(x)|^2\,d\Gamma
$$
$$
+\sum_{i=1}^{N}\int_{\gamma_i}(r_i\n(x))^2
\Bigl|\int_{\gamma_i}\left[K_i\p(x,\xi)p(\xi)+
K_i\q(x,\xi)\frac{\partial p(\xi)} {\partial
\nu}\right]\,d\gamma_{i_\xi}\Bigr|^2 d\gamma_{i_x}
$$
\begin{equation}\label{35}
+\sum_{i=1}^{N}\int_{\gamma_i}(r_i\m(x))^2
\Bigl|\int_{\gamma_i}\left[K_i\s(x,\xi)p(\xi)+
K_i\w(x,\xi)\frac{\partial p(\xi)} {\partial
\nu}\right]\,d\gamma_{i_\xi}\Bigr|^2 d\gamma_{i_x}.
\end{equation}
Next, using relationships (\ref{20})$-$(\ref{re28'}),
we have \label{lla}
$$
0=\int_{\tilde\Omega_R}-(\Delta+k^2)p(x)\overline{z(x)}\,dx
+\sum_{i=1}^{N}\int_{\Omega_i}-(\Delta+k^2)p(x)\overline{z(x)}\,dx
$$
$$
=\int_{\tilde\Omega_R}-\overline{(\Delta+\bar k^2)z(x)}p(x)\,dx -\int_{\Gamma}
\overline{z(x)}\frac {\partial p(x)}{\partial \nu}\,d\Gamma-\Sigma_R(z,p)
$$
$$
-\sum_{i=1}^{N}\int_{\hat\gamma_i}\left(\left(\overline{ \frac{\partial
z}{\partial \nu}}\right)_- p-\overline{z}_-\frac {\partial p}{\partial
\nu}\right)\, d\hat\gamma_i
$$
$$
+\sum_{i=1}^{N}\int_{\Omega_i} -\overline{(\Delta+\bar k^2)z(x)}p(x)\,dx
+\sum_{i=1}^{N}\int_{\hat\gamma_i}\left(\left(\overline{ \frac{\partial
z}{\partial \nu}}\right)_+ p-\overline{z}_+\frac {\partial p}{\partial
\nu}\,\right)\, d\hat\gamma_i
$$
$$
=\int_{\Omega'_R} -\overline{(\Delta+\bar k^2)z(x)}p(x)\,dx -\int_{\Gamma}
\overline{z}\frac {\partial p} {\partial \nu}\,d\Gamma-\Sigma_R(z,p)
$$
$$
-\sum_{i=1}^{N}\int_{\gamma_i} (\overline{z}_+ -\overline{z}_-) \frac {\partial
p}{\partial \nu}\,d\gamma_i +\sum_{i=1}^{N}\int_{\gamma_i}
\left[\left(\overline{\frac {\partial z}{\partial \nu}}\right)_+
-\left(\overline{\frac{\partial z}{\partial \nu}}\right)_-
\right]p\,d\gamma_i
$$
$$
=\int_{\Omega'_R}\chi_{\omega_0}(x) \overline{l_0(x)}p(x)\,dx -\int_{\Gamma} q_1^{-2}|z|^2\,d\Gamma-\Sigma_R(z,p)
$$
\begin{multline*}
-\sum_{i=1}^{N}\int_{\gamma_i}\frac {\partial p(x)}{\partial
\nu} \int_{\gamma_i}K_i\q (\xi,x)(r_i\n(\xi))^2
\int_{\gamma_i}\overline{\Bigl[\Bigr.p(y)}\times\\ \times\overline{ K_i\p (\xi,y)
+\frac{\partial p(y)} {\partial \nu}K_i\q (\xi,y)\Bigr.\Bigr]}
\,d\gamma_{i_y}\,d\gamma_{i_{\xi}}\, d\gamma_{i_x}
\end{multline*}
\begin{multline*}
-\sum_{i=1}^{N}\int_{\gamma_i}\frac {\partial p(x)}{\partial
\nu} \int_{\gamma_i}K_i\w (\xi,x)(r_i\m(\xi))^2
\int_{\gamma_i}\overline{\Bigl[\Bigr.p(y)}\times\\ \times\overline{ K_i\s (\xi,y)
+\frac{\partial p(y)} {\partial \nu}K_i\w
(\xi,y)\Bigr.\Bigr]}\,d\gamma_{i_y}
\,d\gamma_{i_{\xi}}\,d\gamma_{i_x}
\end{multline*}
\begin{multline*}
-\sum_{i=1}^{N}\int_{\gamma_i}p(x) \int_{\gamma_i}K_i\p
(\xi,x)(r_i\n(\xi))^2 \int_{\gamma_i}\overline{\Bigl[\Bigr.p(y)}\times\\
\times \overline{K_i\p (\xi,y) +\frac{\partial p(y)} {\partial \nu}K_i\q
(\xi,y)\Bigr.\Bigr]} \,d\gamma_{i_y}\,d\gamma_{i_{\xi}}\,
d\gamma_{i_x}
\end{multline*}
\begin{multline*}
-\sum_{i=1}^{N}\int_{\gamma_i} p(x) \int_{\gamma_i}K_i\s
(\xi,x)(r_i\m(\xi))^2 \int_{\gamma_i}\overline{\Bigl[\Bigr.p(y)}\times\\
\times \overline{K_i\s (\xi,y) +\frac{\partial p(y)} {\partial \nu}K_i\w
(\xi,y)\Bigr.\Bigr]}\,d\gamma_{i_y}
\,d\gamma_{i_{\xi}}\,d\gamma_{i_x}
\end{multline*}
$$
=\int_{\omega_0} \overline{l_0(x)}p(x)\,dx -\int_{\Gamma} q_1^{-2}|z|^2\,d\Gamma
-\Sigma_R(z,p)
$$
$$
-\sum_{i=1}^{N}\int_{\gamma_i}(r_i\n(\xi))^2\,d\gamma_{i_{\xi}}
\int_{\gamma_i}\overline{\Bigl[\Bigr.p(y) K_i\p (\xi,y)+\frac{\partial p(y)}
{\partial \nu}K_i\q (\xi,y)\Bigr.\Bigr]} \,d\gamma_{i_y}\times
$$
$$
\times \int_{\gamma_i}\Bigl[\Bigr.p(x) K_i\p
(\xi,x)+\frac{\partial p(x)} {\partial \nu}K_i\q
(\xi,x)\Bigr.\Bigr] \,d\gamma_{i_x}
$$
$$
-\sum_{i=1}^{N}\int_{\gamma_i}(r_i\m(\xi))^2\,d\gamma_{i_{\xi}}
\int_{\gamma_i}\overline{\Bigl[\Bigr.p(y) K_i\s (\xi,y)+\frac{\partial p(y)}
{\partial \nu}K_i\w (\xi,y)\Bigr.\Bigr]} \,d\gamma_{i_y}\times
$$
$$\label{lla1}
\times \int_{\gamma_i}\Bigl[\Bigr.p(x) K_i\s
(\xi,x)+\frac{\partial p(x)} {\partial \nu}K_i\w
(\xi,x)\Bigr.\Bigr] \,d\gamma_{i_x}
$$
$$
=\int_{\omega_0} \overline{l_0(x)}p(x)\,dx -\int_{\Gamma}
q_1^{-2}|z|^2\,d\Gamma-\Sigma_R(z,p)
$$\\[-35pt]
\begin{align}\label{36}
\!\!-\!\sum_{i=1}^{N}\int_{\gamma_i}\!(r_i\n(x))^2 &
\Bigl|\int_{\gamma_i}\left[K_i\p(x,\xi)p(\xi)\!+\!
K_i\q(x,\xi)\frac{\partial p(\xi)} {\partial
\nu}\right]d\gamma_{i_\xi}\Bigr|^2 d\gamma_{i_x} \nonumber\\
\!\!-\!\sum_{i=1}^{N}\int_{\gamma_i}\!(r_i\m(x))^2 &
\Bigl|\int_{\gamma_i}\left[K_i\s(x,\xi)p(\xi)\!+\!
K_i\w(x,\xi)\frac{\partial p(\xi)} {\partial
\nu}\right]d\gamma_{i_\xi}\Bigr|^2 d\gamma_{i_x}.
\end{align}
Equality (\ref{34}) follows now from two relationships
(\ref{35}) and (\ref{36}) and the fact that $\Sigma_R(z,p)=o(1)$ when $R\to\infty.$
\end{proof}

An alternative representation for the minimax estimate in terms of the solution to a system of integro-differential equations
is given in the next theorem. This solution is independent of  the specific form of functional  (\ref{9}).
\begin{pred}\label{t6}
The minimax estimate of (\ref{9}) has the form
\begin{equation} \label{37}
\widehat{\widehat {l(\varphi)}}= l(\Hat \varphi)= \int_{\omega_0}\overline{
l_0(x)} \hat \varphi (x)\, dx,
\end{equation}
where function $\hat \varphi (x)$ is determined from the solution to the problem
(\ref{38})$-$(\ref{45'}):
\begin{equation} \label{38}
\hat p \in L^2(\Sigma,H_{\rm
loc}^1(\Omega',\Delta))\footnote{Here
$L^2(\Sigma,H^1_{\rm loc}(\Omega',\Delta))$ denotes a
class of functions the belong to
$L^2(\Sigma,H^1(\Omega'_R,\Delta))$ for any $R>0$.},
\end{equation}
\begin{equation} \label{39}
(\Delta+\bar k^2) \hat p(x,\omega)=0 \quad \mbox{in}\quad \Omega',
\end{equation}
\begin{equation}\label{40}
\frac{\partial\hat p(\cdot,\omega)}{\partial \nu}=0 \quad \mbox{on} \quad
\Gamma,
\end{equation}
\begin{multline*}
[\hat p(x,\omega))]_{\gamma_i}
=-\int_{\gamma_i}\overline{K_i\q(\xi,x)}\left[(r_i\n(\xi))^2y_i\n(\xi,\omega)-\hat
v_i\n(\xi,\omega)\right]d\gamma_{i_{\xi}} -\\
-\int_{\gamma_i}\overline{K_i\w(\xi,x)}\left[(r_i\m(\xi))^2y_i\m(\xi,\omega)-\hat
v_i\m(\xi,\omega)\right]d\gamma_{i_{\xi}},
\end{multline*}\\[-40pt]
\begin{multline*}
\left[\frac{\partial \hat p(x,\omega))}{\partial
\nu}\right]_{\gamma_i}
=\int_{\gamma_i}\overline{K_i\p(\xi,x)}\left[(r_i\n(\xi))^2y_i\n(\xi,\omega)-\hat
v_i\n(\xi,\omega)\right]d\gamma_{i_{\xi}}+\\
+\int_{\gamma_i}\overline{K_i\s(\xi,x)}\left[(r_i\m(\xi))^2y_i\m(\xi,\omega)-\hat
v_i\m(\xi,\omega)\right]d\gamma_{i_{\xi}}
\end{multline*}
\begin{equation}\label{41}
\mbox{on} \quad \gamma_i,\quad i=\overline{1,N},
\end{equation}
\begin{equation}\label{41'}
\frac{\partial \hat p(x,\omega)}{\partial
r}+i\bar k\hat p(x,\omega)=o(1/r),\,\,r=|x|,\,\,r\to \infty,
\end{equation}
\begin{equation} \label{42}
\hat \varphi \in L^2(\Sigma,H_{\rm loc}^1(\Omega',\Delta)),
\end{equation}
\begin{equation}\label{43}
(\Delta+k^2) \hat \varphi(x,\omega)=0 \quad \mbox{in} \quad \Omega',
\end{equation}
\begin{equation}\label{44}
\frac{\partial \hat \varphi(\cdot,\omega)}{\partial\nu}
=q_1^{-2}\hat p(\cdot,\omega) +h_0\quad
\mbox {on} \quad \Gamma,
\end{equation}
\begin{equation} \label{45}
[\hat \varphi(\cdot,\omega)]_{\gamma_i}=0,\quad \left[\frac{\partial \hat
\varphi(\cdot,\omega)} {\partial\nu}\right]_{\gamma_i}=0, \quad i=
\overline{1,N},
\end{equation}
\begin{equation}\label{45'}
\frac{\partial \hat \varphi(x,\omega)}{\partial
r}-ik\hat \varphi(x,\omega)=o(1/r),\,\,r=|x|,\,\,r\to \infty,
\end{equation}
where
\begin{equation}\label{45''}
\hat v_i\n(\xi,\omega)=(r_i\n(\xi))^2
\int_{\gamma_i}\left[K_i\p(\xi,\eta)\hat \varphi(\eta,\omega)+
K_i\q(\xi,\eta)\frac{\partial \hat \varphi(\eta,\omega)} {\partial
\nu}\right]\,d\gamma_{i_\eta},
\end{equation}
\begin{equation}\label{45'''}
\hat v_i\m(\xi,\omega)=(r_i\m(\xi))^2
\int_{\gamma_i}\left[K_i\s(\xi,\eta)\hat \varphi(\eta,\omega)+
K_i\w(\xi,\eta)\frac{\partial \hat \varphi(\eta,\omega)} {\partial
\nu}\right]\,d\gamma_{i_\eta},
\end{equation}
and the right-hand sides in (\ref{41}) are considered for every
realization of random functions $y_i\n(\xi)=y_i\n(\xi,\omega)$ and
$y_i\m(\xi)=y_i\m(\xi,\omega)$ which belong with probability $1$
to the space $L^2(\gamma_i),$ $i=\overline{1,N}.$ Problem
(\ref{38})$-$(\ref{45'}) is uniquely solvable.
\end{pred}
\begin{proof} The proof of this theorem is similar to the proof of Theorems
1.2 and
2.4.
\end{proof}

 \begin{predlllll}{\it Function $\hat \varphi (x,\omega)$ which is determined from the solution to problem (\ref{38})$-$(\ref{45'}) can be taken as a good estimate of the unknown solution $\varphi (x)$ to the initial Neumann problem
(\ref{1})$-$(\ref{3})} (see Remark 1 on p. \pageref{remark}).
\end{predlllll}

Set
\begin{multline} \label{3r1}
\tilde K_i\p(x,\eta)=\int_{\gamma_i}\left[K_i\p(\xi,\eta)
\overline{K_i\q(\xi,x)}(r_i\n(\xi))^2\right.\\ \left.
+K_i\s(\xi,\eta)
\overline{K_i\w(\xi,x)}(r_i\m(\xi))^2 \right]\,d\gamma_{i_{\eta}},
\end{multline}
\begin{multline} \label{3r2}
\tilde K_i\q(x,\eta)=\int_{\gamma_i}\left[K_i\q(\xi,\eta)
\overline{K_i\q(\xi,x)}(r_i\n(\xi))^2
\right.\\ \left.+K_i\w(\xi,\eta)
\overline{K_i\w(\xi,x)}(r_i\m(\xi))^2 \right]\,d\gamma_{i_{\eta}},
\end{multline}
\begin{multline} \label{3r3}
\tilde K_i\s(x,\eta)=-\int_{\gamma_i}\left[K_i\p(\xi,\eta)
\overline{K_i\p(\xi,x)}(r_i\n(\xi))^2
\right.\\ \left.+K_i\s(\xi,\eta)
\overline{K_i\s(\xi,x)}(r_i\m(\xi))^2 \right]\,d\gamma_{i_{\eta}},
\end{multline}
\begin{multline} \label{3r4}
\tilde K_i\w(x,\eta)=-\int_{\gamma_i}\left[K_i\q(\xi,\eta)
\overline{K_i\p(\xi,x)}(r_i\n(\xi))^2
\right.\\ \left.+K_i\w(\xi,\eta)
\overline{K_i\s(\xi,x)}(r_i\m(\xi))^2 \right]\,d\gamma_{i_{\eta}},
\end{multline}
Then Theorem 2.4 can be formulated as follows.
\begin{pred}\label{mX}
The minimax estimate of the value of functional $l(\varphi)$ has the form
\begin{equation} \label{2000a}
\widehat{\widehat{l(\varphi)}}
=\sum_{i=1}^N\int_{\gamma_i}\left(\overline{\hat u_i\n(x)}
y_i\n(x)+\overline{\hat
u_i\m(x)}y_i\m(x)\right)\, d\gamma_i+\hat c,
\end{equation}
where $\hat c=\int_{\Gamma}\bar zg_0\,d\Gamma,$
\begin{align}   \label{20'}
\hat u_i\n(x) &=(r_i\n(x))^2
\int_{\gamma_i}\left[K_i\p(x,\xi)p(\xi)+
K_i\q(x,\xi)\frac{\partial p(\xi)} {\partial
\nu}\right]\,d\gamma_{i_\xi}, \\ \hat u_i\m(x) &=(r_i\m(x))^2
\int_{\gamma_i}\left[K_i\s(x,\xi)p(\xi)+
K_i\w(x,\xi)\frac{\partial p(\xi)} {\partial
\nu}\right]\,d\gamma_{i_\xi},\quad i= \overline{1,N},\nonumber
\end{align}
and functions $z$ and $p$ are determined from the uniquely solvable problem
\begin{equation} \label{3r13}
z \in H^1_{\rm loc}(\Omega',\Delta),
\end{equation}
\begin{equation}\label{3r14}
-(\Delta +\bar
k^2)z(x)=\chi_{\omega_0}(x)l_0(x) \,\, \mbox{in}\,\, \Omega',
\end{equation}
\begin{equation}\label{3r15}
\frac{\partial z
}{\partial\nu}=0\,\,\mbox{on}\,\,\Gamma,
\end{equation}
\begin{equation}\label{3r16}
[z(x)]_{\gamma_i}=\int_{\gamma_i}\left[\tilde K_i\p(x,\eta)p(\eta)
+ \tilde K\q_i(x,\eta)\frac{\partial p(\eta)}{\partial\nu}\right]\,d\gamma_{i_{\eta}},\,\, \mbox
{on}\,\, \gamma_i,\,\,i=\overline{1,N},
\end{equation}
\begin{equation}\label{3r17}
\left[\frac{\partial z(x)
}{\partial\nu}\right]_{\gamma_i}=\int_{\gamma_i}\left[\tilde K_i\s(x,\eta)p(\eta)
+ \tilde K\w_i(x,\eta)\frac{\partial p(\eta)}{\partial\nu}\right]\,d\gamma_{i_{\eta}},\,\, \mbox
{on}\,\, \gamma_i,\,\,i=\overline{1,N},
\end{equation}
\begin{equation}\label{3r18}
\frac{\partial z}{\partial
r}+i\bar kz=o(1/r),\,\,r=|x|,\,\,r\to \infty,
\end{equation}
\begin{equation} \label{3r13'}
p \in H^1_{\rm loc}(\Omega',\Delta),
\end{equation}
\begin{equation}\label{3r19}
(\Delta+k^2)p(x)=0\,\, \mbox{in}\,\,\Omega',
\end{equation}
\begin{equation}\label{3r20}
\frac{\partial p}{\partial\nu}=Q^{-1}z_0\,\,\mbox{on}\,\,\Gamma,
\end{equation}
\begin{equation}\label{3r21}
[p(x)]_{\gamma_i}=0,\quad \left[\frac{\partial p(x)
}{\partial\nu}\right]_{\gamma_i}=0,\,\, \mbox
{on}\,\, \gamma_i,\,\,i=\overline{1,N},
\end{equation}
\begin{equation}\label{3r22}
\frac{\partial p}{\partial
r}-ikp=o(1/r),\,\,r=|x|, \,\,r\to \infty.
\end{equation}
\end{pred}

\subsection{Integral equation systems whose solutions are used to express minimax estimates}
In the previous section, we have obtained the integro-differential equations whose solutions are used to express minimax estimates.
In this section, we use the developed potential theory in Sobolev spaces and reduce these integro-differential equations to integral equations over an unclosed surface which is a union of the boundary of domain  $\Omega$
and surfaces $\gamma_i,$ $i=\overline{1,N},$ on which observations are made. This reduction allows one to decrease the dimensionality of the problem of finding minimax estimates.

We define first, in addition to the single- and double-layer potentials introduced in the previous sections, the corresponding boundary integral operators $S_k,$ $K_k,$ $K'_k,$
and $T_k$:
\begin{equation}
\label{eq8} (S_k\varphi)(x) := 2\int_\Gamma{\Phi_k(x,y)
\varphi(y)}\,d\Gamma_y, \,\,x \in \Gamma,
\end{equation}
\begin{equation}\label{eq8'}
(K_k\varphi)(x):= 2\int_\Gamma{\frac{\partial
\Phi_k(x,y)}{\partial \nu_y} \varphi(y)}\,d\Gamma_y, \,\,x \in
\Gamma,
\end{equation}
\begin{equation}
\label{eq9} (K'_k\psi)(x) := 2\int_\Gamma{\frac{\partial
\Phi_k(x,y)}{\partial \nu_x} \psi(y)}\,d\Gamma_y,\,\,x \in \Gamma,
\end{equation}
\begin{equation}\label{eq9'}
(T_k\psi)(x) := 2\frac{\partial}{\partial
\nu_x}\int_\Gamma{\frac{\partial \Phi_k(x,y)}{\partial \nu_y}
\psi(y)}\,d\Gamma_y,\,\,x \in \Gamma;
\end{equation}
in the three-dimensional case, their kernels are determined by the formulas
\begin{equation}
\label{uop}
\frac{\partial
\Phi_k(x,y)}{\partial \nu_y}=\Phi_k(x,y)\frac{(x-y,\nu_y)}{|x-y|^2}\left(1-ik|x-y|\right),
\end{equation}
\begin{equation}
\label{uop1}
\frac{\partial
\Phi_k(x,y)}{\partial \nu_x} =\Phi_k(x,y)\frac{(x-y,\nu_x)}{|x-y|^2}\left(ik|x-y|-1\right),
\end{equation}
\begin{multline}
\label{uop2}
\frac{\partial^2 \Phi_k(x,y)}{\partial
\nu_x\partial \nu_y}=\Phi_k(x,y)\left[\frac{(\nu_x,\nu_y)}{|x-y|^2}\left(1-ik{|x-y|}\right)
\right.\\\left.+\frac{(x-y,\nu_x)(x-y,\nu_y)}{|x-y|^4}\left(-3+3ik|x-y|+k^2|x-y|^3  \right)\right];
\end{multline}
in the two-dimensional case,
\begin{equation}
\label{uop'}
\frac{\partial
\Phi_k(x,y)}{\partial \nu_y}=\frac{ik(x-y,\nu_y)}{4|x-y|}H_1^{(1)}(k|x-y|),
\end{equation}
\begin{equation}
\label{uop1'}
\frac{\partial
\Phi_k(x,y)}{\partial \nu_x} =-\frac{ik(x-y,\nu_x)}{4|x-y|}H_1^{(1)}(k|x-y|),
\end{equation}
\begin{multline}
\label{uop2'}
\frac{\partial^2 \Phi_k(x,y)}{\partial
\nu_x\partial \nu_y}=\frac{ikH_1^{(1)}(k|x-y|)}{4|x-y|^3}\Bigl[(\nu_x,\nu_y)|x-y|^2
\Bigr.\\ \Bigl.-(x-y,\nu_x)(x-y,\nu_y)(|x-y|+1)\Bigr]+\frac{ik^2H_0^{(1)}(k|x-y|)}{4|x-y|^2}
(x-y,\nu_x)(x-y,\nu_y),
\end{multline}
where $H_1^1(z)$ denotes the order-one Hankel function of the first kind and $(\cdot,\cdot)$ the inner product in $\mathbb R^n,$ $n=2,3.$

Note that, for example, in the three-dimensional case, the kernels of integral operators
\eqref{eq8}$-$\eqref{eq9} have a weak singularity and the integral in the right-hand side of \eqref{eq9'} is understood as a Cauchy singular integral.

Let us formulate the properties of the operators introduced above
that are essential for the reduction of problem (\ref{3r13})$-$(\ref{3r22}) to a system of surface integral equations.

If $\Gamma$ is a $C^2$-surface, then the following operators are continuous at  $|s|\leq 1/2$:
$$
S_k:\, H^{-1/2+s}(\Gamma) \to H^{1/2+s}(\Gamma),
$$
$$
K_k:\, H^{1/2+s}(\Gamma) \to H^{3/2+s}(\Gamma),
$$
$$
K'_k:\, H^{-1/2+s}(\Gamma) \to H^{1/2+s}(\Gamma),
$$
\begin{equation}\label{rt}
T_k:\, H^{1/2+s}(\Gamma) \to H^{-1/2+s}(\Gamma)
\end{equation}
(similar statements are valid for the operators $S_{\!-\bar
k},$ $K_{\!-\bar k},$ $K'_{\!-\bar k},$ and $T_{\!-\bar k}$).

Also, operators  $K_k$ and $K_{\!-\bar k}'$ acting, respectively, from
$H^{1/2}(\Gamma)$ to $H^{1/2}(\Gamma)$ and from $H^{-1/2}(\Gamma)$ to
$H^{-1/2}(\Gamma)$ are compact according to (\ref{rt}) and the following equalities hold:
$$
\!(K_k\varphi,\psi)_{L^2(\Gamma)}=(\varphi,K'_{\!-\bar
k}\psi)_{L^2(\Gamma)}\,\,\forall \varphi\in H^{1/2}(\Gamma),
\psi\in H^{-1/2}(\Gamma),
$$
\begin{equation}
\label{rt1}
(S_k\varphi,\psi)_{L^2(\Gamma)}=(\varphi,S_{\!-\bar
k}\psi)_{L^2(\Gamma)}\,\,\,\,\,\forall \varphi\in
H^{-1/2}(\Gamma), \psi\in H^{-1/2}(\Gamma),
\end{equation}
$$
\!\!(T_k\varphi,\psi)_{L^2(\Gamma)}=(\varphi,T_{\!-\bar
k}\psi)_{L^2(\Gamma)}\,\,\,\,\forall \varphi\,\in H^{1/2}(\Gamma),
\psi\in H^{1/2}(\Gamma),
$$
where
$$
(f,g)_{L^2(\Gamma)}=\int_{\Gamma}f\bar g\,d\Gamma=:g(f)
$$
for every $f\in H^{1/2}(\Gamma)$ and $g\in H^{-1/2}(\Gamma).$

Denote by $\left(\mathcal V_{\Gamma}^k\psi\right)^{+}(x)$ and
$\left(\mathcal W_{\Gamma}^k\psi\right)^{+}(x)$ the restriction on
the domain $\mathbb R^3\setminus \bar\Omega$ of the single- and
double-layer potentials \eqref{simple} and \eqref{double} with a
density $\psi\in H^{s}(\Gamma),$ $s\in \mathbb R.$
The traces on $\Gamma=\partial\Omega$ of these functions and their derivatives
satisfy the relations \cite{Chen}, pp. 224, 225:
\begin{equation}\label{jump1}
\left(\mathcal V_{\Gamma}^k\psi\right)^{+}=\frac{1}{2}S_k\psi \quad\mbox{in}\quad H^{s+1}(\Gamma),
\end{equation}
\begin{equation}\label{jump2}
\left(\mathcal W_{\Gamma}^k\psi\right)^{+}=\frac{1}{2}(\psi+K_k\psi) \quad\mbox{in}\quad H^{s}(\Gamma),
\end{equation}
\begin{equation}\label{jump3}
\frac{\partial}{\partial\nu}\left(\mathcal V_{\Gamma}^k\psi\right)^{+}=-\frac{1}{2}(\psi-K'_k\psi) \quad\mbox{in}\quad H^{s}(\Gamma).
\end{equation}
\begin{equation}\label{jump4}
\frac{\partial}{\partial\nu}\left(\mathcal W_{\Gamma}^k\psi\right)^{+}=\frac{1}{2}T_k\psi \quad\mbox{in}\quad H^{s-1}(\Gamma).
\end{equation}
Similar relations are valid if we replace $k$ by $-\bar k.$

Denote by $D(\Omega)$ a countable set of positive values of wave number $k$
with a limiting point at infinity such that the homogeneous internal Nuemann problem
\begin{equation}\label{Sp}
\Delta v+k^2v=0\,\,\mbox{in}\,\,\Omega,
\end{equation}
\begin{equation}\label{Sp1}
v=0  \mbox{on} \Gamma
\end{equation}
has nontrivial solutions. Then \cite{Colton1}
\begin{equation}\label{pit}
N(I-K_k')=\left\{\left.\frac{\partial v}{\partial
\nu}\right|_{\Gamma}:\, \Delta v+k^2v=0
\,\,\mbox{in}\,\,\Omega,\,\,v=0 \,\ \mbox{on}\,\, \Gamma\right\},
\end{equation}
where $N(I-K_k')$ denotes the null-space of $I-K_k'.$

It is known that the solution $v_1(x)$ of the problem
\begin{equation}\label{1s}
\Delta v_1(x)+k^2v_1(x)=0\,\, \mbox{in}\,\,\mathbb R^3\setminus
\bar\Omega,
\end{equation}
\begin{equation}\label{2s}
\frac{\partial v_1}{\partial\nu}=f_1\,\,\mbox{on}\,\,\Gamma,
\end{equation}
\begin{equation}\label{3s}
\frac{\partial v_1}{\partial
r}-ikv_1=o(1/r),\,\,r=|x|,\,\,r\to \infty\,\,\mbox{if}\,\, \mbox{\rm
Im\,}k\geq 0,
\end{equation}
can be represented as
\begin{equation}\label{eq5'}
v_1(x)=\mathcal W^k_\Gamma\varphi_1(x)-\mathcal V^k_\Gamma f_1(x),\,\,
 x\in\mathbb R^3\setminus
\bar\Omega,
\end{equation}
where the function $\varphi_1:=v_1|_{\Gamma}$ which is the trace of this solution on $\Gamma$ can be determined directly or as a solution to the integral equation
\begin{equation}
\label{eq7}
\varphi_1(x) - K_k \varphi_1(x)=-S_k f_1(x),\,\, x\in\Gamma,
\end{equation}
when $k\notin D(\Omega);$ generally
(that is, for any $k\neq 0,\,\mbox{Im} \, k\geq 0$) it can be found from the equation
\begin{equation}\label{eq77}
\varphi_1(x)-K_k\varphi_1(x)-i\eta T_k\varphi_1(x)=-S_k f_1(x)-i\eta(f_1(x)
+K_k'f_1(x)),\,\, x\in\Gamma,
\end{equation}
in which a number $\eta\in \mathbb R,$ $\eta \neq 0,$ is chosen so that
\begin{equation}\label{ineq1q}
\eta\mbox{Re}\,k>0.
\end{equation}
The existence of solutions to these integral equations at any $f_1\in
H^{-1/2}(\Gamma)$ follows from the unique solvability of BVP (\ref{1s})$-$(\ref{3s}).
The  solution to (\ref{eq7}) is unique because
\begin{equation}\label{Fr}
\mbox{dim}\, N(I-K_k)=\mbox{dim}\, N(I-K_k')=\mbox{dim}\,
N(I-K_{-\bar k}')=0
\end{equation}
due to (\ref{pit}) and the Fredholm alternative; the uniqueness for (\ref{eq77}) is a consequence of the fact that the operator $I-K_k-i\eta T_k:\, H^{1/2}(\Gamma)\to
H^{-1/2}(\Gamma)$ defined by the left-hand side of
(\ref{eq77}) is an isomorphism \cite{Chen}.

A reasoning similar to that in \cite{Chen}$-$\cite{Colton2} enables us to prove that the solution $v_2(x)$ to the problem
\begin{equation}\label{1s'}
\Delta v_2(x)+\bar k^2v_2(x)=0\,\, \mbox{in}\,\,\mathbb R^3\setminus
\bar\Omega,
\end{equation}
\begin{equation}\label{2s'}
\frac{\partial v_2}{\partial\nu}=f_2\,\,\mbox{on}\,\,\Gamma,
\end{equation}
\begin{equation}\label{3s'}
\frac{\partial v_2}{\partial
r}+i\bar kv_2=o(1/r),\,\,r=|x|,\,\,r\to \infty,
\end{equation}
can be represented as
\begin{equation}\label{eq5a}
v_2(x)=\mathcal W^{-\bar k}_\Gamma\varphi_2(x)-\mathcal V^{-\bar k}_\Gamma f_2(x),\,\,
 x\in\mathbb R^3\setminus
\bar\Omega,
\end{equation}
where the function $\varphi_2:=v_2|_{\Gamma}$ which is a trace of this solution on $\Gamma$ can be determined directly or as a solution to the integral equation
\begin{equation}
\label{eq7a} \varphi_2(x) - K_{-\bar k} \varphi_2(x)=-S_{-\bar k} f_2(x),\,\, x\in\Gamma,
\end{equation}
for $k\notin D(\Omega)$;  generally, this function can be determined from the integral equation
(for any $k\neq 0,\,\mbox{Im} \, k\geq 0$)
\begin{multline}\label{eq77'}
\varphi_2(x)-K_{-\bar k}\varphi_2(x)+i\eta T_{-\bar k}\varphi_2(x)\\=-S_{-\bar k} f_2(x)+i\eta(f_2(x)
+K_{-\bar k}'f_2(x)),\,\, x\in\Gamma,
\end{multline}
where $\eta$ satisfies conditon \eqref{ineq1q}.

Indeed, the solution $v_2\in H^1(\mathbb R^3\setminus\bar\Omega)$ of the Helmholtz equation (\ref{1s'}) that satisfies radiation conditions (\ref{3s'})
and the property $\frac{\partial v_2}{\partial\nu}\in L^1(\Gamma)$
admits the integral representation
\begin{equation}\label{eq5a'}
v_2(x)=\int_\Gamma\frac{\partial
\Phi_{-\bar k}(x,y)}{\partial \nu_y}v_2(y)\,d\Gamma_y-\int_\Gamma\Phi_{-\bar k}(x,y)
\frac{\partial v_2(y)}{\partial\nu}\,d\Gamma_y,\quad
 x\in\mathbb R^3\setminus
\bar\Omega;
\end{equation}
replacing in this formula $\left.\frac{\partial v_2}{\partial\nu}\right|_{\Gamma}$ by $f_2,$ we obtain equality \eqref{eq5a}
for the solution $v_2(x)$ to problem \eqref{1s'}$-$\eqref{3s'}.

Calculating the traces on $\Gamma$ for both sides of \eqref{eq5a} and using relationships
\eqref{jump1} and \eqref{jump2} with $k$ replaced by $-\bar k$, obtain a boundary integral equation for $\varphi_2$
\begin{equation}\label{sled}
\varphi_2-K_{-\bar k}\varphi_2=-S_{-\bar k}f_2\quad
\mbox{on}\quad\Gamma;
\end{equation}
taking into account \eqref{jump3}$-$\eqref{jump4}, we can obtain
another integral equation for $\varphi_2$
\begin{equation}\label{sled1}
 T_{-\bar k}\varphi_2=f_2+K'_{-\bar k}f_2=0\quad \mbox{on}\quad\Gamma.
\end{equation}
Multiplying both sides of \eqref{sled1}
by a number $\eta$ and adding to \eqref{sled1}, we find that $\varphi_2$ satisfies integral equation \eqref{eq77'}.

The uniqueness of solution of this integral equation follows from
the fact that operator $K_{-\bar k}$ is adjoint to  $K'_k$ and
operator $T_{-\bar k}$ is adjoint to  $T_k$; therefore,
$I-K_{-\bar k}+i\eta T_{-\bar k}$ is adjoint to  $I-K'_k-i\eta
T_k.$ Also,
$$
\mbox{Im\,}(I-K_{-\bar k}+i\eta T_{-\bar k})=\mbox{Ker\,}(I-K'_k-i\eta T_k)^\bot.
$$
In \cite{Colton1} it is proved (see Theorem 3.34) that under condition \eqref{ineq1q}
$\mbox{Ker\,}(I-K'_k-i\eta T_k)=\emptyset.$ Thus, if \eqref{ineq1q} holds, then integral equation \eqref{eq77'} is uniquely solvable.


In order to reduce problem (\ref{3r13})$-$(\ref{3r22}) in
unbounded domain $\mathbb R^3\setminus \bar\Omega$ to a system of
boundary integral equations let us apply the results formulated
above in Subsection 2.5.

Introduce functions
 $\varrho_i\n(\hat u_i\n,\hat u_i\m)\in
H_{0\,0}^{1/2}(\gamma_i)$ and $\varrho_i\m(\hat u_i\n,\hat u_i\m)\in
L^{2}(\gamma_i)$ defined on $\gamma_i,$ $i=\overline{1,N}:$
\begin{multline}\label{r16'a}
\varrho_i\n(\hat u_i\n,\hat u_i\m)=\varrho_i\n(\cdot;\hat u_i\n,\hat u_i\m)\\:=\int_{\gamma_i}\left[\overline{K_i\q(\xi,\cdot)}\hat
u_i\n(\xi)+ \overline{K_i\w(\xi,\cdot)}\hat u_i\m(\xi)\right]\,d\gamma_{i_{\xi}},
\end{multline}
\begin{multline}\label{r16''a}
\varrho_i\m(\hat u_i\n,\hat u_i\m)=\varrho_i\n(\cdot;\hat u_i\n,\hat u_i\m)\\:
= -\int_{\gamma_i}\left[ \overline{K_i\p(\xi,\cdot)}\hat u_i\n(\xi)+
\overline{K_i\s(\xi,\cdot)}\hat u_i\m(\xi)\right]\,d\gamma_{i_{\xi}},
\end{multline}

Introduce also a function  $z_{in}=z_{in}(\cdot;\hat u)\in H^1_{\rm loc}(\mathbb R^n\setminus(\cup_{i=1}^N\bar\gamma_i), \Delta)$
which solves the BVP
\begin{equation}\label{li7}
-(\Delta+\bar k^2)z_{in}=\chi_{\omega_0}(x)l_0\,\, \mbox{in}\,\,\mathbb R^n\setminus(\cup_{i=1}^N\bar\gamma_i),
\end{equation}
\begin{equation}\label{li8}
[z_{in}]_{\gamma_i}=\varrho\n_i(\hat u_i\n,\hat u_i\m),\quad\left[\frac{\partial z_{in}
}{\partial\nu}\right]_{\gamma_i}=\rho\m_i(\hat u_i\n,\hat u_i\m)\,\,\mbox{on}\,\,\gamma_i,\,\,i=\overline{1,N},
\end{equation}
\begin{equation}\label{li9}
\frac{\partial z_{in}}{\partial
r}+i\bar kz_{in}=o(1/r),\,\,r=|x|,\,\,r\to \infty;
\end{equation}
in line with Theorem \ref{ska}, this function, in the domain $\mathbb R^n\setminus(\cup_{i=1}^N\bar\gamma_i)$, is determined according to
\begin{equation*}
z_{in}=\Sigma_{i=1}^N\left(\mathcal V^{-\bar k}_{\gamma_i}\varrho\m_i(\hat u_i\n,\hat u_i\m)-\mathcal W^{-\bar k}_{\gamma_i}\varrho\n_i(\hat u_i\n,\hat u_i\m)\right)+N_{-\bar k}l_0
\end{equation*}
\begin{multline}\label{lkj}
=2\Sigma_{i=1}^N\left(\int_{\gamma_i}\Phi_{-\bar k}(\cdot,y)\varrho\m_i(y;\hat u_i\n,\hat u_i\m)\,d\gamma_{i_y}-\int_{\gamma_i}\frac{\partial\Phi_{-\bar k}(\cdot,y)}{\partial \nu_{i_y}
}\varrho\n_i(y;\hat u_i\n,\hat u_i\m)\,d\gamma_{i_y}\right)\\+\int_{\omega_0}\Phi_{-\bar k}(\cdot,y)l_0(y)\,dy.
\end{multline}

Then the solution to problem \eqref{3r13}$-$\eqref{3r18} can be represented in
$\Omega'$ as
\begin{equation}
\label{li3}
z=z_s+z_{in},
\end{equation}
where the function $z_s=z_s(\cdot;\hat u),$
\begin{equation*}
\hat u=(\hat u_1\n, \ldots ,
\hat u_{N}\n, \hat u_1\m, \ldots ,\hat u_{N}\m)\in \left(L^2(\gamma_1)
\times\ldots \times L^2(\gamma_N)\right)^2,
\end{equation*}
is determined from the solution to the following problem
\begin{equation} \label{3r13aa}
z_s\in H^1_{\rm loc}(\Omega',\Delta),
\end{equation}
\begin{equation}\label{li4}
(\Delta+\bar k^2)z_s=0\,\, \mbox{in}\,\,\Omega',
\end{equation}
\begin{equation}\label{li5}
\frac{\partial z_s}{\partial\nu} =
-\frac{\partial z_{in}}{\partial\nu}\,\,\mbox{on}\,\,\Gamma,
\end{equation}
\begin{equation}\label{li6}
\frac{\partial z_s}{\partial
r}+i\bar kz_s=o(1/r),\,\,r=|x|,\,\,r\to \infty,
\end{equation}
where
\begin{equation*}
\frac{\partial z_{in}}{\partial \nu}=\Sigma_{i=1}^N\left(\frac{\partial \mathcal V^{-\bar k}_{\gamma_i}\varrho_i\m(\hat u_i\n,\hat u_i\m)}{\partial \nu}
-\frac{\partial \mathcal W^{-\bar k}_{\gamma_i}\varrho_i\n(\hat u_i\n,\hat u_i\m)}{\partial \nu}\right)+\frac{\partial N_{-\bar k}l_0}{\partial \nu}
\end{equation*}
\begin{multline}\label{li10}
\!=\!2\Sigma_{i=1}^N\left(\!\int_{\gamma_i}\frac{\partial \Phi_{-\bar k}(\cdot,y)}{\partial \nu}\varrho_i\m(y;\hat u_i\n,\hat u_i\m)d\gamma_{i_y}\!-\!\int_{\gamma_i}\frac{\partial^2\Phi_{-\bar k}(\cdot,y)}{\partial \nu\partial \nu_y
}\varrho_i\n(y;\hat u_i\n,\hat u_i\m)d\gamma_{i_y}\!\right)\\
+\int_{\omega_0}\frac{\partial \Phi_{-\bar k}(\cdot,y)}{\partial \nu}l_0(y)\,dy.
\end{multline}
According to (\ref{eq5a}) and (\ref{eq77'}), where
$v_2=z_s$, $f_2=\left.-\frac{\partial z_{in}}{\partial \nu}\right|_\Gamma,$ and $\varphi_2=z_s|_{\Gamma},$
the trace of $z_s$ on $\Gamma$ denoted by
\begin{equation}\label{eq7c}
 \psi:=z_s|_{\Gamma}=\left.z\right|_{\Gamma}-\left.z_{in}\right|_{\Gamma}
\end{equation}
(see \eqref{li3}) satisfies the integral equation
\begin{equation}\label{eq77a}
\psi-K_{-\bar k}\psi+i\eta T_{-\bar k}\psi=S_{-\bar k}\left. \frac{\partial z_{in}}{\partial \nu}\right|_\Gamma-i\eta(\left.\frac{\partial z_{in}}{\partial \nu}\right|_\Gamma
+K_{-\bar k}'\left.\frac{\partial z_{in}}{\partial \nu}\right|_\Gamma).
\end{equation}
Introduce also the notations
\begin{equation}\label{eq7b}
\chi:=p|_{\Gamma}, \quad\varphi\n_i:=p|_{\gamma_i},\quad \varphi\m_i:=\left.\frac{\partial p} {\partial
\nu}\right|_{\gamma_i},\quad i=\overline{1,N}.
\end{equation}
Then, by virtue of \eqref{re27} and \eqref{li3} $\frac{\partial p}{\partial\nu}=q_1^{-2}z=q_1^{-2}(\psi+z_{in})$ on $\Gamma$.

Taking into account relationships \eqref{eq5'}, \eqref{eq77} in which we set $f_1=q_1^{-2}z=q_1^{-2}(\psi+z_{in}),$ $\varphi_1=\chi,$ and $v_1=p,$ we obtain an integral representation for the solution $p$ of BVP
\eqref{re26}$-$\eqref{re28''} in the domain $\mathbb R^3\setminus
\bar\Omega $
\begin{equation}\label{a333}
p=\mathcal W^k_\Gamma\chi-\mathcal
V^k_\Gamma q_1^{-2}(\psi+z_{in}\left.\right|_\Gamma),
i=\overline{1,N},
\end{equation}
where $\chi=p|_{\Gamma}$ satisfies on $\Gamma$ a boundary integral equation
\begin{equation}\label{a22}
(I-K_k-i\eta T_k)\chi=-(S_k
+i\eta(I
+K_k'))q_1^{-2}(\psi+z_{in}\left.\right|_\Gamma),
\end{equation}
From \eqref{a333} it follows that
\begin{equation}\label{a222}
\left.\frac{\partial p}{\partial\nu}\right|_{\gamma_i}=
\left.\frac{\partial\mathcal W^k_\Gamma\chi}{\partial\nu}\right|_{\gamma_i}-\left.\frac{\partial\mathcal
V^k_\Gamma q_1^{-2}(\psi+z_{in}\left.\right|_\Gamma)}{\partial\nu}\right|_{\gamma_i},\,\,
i=\overline{1,N}.
\end{equation}
and
\begin{equation}\label{a33a}
\left.p\right|_{\gamma_i}=\mathcal W^k_\Gamma\chi\left.\right|_{\gamma_i}-\mathcal
V^k_\Gamma q_1^{-2}(\psi+z_{in}\left.\right|_\Gamma)\left.\right|_{\gamma_i},
i=\overline{1,N},
\end{equation}
Since $\varphi\n_i=p|_{\gamma_i}$ and $\varphi\m_i=\left.\frac{\partial p} {\partial
\nu}\right|_{\gamma_i},$ $i=\overline{1,N},$ the latter equalities combined with (\ref{eq77a}) and (\ref{a22}) mean that functions $\psi,$ $\chi,$ $\varphi\n_i,$ $\varphi\m_i$
and $\hat u_i\n,$ $\hat u_i\m$ defined by (\ref{eq7c}), (\ref{eq7b}), and (\ref{20}) solve the following integral equation system
\begin{equation}\label{a1i}
(I-K_{\!-\bar k}+i\eta T_{\!-\bar k})\psi=(S_{-\bar k}
-i\eta(I +K_{-\bar k}')) \left.\frac{\partial z_{in}}{\partial
\nu}\right|_\Gamma,
\end{equation}
\begin{equation}\label{a2i}
(I-K_k-i\eta T_k)\chi=-(S_k
+i\eta(I
+K_k'))q_1^{-2}(\psi+z_{in}\left.\right|_\Gamma),
\end{equation}
\begin{equation}\label{a3i}
\varphi\n_i=\mathcal W^k_\Gamma\chi\left.\right|_{\gamma_i}-\mathcal
V^k_\Gamma q_1^{-2}(\psi+z_{in}\left.\right|_\Gamma)\left.\right|_{\gamma_i},
i=\overline{1,N},
\end{equation}
\begin{equation}\label{a4i}
\varphi\m_i=\left.\frac{\partial\mathcal W^k_\Gamma\chi}{\partial\nu}\right|_{\gamma_i}-\left.\frac{\partial\mathcal
V^k_\Gamma q_1^{-2}(\psi+z_{in}\left.\right|_\Gamma)}{\partial\nu}\right|_{\gamma_i}
,\,\,
i=\overline{1,N},
\end{equation}
\begin{align}   \label{20a}
\hat u_i\n &=(r_i\n)^2
\int_{\gamma_i}\left[K_i\p(\cdot,\xi)\varphi\n_i(\xi)+
K_i\q(\cdot,\xi)\varphi\m_i(\xi)\right]\,d\gamma_{i_\xi}, \\ \hat u_i\m &=(r_i\m)^2
\int_{\gamma_i}\left[K_i\s(\cdot,\xi)\varphi\n_i(\xi)+
K_i\w(\cdot,\xi)\varphi\m_i(\xi)\right]\,d\gamma_{i_\xi},\quad i= \overline{1,N}. \label{21a}
\end{align}
Resolve this system with respect to functions $\hat u_i\n$ and $\hat u_i\m.$
Replace in \eqref{r16'a} and \eqref{r16''a}  $\hat u_i\n$ and $\hat u_i\m$ by their expressions \eqref{20a} and \eqref{21a} to obtain
$\varrho_i\n(\hat u_i\n,\hat u_i\m)=\rho\n_i(\varphi\n_i,\varphi\m_i)$ and $\varrho_i\m(\hat u_i\n,\hat u_i\m)=\rho\m_i(\varphi\n_i,\varphi\m_i)$
where
\begin{equation}\label{rho}
\rho\n_i(\varphi\n_i,\varphi\m_i):=\int_{\gamma_i}\left[\tilde K_i\p(\cdot,\eta)\varphi\n_i(\eta)
+ \tilde K\q_i(\cdot,\eta)\varphi\m_i(\eta)\right]\,d\gamma_{i_{\eta}},
\end{equation}
\begin{equation}\label{rho1}
\rho\m_i(\varphi\n_i,\varphi\m_i):=\int_{\gamma_i}\left[\tilde K_i\s(\cdot,\eta)\varphi\n_i(\eta)
+ \tilde K\w_i(\cdot,\eta)\varphi\m_i(\eta)\right]\,d\gamma_{i_{\eta}},\,\, \end{equation}
and $\tilde K_i^{(i,j)}(\cdot,\cdot),$ $i,j=1,2,$ are determined according to
\eqref{3r1}$-$\eqref{3r4}.
Next, replacing in (\ref{lkj}) functions $\varrho_i\n(\hat u_i\n,\hat u_i\m)$ by $\rho\n_i(\varphi\n_i,\varphi\m_i)$ and $\varrho_i\m(\hat u_i\n,\hat u_i\m)$ by $\rho\m_i(\varphi\n_i,\varphi\m_i)$ we have
\begin{equation}\label{gfu}
z_{in}=\Sigma_{i=1}^N\left(\mathcal V^{-\bar k}_{\gamma_i}\rho\m_i(\varphi\n_i,\varphi\m_i)-\mathcal W^{-\bar k}_{\gamma_i}\rho\n_i(\varphi\n_i,\varphi\m_i)\right)+N_{-\bar k}l_0.
\end{equation}
Substituting this expression into (\ref{a1i})$-$(\ref{a4i}) we conclude that $\psi,$ $\chi,$ $\varphi\n_i,$ and $\varphi\m_i$ satisfy the integral equation system
\begin{multline}\label{a1a}
(I-K_{\!-\bar k}+i\eta T_{\!-\bar k})\psi=(S_{-\bar k}
-i\eta(I +K_{-\bar k}'))\\
\left(\Sigma_{i=1}^N\left(\left.\frac{\partial \mathcal V^{-\bar k}_{\gamma_i}\rho\m_i(\varphi\n_i,\varphi\m_i)} {\partial
\nu}\right|_{\Gamma}-\left.\frac{\partial\mathcal W^{-\bar k}_{\gamma_i}
\rho\n_i(\varphi\n_i,\varphi\m_i)}{\partial
\nu}\right|_{\Gamma}\right)+\left.\frac{\partial N_{-\bar k}l_0}{\partial
\nu}\right|_{\Gamma}\right),
\end{multline}
\begin{multline}\label{a2a}
(I-K_k-i\eta T_k)\chi=-(S_k +i\eta(I
+K_k'))\\
q_1^{-2}\left(\psi+\Sigma_{i=1}^N\left(\left.\mathcal V^{-\bar k}_{\gamma_i}\rho\m_i(\varphi\n_i,\varphi\m_i)\right|_{\Gamma}-\left.\mathcal W^{-\bar k}_{\gamma_i}\rho\n_i(\varphi\n_i,\varphi\m_i)\right|_{\Gamma}\right)+\left.N_{-\bar k}l_0\right|_{\Gamma}\right),
\end{multline}\\[-30pt]
\begin{multline}\label{a3a}
\varphi\n_i=\left.\mathcal W_{\Gamma}^k\chi\right|_{\gamma_i}
-\mathcal
V_{\Gamma}^kq_1^{-2}\Bigl(\psi+\sum_{i=1}^N\Bigl(\left.\mathcal V^{-\bar k}_{\gamma_i}\rho\m_i(\varphi\n_i,\varphi\m_i)\right|_{\Gamma}\Bigr.\Bigr.\\ \Bigl.\Bigl.-\left.\mathcal W^{-\bar k}_{\gamma_i}\rho\n_i(\varphi\n_i,\varphi\m_i)\right|_{\Gamma}\Bigr)
+\left.N_{-\bar k}l_0\right|_{\Gamma}\left.\Bigr)\right|_{\gamma_i},
\,\,i=\overline{1,N},
\end{multline}\\[-40pt]
\begin{equation*}
\varphi\m_i=\left.\frac{\partial\mathcal W_{\Gamma}^k\chi}{\partial\nu}\right|_{\gamma_i}
\end{equation*}
\begin{multline}\label{a4a}
-\left.\frac{\partial\mathcal
V_{\Gamma}^kq_1^{-2}\left(\psi+\Sigma_{i=1}^N\left(\left.\mathcal V^{-\bar k}_{\gamma_i}\rho\m_i(\varphi\n_i,\varphi\m_i)\right|_{\Gamma}-\left.\mathcal W^{-\bar k}_{\gamma_i}\rho\n_i(\varphi\n_i,\varphi\m_i)\right|_{\Gamma}\right)+N_{-\bar k}l_0\left.\right|_{\Gamma}\right)}{\partial\nu}\right|_{\gamma_i},
\\ \,\,i=\overline{1,N}.
\end{multline}

Summing up the above reasoning and taking into consideration Theorems \ref{t4} and \ref{mX} we arrive at the following result.
\begin{pred}\label{t3.2'}
The minimax estimate of
$l(\varphi)$
has the form
\begin{equation} \label{2001iu}
\widehat{\widehat{l(\varphi)}}
=\sum_{i=1}^N\int_{\gamma_i}\left(\overline{\hat u_i\n(x)}
y_i\n(x)+\overline{\hat
u_i\m(x)}y_i\m(x)\right)\, d\gamma_i+\hat c,
\end{equation}
where
\begin{equation} \label{2002iu}
\hat c=\int_{\Gamma}\overline{\left[ \psi
+\sum_{i=1}^N\left(\left.\mathcal V^{-\bar k}_i\rho\m_i(\varphi\n_i,\varphi\m_i)\right|_{\Gamma}
-\left.\mathcal W^{-\bar k}_i\rho\n_i(\varphi\n_i,\varphi\m_i)
\right|_{\Gamma}\right)+N_{-\bar k}
l_0\left.\right|_{\Gamma}\right]}g_0\,d\Gamma,
\end{equation}
\begin{align}
\hat u_i\n(x) &=(r_i\n(x))^2
\int_{\gamma_i}\left[K_i\p(x,\xi)\varphi\n_i(\xi)+
K_i\q(x,\xi)\varphi\m_i(\xi)\right]\,d\gamma_{i_\xi},
\label{2003iu}\\ \hat u_i\m(x) &=(r_i\m(x))^2
\int_{\gamma_i}\left[K_i\s(x,\xi)\varphi\n_i(\xi)+
K_i\w(x,\xi)\varphi\n_i(\xi) \right]\,d\gamma_{i_\xi},
\quad i= \overline{1,N}.\nonumber
\end{align}
The auxiliary function
$$
\psi:=\left.z\right|_{\Gamma}-\Sigma_{i=1}^N\left(\left.\mathcal V^{-\bar k}_{\gamma_i}\rho\m_i(\varphi\n_i,\varphi\m_i)\right|_{\Gamma}-\left.\mathcal W^{-\bar k}_{\gamma_i}\rho\n_i(\varphi\n_i,\varphi\m_i)\right|_{\Gamma}\right)-\left.N_{-\bar k}l_0\right|_{\Gamma}
$$
defined on $\Gamma$ and functions $\chi:=p|_{\Gamma},$  $\varphi\n_i:=p|_{\gamma_i},$ and $\varphi\m_i:=\frac{\partial p} {\partial \nu_{\gamma_i}},$ $i=\overline{1,N},$
are determined from the solution to the integral equation system \eqref{a1a}$-$\eqref{a4a},
in which
$\eta$ is an arbitrary real number such that $\eta\mbox{\rm Re}\,k>0.$
This system is uniquely solvable for all values of wave numbers $k$,  $\mbox{\rm Im\,} k\geq 0.$

The estimation error
$\sigma=l(p)^{1/2},$
where
\begin{multline}\label{a3as'}
p=\mathcal W_{\Gamma}^k\chi-\mathcal
V_{\Gamma}^kq_1^{-2}\Bigl(\psi+\Sigma_{i=1}^N\left(\left.\mathcal V^{-\bar k}_{\gamma_i}\rho\m_i(\varphi\n_i,\varphi\m_i)\right|_{\Gamma}\Bigr.\right.\\
\left.\Bigl.-\left.\mathcal W^{-\bar k}_{\gamma_i}\rho\n_i(\varphi\n_i,\varphi\m_i)\right|_{\Gamma}\right)+N_{-\bar k}l_0\left.\right|_{\Gamma}\Bigr)
\end{multline}
in the domain $\mathbb R^3\setminus\bar\Omega$.
\end{pred}
\begin{proof}
It is necessary to prove only the unique solvability of system \eqref{a1a}$-$\eqref{a4a} which follows from the unique solvability of the system of integro-differential equations \eqref{20}$-$\eqref{re28''}.

Indeed, let the integral equation system \eqref{a1a}$-$\eqref{a4a}
has another solution
$\tilde\psi,$  $\tilde\chi,$  $\tilde\varphi\n_i,$
$\tilde\varphi\m_i,$ $i=\overline{1,N}.$
Introduce functions  $\tilde{\hat u}_i\n(x),$
$\tilde{\hat u}_i\m(x),$
$\tilde p,$ and $\tilde z_{in}$ by formulas
\eqref{20a}, \eqref{21a}, \eqref{a3as'}, and \eqref{gfu}
in which $\psi,$  $\chi,$
$\varphi\n_i,$ and $\varphi\m_i,$ $i=\overline{1,N}$
are replaced by $\tilde\psi,$
$\tilde\chi,$  $\tilde\varphi\n_i,$ and $\tilde\varphi\m_i,$
$i=\overline{1,N},$
and functions $\tilde z_s$ and $\tilde z$
by $\tilde z_s=\mathcal W_{\Gamma}^{-\bar k}\tilde\psi-\mathcal V_{\Gamma}^{-\bar k}\tilde z_{in}|_{\Gamma}$ and
$\tilde z=\tilde z_s+\tilde z_{in}.$
Then from \eqref{a3a}$-$\eqref{a4a} it follows that
$\tilde\varphi\n_i=\tilde p|_{\gamma_i}$ and
$\tilde\varphi\m_i
=\left.\frac{ \partial \tilde p}{\partial \nu}\right|_{\gamma_i},$
$i=\overline{1,N}.$
Theorem \ref{ska} and the equalities
$\varrho_i\n(\tilde{\hat u}_i\n,\tilde{\hat u}_i\m)=\rho\n_i(\tilde\varphi\n_i,\tilde\varphi\m_i),$
$\varrho_i\m(\tilde{\hat u}_i\n,\tilde{\hat u}_i\m)=\rho\m_i(\tilde\varphi\n_i,\tilde\varphi\m_i)$
imply that
$\tilde z$ and $\tilde p$ will also satisfy
integro-differential equation system \eqref{20}--\eqref{re28''}
which is uniquely solvable. The latter statement and the fact that operators $I-K_k-i\eta T_k,\,I-K_{-\bar k}
+i\eta T_{-\bar k} :\, H^{1/2}(\Gamma)\to
H^{-1/2}(\Gamma)$ are isomorphic mappings yield
$\tilde\varphi\n_i=\varphi\n_i$,
$\tilde\varphi\m_i=\varphi\m_i,$ $i=\overline{1,N},$ and
$\tilde\psi=\psi,$  $\tilde\chi=\chi.$
\end{proof}

Using Theorem \ref{t6},
the notations
$$
D_{-\bar k}:=\sum_{i=1}^N\left(\mathcal V^{-\bar k}_id_i^{(2)}-\mathcal W^{-\bar k}_id_i^{(1)}\right),
$$
$$
d_i^{(1)}:=-\int_{\gamma_i}\overline{K_i\q(\xi,\cdot)}(r_i\n(\xi))^2y\n_i(\xi,\omega)\,d\gamma_i
-\int_{\gamma_i}\overline{K_i\w(\xi,\cdot)}(r_i\m(\xi))^2y\m_i(\xi,\omega)\,d\gamma_i,
$$
$$
d_i^{(2)}:=\int_{\gamma_i}\overline{K_i\p(\xi,\cdot)}(r_i\n(\xi))^2y\n_i(\xi,\omega)\,d\gamma_i
+\int_{\gamma_i}\overline{K_s\w(\xi,\cdot)}(r_i\m(\xi))^2y\m_i(\xi,\omega)\,d\gamma_i,
$$
$\tilde\psi:=\hat p_s|_{\Gamma}-D_{-\bar k}\left.\right|_\Gamma$,
 $\tilde\chi:=\hat\varphi|_{\Gamma},$
$\tilde\varphi\n_i:=\hat\varphi|_{\gamma_i},$
$\tilde\varphi\m_i:=\left.\frac{\partial \hat\varphi(\xi)}
{\partial \nu}\right|_{\gamma_i},$ $i=\overline{1,N},$ and the
reasoning that led to the proof of Theorem \ref{t3.2'}, we can
prove the following
\begin{pred}\label{t3.2a'}
The minimax estimate of
$l(\varphi)$
has the form
\begin{equation} \label{2001}
\widehat{\widehat{l(\varphi)}}
=l(\hat\varphi),
\end{equation}
where
\begin{multline}\label{a3ass}
\hat\varphi=\mathcal W_{\Gamma}^k\tilde\chi-\mathcal
V_{\Gamma}^k\Bigl[q_1^{-2}\Bigl(\tilde\psi+\sum_{i=1}^N\left(\left.\mathcal V^{-\bar k}_{\gamma_i}\rho\m_i(\tilde\varphi\n_i,\tilde\varphi\m_i)\right|_{\Gamma}\Bigr.\Bigr.\right.\\
\left.\Bigl.\Bigl.-\left.\mathcal W^{-\bar k}_{\gamma_i}\rho\n_i(\tilde\varphi\n_i,\tilde\varphi\m_i)\right|_{\Gamma}\right)+D_{-\bar k}\left.\right|_{\Gamma}\Bigr)+h_0\Bigr],
\end{multline}
in the domain $\mathbb R^3\setminus\bar\Omega$ and
functions $\tilde\psi,$ $\tilde\chi,$  $\tilde\varphi\n_i,$ and $\tilde\varphi\m_i,$ $i=\overline{1,N},$ are determined from the solution to the uniquely solvable integral equation system
\begin{multline}\label{a1ag}
(I-K_{\!-\bar k}+i\eta T_{\!-\bar k})\tilde\psi=(S_{-\bar k}
-i\eta(I +K_{-\bar k}')) \\
\left(\Sigma_{i=1}^N\left(\left.\frac{\partial \mathcal V^{-\bar k}_{\gamma_i}\rho\m_i(\tilde\varphi\n_i,\tilde\varphi\m_i)} {\partial
\nu}\right|_{\Gamma}-\left.\frac{\partial\mathcal W^{-\bar k}_{\gamma_i}
\rho\n_i(\tilde\varphi\n_i,\tilde\varphi\m_i)}{\partial
\nu}\right|_{\Gamma}\right)+\left.\frac{\partial D_{-\bar k}}{\partial
\nu}\right|_{\Gamma}\right),
\end{multline}
\begin{equation}\label{a2ag}
(I-K_k-i\eta T_k)\tilde\chi=-(S_k +i\eta(I
+K_k'))
\end{equation}
\begin{equation*}
\left[q_1^{-2}\left(\tilde\psi+\sum_{i=1}^N\left(\mathcal V^{-\bar k}_i\rho\m_i(\tilde\varphi\n_i,\tilde\varphi\m_i)\left.\right|_\Gamma-\mathcal W^{-\bar k}_i\rho\n_i(\tilde\varphi\n_i,\tilde\varphi\m_i)\left.\right|_\Gamma\right)+D_{-\bar k}
\left.\right|_{\Gamma}\right)
+h_0\right],
\end{equation*}
\begin{multline}\label{a3ag}
\tilde\varphi\n_i=\mathcal W_k\tilde\chi\left.\right|_{\gamma_i}-\mathcal
V_k\Bigl[q_1^{-2}\Bigl(\tilde\psi+\sum_{i=1}^N\left(\mathcal V^{-\bar k}_i\rho\m_i(\tilde\varphi\n_i,\tilde\varphi\m_i)\left.\right|_\Gamma-\mathcal W^{-\bar k}_i\rho\n_i(\tilde\varphi\n_i,\tilde\varphi\m_i)\left.\right|_\Gamma\right)\Bigr.\Bigr.
\\ \Bigl.\Bigl.\left.+D_{-\bar k}\left.\right|_\Gamma\Bigr)+h_0\Bigr]\right|_{\gamma_i},
\,\,
i=\overline{1,N},
\end{multline}
\begin{equation*}\label{a4ag}
\left.\tilde\varphi\m_i=\frac{\partial\mathcal W_k\tilde\chi}{\partial\nu}\right|_{\gamma_i}
\end{equation*}
\begin{multline}\label{a5ag}
\!\!\!\!-\!\left.\frac{\partial\mathcal V_k\Bigl[q_1^{-2}\Bigl(\tilde\psi\!+\!\sum_{i=1}^N\left(\mathcal V^{-\bar k}_i\rho\m_i(\tilde\varphi\n_i,\tilde\varphi\m_i)\left.\right|_\Gamma\!-\!\mathcal W^{-\bar k}_i\rho\n_i(\tilde\varphi\n_i,\tilde\varphi\m_i)\left.\right|_\Gamma\right)
\!+\!D_{-\bar k}\left.\right|_\Gamma\Bigr)\!+\!h_0\Bigr]}{\partial\nu}\right|_{\gamma_i},
\\i=\overline{1,N},
\end{multline}
where $\eta$ is an arbitrary real number such that $\eta\mbox{\rm Re}\,k>0$
and functions
$
\rho\n_i(\tilde\varphi\n_i,\tilde\varphi\m_i)
$
and
$
\rho\m_i(\tilde\varphi\n_i,\tilde\varphi\m_i)
$
are determined by \eqref{rho} and \eqref{rho1} in which $\varphi\n_i$ and $\varphi\m_i$
should be replaced, respectively, by $\tilde\varphi\n_i$ and $\tilde\varphi\m_i.$
\end{pred}

Note that integral equation systems \eqref{a1a}$-$\eqref{a4a} and
\eqref{a1ag}$-$\eqref{a5ag} are singular.

If $k\notin D(\Omega),$ then repeating the reasoning used in the
proof of Theorems \ref{t3.2'} and \ref{t3.2a'} where equalities of
the form \eqref{eq77} and \eqref{eq77'} are replaced,
respectively, by those of the form \eqref{eq7} and \eqref{eq7a} we
see
that the minimax estimate of $l(\varphi)$ may be found from \eqref{2001iu}$-$\eqref{2003iu} or \eqref{2001}, \eqref{a3ass},
where the functions $\psi,$ $\chi,$  $\varphi\n_i,$ $\varphi\m_i,$ $i=\overline{1,N},$
and $\tilde\psi,$ $\tilde\chi,$  $\tilde\varphi\n_i,$ $\tilde\varphi\m_i,$ $i=\overline{1,N},$
are determined from the solutions of the weakly singular integral equation systems \eqref{a1a}$-$\eqref{a4a} and
 \eqref{a1ag}$-$\eqref{a5ag} with $\eta=0.$

\begin{predlllll}
The assumption that surfaces $\gamma_i,$ $i=\overline{1,N},$ are
pairwise non-overlapping is not essential. Slightly changing the
proof, one can extend all results of this chapter to the case when
surfaces $\gamma_i$ intersect on a finite system of contours.

\end{predlllll}

\begin{predlllll} The method proposed in this chapter enables one to solve the problem of
the minimax estimation of the value of a functional defined on
$\Phi(x,t):=\mbox{Re} \left[e^{-i\omega t}\varphi(x)\right]$ of
the form
$$
L(\Phi):=\int_{t_0}^{T}\!\!\!\int_{\omega_0}l_0(x,t)\Phi(x,t)\,dx\,dt
$$
from the observations
\begin{equation*}
y_{i}\n(x,t)=\int_{\gamma_i}
K_{i}\p(x,\xi)\Phi(\xi,t)\,d\gamma_{i_\xi}+\int_{\gamma_i}
K_{i}\q(x,\xi)\frac{\partial\Phi(\xi,t)}{\partial\nu}\,d\gamma_{i_\xi}
+ \eta_{i}\n(x,t),
\end{equation*}
\begin{equation*}
y_{i}\m(x,t)=\int_{\gamma_i}
K_{i}\s(x,y)\Phi(\xi,t)\,d\gamma_{i_\xi}+\int_{\gamma_i}
K_{i}\w(x,y)\frac{\partial\Phi(\xi,t)}{\partial\nu}\,d\gamma_{i_\xi}
+\eta_{i}\m(x,t),
\end{equation*}
$i=\overline{1,N},$ in a time interval from $t=t_0$ to $t=T.$
Here we assume that for $f \in G_0,$ $\omega>0,$
$l_0\in L^2(\omega_0\times (t_0,T))$ is a given function, and
$\eta\n_i(x,t)$ and $\eta\m_i(x,t)$ are observations errors which are realizations of random fields defined on
$\gamma_i\times(t_0,T)$ that are continuous in the mean-square sense and have zero expectation and unknown second moments $\mathbf E|\eta\n_i(x,t)|^2$ and $\mathbf E|\eta\m_i(x,t)|^2$ satisfying the inequality
\begin{multline*}
\sum_{i=1}^{N}\! \int_{t_0}^T\int_{\gamma_i}\!\! \mathbf E|\eta\n_i(x,t)|^2 \left(r\n_i(x,t)\right)^2 d\gamma_i\,dt\\+
\sum_{i=1}^{N}\! \int_{t_0}^T\int_{\gamma_i}\!\! \mathbf E|\eta\m_i(x,t)|^2 \left(r\m_i(x,t)\right)^2 d\gamma_i\,dt\leq 1,
\end{multline*}
where $r_i \n(x,t), r_i \m(x,t)$ are given functions continuous on $\gamma_i\times (t_0,T),$ $i= \overline{1,N},$ that do not vanish on these sets.
\end{predlllll}

\newpage
\makeatletter
\renewcommand{\section}{\@startsection{section}{1}%
{\parindent}{3.5ex plus 1ex minus .2ex}%
{2.3ex plus.2ex}{\normalfont\Large{\bf PART\ \ }}} \makeatother

\begin{center}
\section[
Minimax estimation of the solutions to the boundary value  problems
from point observations
]{}
\end{center}

\begin{quote}
{\bf Minimax estimation of the solutions to the boundary value  problems
from point observations
}
\end{quote}

In the previous chapters we looked for estimates of unknown
solutions (and the right-hand sides of equations entering the
statements of the corresponding problems) from the observations of
these solutions distributed on a system of subdomains or surfaces.
In this chapter, we consider similar problems in the case of point
observations and propose constructive minimax estimation methods.

Let $x_k',$ $k=\overline{1,N}$ and $x_k,$ $k=\overline{1,m}$ be
given systems of points belonging to domain $\mathbb
R^3\setminus\bar\Omega.$ The problem is as follows: to estimate
the expression
\begin{equation}\label{2p3}
l(\varphi)=\sum_{i=1}^{m}\bar a_i\varphi(x_i),
\end{equation}
from the observations of the form
\begin{equation}\label{1p3}
y_k=\varphi(x_k')+\eta_k,\,\,k=\overline{1,N},
\end{equation}
that correspond to the system state $\varphi$ described by problem
(\ref{1})$-$(\ref{3'}) in the class of estimates
\begin{equation}\label{3p3}
\widehat{l(\varphi)}=\sum_{k=1}^N \bar u_i
y_k+c,
\end{equation}
linear with respect to observations (\ref{1p3}) under the
following assumptions:
$h\in G_0$ and $\eta:=(\eta_{1},\ldots,\eta_{N})\in G_1,$ where
the set $G_0$ is given by formula \eqref{7},
 $\eta_i$ are errors of observations \eqref{1p3} that are realizations of random
 quantities
 $\eta_i=\eta_i(\omega)$, $G_1$ is the set
of random vectors $ \tilde \eta=(\tilde \eta_{1},\ldots,\tilde
\eta_{N})$ with zero expectations and finite second moments
satisfying the condition
\begin{equation}\label{jdm}
\sum_{i=1}^{N}r^2_i\mathbf E|\tilde\eta_i|^2 \leq 1,
\end{equation}
and $u_i\in \mathbb C,$ $i= \overline{1,N},$, $a_i\in \mathbb C,$
$i= \overline{1,m},$ $q_i\in \mathbb R,$ $i= \overline{1,N},$  and
$r_i\neq 0$ are given numbers.

Set $ u:=(u_1,\ldots,u_N)\in \mathbb R^N. $

{\bf Definition.} {\it The estimate
$$
\widehat{\widehat
{l(\varphi)}}=\sum_{i=1}^N\overline{\hat u_i}
y_i+\hat c,
$$
in which numbers $\hat u_i$ and $\hat c$ are determined from the
condition
\begin{equation} \label{11f}
\sup_{\tilde h \in G_0,\, \tilde \eta
\in G_1}\mathbf E|l(\tilde\varphi)-\widehat
{l(\tilde\varphi)}|^2 \to \inf_{u\in \mathbb R^N,\,c\in \mathbb C} ,
\end{equation}
where
\begin{equation} \label{llxf}
\widehat
{l(\tilde\varphi)}=\sum_{i=1}^N\overline{u_i}
\tilde y_i +c,
\end{equation}
\begin{equation} \label{4iuf}
\tilde y_{i}=\tilde \varphi(x_i)+ \tilde \eta_{i},\quad i= \overline{1,N},
\end{equation}
and $\tilde\varphi(x)$ is the solution to the Neumann BVP at  $h=\tilde h,$
will be called the minimax estimate of expression (\ref{2p3}).

The quantity
\begin{equation} \label{12df}
\sigma:=\{\sup_{\tilde h \in G_0,\,
\tilde \eta \in G_1}\mathbf E|l(\tilde \varphi)-\widehat{\widehat
{l(\tilde\varphi)}}|^2\}^{1/2}
\end{equation}
will be called the error of the minimax estimation of $l(\varphi).$}

Based on the proof similar to  that of Lemma \ref{lem2.1} (in
fact, much simpler) we can show that the following statement is
valid in the case of point observations.
\begin{predl}
Finding the minimax estimate of functional
$l(\varphi)$ is equivalent to the problem of optimal control of the system described by BVP
$$
z(\cdot;u)\in \mathcal D'(\mathbb R^3\setminus\bar\Omega),
$$
$$
\Delta z(x;u)+\bar
k^2z(x;u)=\sum_{i=1}^{m}a_i\delta(x-x_i)-\sum_{k=1}^{N}u_k\delta
(x-x'_k)\,\, \mbox{in}\,\,\mathbb R^3\setminus\bar\Omega,
$$
$$
\frac{\partial z(\cdot;u)
}{\partial\nu}=0\,\,\mbox{on}\,\,\Gamma,
$$
$$
\frac{\partial z(\cdot;u)}{\partial
r}+i   \bar kz(\cdot;u)=o(1/r),\,\,r=|x|,\,\,r\to \infty.
$$
with the cost function
\begin{equation}\label{17f} I(u)
=\int_{\Gamma}q_1^{-2}(x)
z^2(x;u)\,d\Gamma+\sum_{i=1}^{N}r^{-2}_i|u_i|^2 \to \min_{u\in \mathbb R^N}.
\end{equation}
\end{predl}

 Starting from this lemma and proceeding with the reasoning that led from Lemma
  \ref{lem2.1} to Theorems \ref{t4} and \ref{t6}, we arrive at the
  following result
\begin{pred}\label{3.1}
The minimax estimate of $l(\varphi)$ has the form
\begin{equation}\label{g2rta}
\widehat{\widehat{l(\varphi)}}=\sum_{i=1}^N\overline{\hat u_i}
y_i+\hat c=l(\hat\varphi),
\end{equation}
where
\begin{equation}\label{g2rt1p3}
\hat u_k=r_k^2p(x_k'),
\quad k=\overline{1,N},\quad\hat
c=\int_\Gamma\overline{z}
h_0\,d\Gamma,
\end{equation}
the functions $z,p,\in \mathcal D'(\mathbb
R^3\setminus\bar\Omega)$ and $\hat
\varphi=\hat\varphi(\cdot,\omega)\in \mathcal D'(\mathbb
R^3\setminus\bar\Omega)$ are determined, respectively, from the
solution to the following problems:
\begin{equation}\label{1urr3}
-(\Delta +\bar
k^2)z(x)=\sum_{i=1}^{m}a_i\delta(x-x_i)-\sum_{k=1}^{N}r_k^2p(x'_k)\delta
(x-x'_k)\,\, \mbox{in}\,\,\Omega,
\end{equation}
\begin{equation}\label{2urr3}
\frac{\partial z
}{\partial\nu}=0\,\,\mbox{on}\,\,\Gamma,
\end{equation}
\begin{equation}\label{3urr3}
\frac{\partial z}{\partial
r}+i\bar kz=o(1/r),\,\,r=|x|,\,\,r\to \infty,
\end{equation}
\begin{equation}\label{4urr3}
\Delta p(x)+k^2p(x)=0\,\, \mbox{in}\,\,\Omega,
\end{equation}
\begin{equation}\label{5urr3}
\frac{\partial p}{\partial\nu}=q_1^{-2}z,\,\,\mbox{on}\,\,\Gamma,
\end{equation}
\begin{equation}\label{6urr3}
\frac{\partial p}{\partial
r}-ikp=o(1/r),\,\,r=|x|, \,\,r\to \infty.
\end{equation}
and
\begin{equation}\label{1urr3p}
-(\Delta +\bar
k^2)\hat p(x)=\sum_{i=1}^{N}r_i^2[y(x_k')-\varphi(x'_k)]\delta
(x-x'_k)\,\, \mbox{in}\,\,\Omega,
\end{equation}
\begin{equation}\label{2urr3p}
\frac{\partial \hat p
}{\partial\nu}=0\,\,\mbox{on}\,\,\Gamma,
\end{equation}
\begin{equation}\label{3urr3p}
\frac{\partial \hat p}{\partial
r}+i\bar k\hat p=o(1/r),\,\,r=|x|,\,\,r\to \infty,
\end{equation}
\begin{equation}\label{4urr3p}
\Delta \hat\varphi(x)+k^2\hat\varphi(x)=0\,\, \mbox{in}\,\,\Omega,
\end{equation}
\begin{equation}\label{5urr3p}
\frac{\partial \hat\varphi}{\partial\nu}=q_1^{-2}\hat p+
h_0\,\,\mbox{on}\,\,\Gamma,
\end{equation}
\begin{equation}\label{6urr3p}
\frac{\partial \hat\varphi}{\partial
r}-ik\hat\varphi=o(1/r),\,\,r=|x|, \,\,r\to \infty.
\end{equation}
Problem (\ref{1urr3})$-$(\ref{6urr3p}) is uniquely solvable.
The following estimate is valid for the error $\sigma$ of the minimax estimation of
$l(\varphi)$
\begin{equation} \label{34'}
\sigma=[l(p)]^{1/2}= \left(\sum_{i=1}^{m}\bar a_ip(x_i)\right)^{1/2}.
\end{equation}
\end{pred}

In conclusion, we formulate the statements similar to Theorems
\ref{t3.2'} and \ref{t3.2a'} that enable one to reduce, in line
with the algorithm applied in the proof of Theorem \ref{3.1}, the
determination of minimax estimates to a problem of less
dimensionality.
\begin{pred}
The minimax estimate of $l(\varphi)$ has the form
\begin{equation}\label{pp3}
\widehat{\widehat{l(\varphi)}}=\sum_{k=1}^N \bar {\hat u}_k
y_k+\hat c, \,\,
\end{equation}
where
\begin{gather} \label{p1p3}
\hat u_k=r_k^2p(x'_k), \,\,k=\overline{1,N},\\
\hat c=\int_{\Gamma}\overline{\left[\psi+
\left.\left(\sum_{j=1}^ma_j\Phi_{-\bar
k}(\cdot-x_j)
-\sum_{l=1}^{N}q_l^2p(x'_l)\Phi_{-\bar
k}(\cdot-x'_l)\right)\right|_{\Gamma}\right]}g_0\,d\Gamma,\label{p1p3'}
\end{gather}
and functions
$$
\psi:=\left.z\right|_{\Gamma}-\left.\left(\sum_{j=1}^ma_j\Phi_{-\bar
k}(\cdot-x_j)
-\sum_{l=1}^{N}q_l^2p(x'_l)\Phi_{-\bar
k}(\cdot-x'_l)\right)\right|_{\Gamma},
$$
$\chi:=p|_{\Gamma},$   and numbers $p(x_l'),$ $l=\overline{1,N},$
are determined from the solution of the following equation system:
$$
(I-K_{\!-\bar k}'+i\eta T_{\!-\bar k})\psi=(-S_{-\bar k}
+i\eta(I +K_{-\bar k}'))
$$
\begin{equation}\label{a1p3}
\left(-\sum_{j=1}^ma_j\left.\frac{\partial\Phi_{-\bar k}(\cdot-x_j)}{\partial \nu}\right|_\Gamma\right.
\left.+\sum_{l=1}^{N}r_l^2p(x'_l)
\left.\frac{\partial \Phi_{-\bar k}(\cdot-x_l')}{\partial \nu}\right|_\Gamma
\right),
\end{equation}
\begin{equation}\label{a2p3}
(I-K_k-i\eta T_k)\chi
\end{equation}
\begin{equation*}
=(-S_k -i\eta(I
+K_k'))q_1^{-2}\left(\psi+\sum_{j=1}^ma_j\left.\Phi_{-\bar
k}(\cdot-x_j)\right|_\Gamma
-\sum_{l=1}^{N}r_l^2p(x'_l)\left.\Phi_{-\bar
k}(\cdot-x'_l)\right|_\Gamma\right),
\end{equation*}
\begin{equation}\label{a3p3}
p(x_i')=\mathcal W_k\chi(x_r')
\end{equation}
\begin{equation*}
\!-\!\left.\mathcal
V_kq_1^{-2} \left(\psi\!+\!\sum_{j=1}^ma_j\left.\Phi_{-\bar
k}(\cdot-x_j)\right|_\Gamma
\!-\!\sum_{l=1}^{N}r_l^2p(x'_l)\left.\Phi_{-\bar
k}(\cdot-x'_l)\right|_\Gamma\right)\right|_{x_i'} ,\,\,
i=\overline{1,N}.
\end{equation*}
in which
$\eta$ is an arbitrary real number such that $\eta\mbox{\rm Re}\,k>0.$
This system is uniquely solvable for all values of wave numbers $k$,  $\mbox{\rm Im\,} k\geq 0.$
\end{pred}
Equation \eqref{a3p3} may be rewritten in the form
\begin{equation}\label{a3'p3}
[1-\alpha_{ii}]p(x_i')
+\sum_{l=1,\,l\neq
i}^{N}\alpha_{li}p(x'_l)\\ =\mathcal
W_k\chi(x_i')-\mathcal
V_kq_1^{-2}\psi(x_i')+\beta_i,\,\, i=\overline{1,N},
\end{equation}
where
$$
\alpha_{ls}=\frac 1{2\pi}r_l^2\int_\Gamma\frac {e^{ik|x_s'-y|}}{|x_s'-y|}q_1^{-2}(y)
\frac {e^{-i\bar k|y-x_l'|}}{|y-x_l'|}\,d\Gamma_y,\quad l,s=\overline{1,N},
$$
$$
\beta_s=\frac 1{2\pi}\sum_{j=1}^ma_j\int_\Gamma\frac {e^{ik|x_s'-y|}}{|x_s'-y|}q_1^{-2}(y)
\frac {e^{-i\bar k|y-x_j|}}{|y-x_j|}\,d\Gamma_y,\quad s=\overline{1,N}.
$$

\begin{pred}
The minimax estimate of $l(\varphi)$ has the form
\begin{equation}\label{pp33}
\widehat{\widehat{l(\varphi)}}=l(\hat\varphi)=\sum_{i=1}^{m}\bar a_i\hat\varphi(x_i),
\end{equation}
where
\begin{multline}\label{a3asss}
\hat\varphi=\mathcal W_{\Gamma}^k\tilde\chi-\mathcal
V_{\Gamma}^k\Bigl[q_1^{-2}\Bigl(\tilde\psi
+\sum_{l=1}^{N}r_l^2(y_l-\hat\varphi(x'_l))\left.\Phi_{-\bar
k}(\cdot-x'_l)\right|_\Gamma+h_0\Bigr]\,\,\mbox{in}\,\,\mathbb R^3\setminus\bar\Omega
\end{multline}
and functions
$$
\tilde\psi:=\left.\hat p\right|_{\Gamma}
-\sum_{l=1}^{N}r_l^2(y_l-\hat\varphi(x'_l))\left.\Phi_{-\bar
k}(\cdot-x'_l)\right|_{\Gamma},
$$
$\tilde\chi:=\hat\varphi|_{\Gamma},$   and numbers $\hat\varphi(x_l'),$ $l=\overline{1,N},$
are determined from the solution of the following equation system:
\begin{multline}\label{a1p3n}
(I-K_{\!-\bar k}'+i\eta T_{\!-\bar k})\tilde\psi\\=(-S_{-\bar k}
+i\eta(I +K_{-\bar k}'))
\left(\sum_{l=1}^{N}r_l^2(\hat\varphi(x'_l)-y_l)
\left.\frac{\partial \Phi_{-\bar k}(\cdot-x_l')}{\partial \nu}\right|_\Gamma\right),
\end{multline}
\begin{equation}\label{a2p3n}
(I-K_k-i\eta T_k)\tilde\chi
\end{equation}
\begin{equation*}
=(-S_k -i\eta(I
+K_k'))\left[q_1^{-2}\left(\tilde\psi+\sum_{l=1}^{N}r_l^2(y_l-\hat\varphi(x'_l))\left.\Phi_{-\bar
k}(\cdot-x'_l)\right|_\Gamma\right)+h_0\right],
\end{equation*}
\begin{equation}\label{a3p3n}
\hat\varphi(x_i')=\mathcal W_k\tilde\chi(x_r')
\end{equation}
\begin{equation*}
\!-\!\left.\mathcal
V_k\left[q_1^{-2} \left(\tilde\psi+\sum_{l=1}^{N}r_l^2(y_l-\hat\varphi(x'_l))\left.\Phi_{-\bar
k}(\cdot-x'_l)\right|_\Gamma\right)+h_0\right]\right|_{x_i'} ,\,\,
i=\overline{1,N}.
\end{equation*}
in which
$\eta$ is an arbitrary real number such that $\eta\mbox{\rm Re}\,k>0.$
This system is uniquely solvable for all values of wave numbers $k$,  $\mbox{\rm Im\,} k\geq 0.$
\end{pred}

\begin{predlllll} Similar results can be obtained if the estimated functional has the form \eqref{9}.
\end{predlllll}
\newpage

\makeatletter
\renewcommand{\section}{\@startsection{section}{1}%
{\parindent}{3.5ex plus 1ex minus .2ex}%
{2.3ex plus.2ex}{\normalfont\Large{\bf\ \ }}} \makeatother

\addcontentsline{toc}{section}{References}
\renewcommand{\refname}
{\Large \bf \begin{center} References
\end{center}}

\end{document}